\documentclass[aos,preprint]{imsart} 

\usepackage{amsthm,amsmath,natbib}
\usepackage{amsfonts,amsbsy,amssymb}
\usepackage[mathscr]{euscript}
\usepackage{bbm}
\usepackage{mathtools}
\usepackage{graphicx}
\usepackage{xcolor}
\usepackage{multirow}
\usepackage{lipsum}
\usepackage{comment}
\usepackage{array}
\usepackage[nomessages]{fp} 
\usepackage{multibib} 
\newcites{appendix}{References}

\usepackage{booktabs}

\RequirePackage[colorlinks,citecolor=blue,urlcolor=blue]{hyperref}
\RequirePackage{natbib}

\startlocaldefs
\newtheorem{theorem}{Theorem}[section]
\newtheorem{algorithm}{Algorithm}[section]
\newtheorem{proposition}{Proposition}[section]
\newtheorem{corollary}{Corollary}[section]
\newtheorem{lemma}{Lemma}[section]
\newtheorem{definition}{Definition}[section]
\DeclareMathOperator{\E}{\mathbb{E}}
\DeclareMathOperator{\Prob}{\mathbb{P}}

\DeclareMathOperator{\Cov}{\operatorname{Cov}}
\DeclareMathOperator{\V}{\mathbb{V}}

\DeclareMathOperator{\supp}{supp}

\endlocaldefs
\begin{document}
\begin{frontmatter}
  \title{On the Optimal Reconstruction of Partially Observed Functional Data}
  \runtitle{Reconstructing Partially Observed Functions}
\author{
  \fnms{Alois}
  \snm{Kneip}\ead[label=e1]{akneip@uni-bonn.de}
}
\and
\author{
  \fnms{Dominik}
  \snm{Liebl}
  \corref{}\ead[label=e2]{dliebl@uni-bonn.de}
}
\address{
  Alois Kneip and Dominik Liebl\\
  Statistische Abteilung\\
  University of Bonn\\
  Adenauerallee 24-26\\
  53113 Bonn, Germany\\
  \printead{e1} \\
  \printead{e2}
}

\affiliation{University of Bonn}
\runauthor{A.~Kneip \& D.~Liebl}

\begin{abstract}
We propose a new reconstruction operator that aims to recover the missing parts of a function given the observed parts. This new operator belongs to a new, very large class of functional operators which includes the classical regression operators as a special case. We show the optimality of our reconstruction operator and demonstrate that the usually considered regression operators generally cannot be optimal reconstruction operators.
Our estimation theory allows for autocorrelated functional data and considers the practically relevant situation in which each of the $n$ functions is observed at $m_i$, $i=1,\dots,n$, discretization points. We derive rates of consistency for our nonparametric estimation procedures using a double asymptotic.
For data situations, as in our real data application where $m_i$ is considerably smaller than $n$, we show that our functional principal components based estimator can provide better rates of convergence than conventional nonparametric smoothing methods.
\end{abstract}


\end{frontmatter}

\section{Introduction}\label{sec:intro}
\begin{figure}[!ht]
\centering
\includegraphics[width=.5\textwidth]{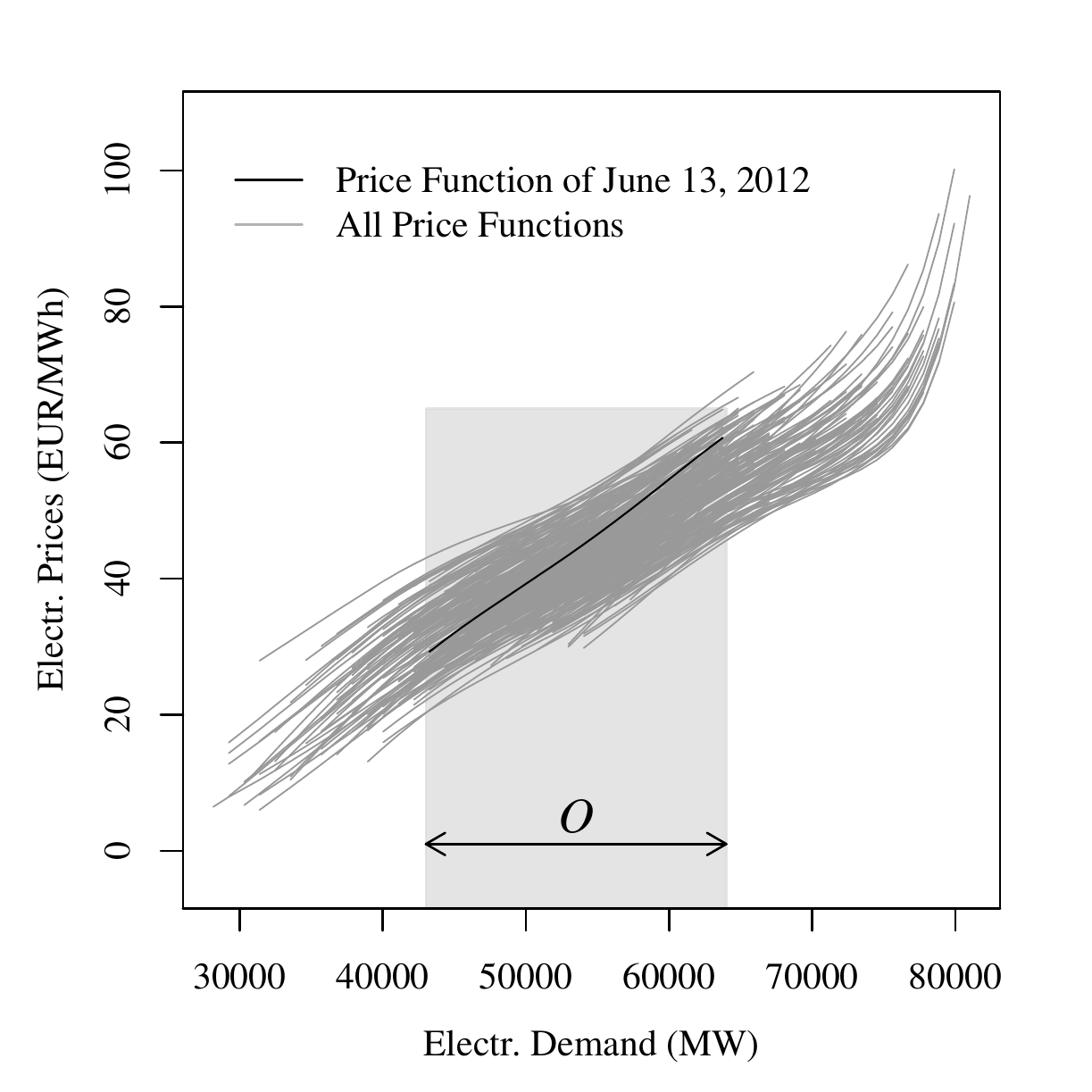}
\caption[]{Partially observed electricity price functions $X_i^O(u)$ with $u\in O_i\subseteq[a,b]$.}
\label{fig:presm_scatter}
\end{figure}

Our work is motivated by a data set from energy economics which is shown in Figure \ref{fig:presm_scatter}. The data consist of partially observed price functions. Practitioners use these functions, for instance, to do comparative statics, i.e., a ceteris-paribus analysis of price effects with respect to changes in electricity demand \citep[cf.][]{weigt2009germany,hirth2013market}. The possibilities of such an analysis, however, are limited by the extent to which we can observe the price functions. This motivates the goal of our work, which is to develop a reconstruction procedure that allows us to recover the total functions from their partial observations.

Let $X_1,\dots,X_n$ be an identically distributed, possibly weakly dependent sample of continuous random functions, where each function $X_i$ is an element of the separable Hilbert space $\mathbb{L}^2([a,b])$ with $[a,b]\subset\mathbb{R}$ and $\E(||X_i||_2^4)<\infty$, where $||X_i||_2^2=\int_a^b(X_i(x))^2dx$. We denote the observed and missing parts of $X_i$ by $X_i^{O_i}$ and $X_i^{M_i}$, where
\begin{equation*}
\begin{array}{ll}
X_i^{O_i}(u)\,:=X_i(u)&\text{ for}\quad u\in O_i\subseteq[a,b]\quad\text{and}\\
X_i^{M_i}(u)  :=X_i(u)&\text{ for}\quad u\in M_i=[a,b]\setminus O_i,
\end{array}
\end{equation*}
and where $O_i=[A_i,B_i]\subseteq[a,b]$ is a random subinterval, independent from $X_i$ with $B_i-A_i>0$ almost surely. In our theoretical part (Section \ref{sec:opt_pred}) we also allow for the general case, where $O_i$ consists of multiple subintervals of $[a,b]$. In what follows we use ``$O$'' and ``$M$'' to denote a given realization of $O_i$ and $M_i$. In addition, we use the following shorthand notation for conditioning on $O_i$ and $M_i$:
\begin{equation*}
\begin{array}{l}
X_i^{O}(u)\,:=X_i^{O_i}(u)\;|(O_i=O)\\
X_i^{M}(u):=X_i^{M_i}(u)|(M_i=M);\\
\end{array}
\end{equation*}
typical realizations of $X_i^{O}$ and $O$ are shown in Figure \ref{fig:presm_scatter}. We denote the inner product and norm of $\mathbb{L}^2(O)$ as $\langle .,.\rangle_2$ and $||.||_2$; their dependency on $O$ will be made obvious by writing, for instance, $\langle x^O,y^O\rangle_2$ and $||x^O||_2^2=\langle x^O,x^O\rangle_2$ for all $x^O,y^O\in \mathbb{L}^2(O)$, where $\langle x^O,y^O\rangle_2=\int_{O} x(u)y(u)du$. Throughout the introduction and Section \ref{sec:opt_pred}, we consider centered random functions, that is, $\E(X_i(u))=\mu(u)$ with $\mu(u)=0$ for all $u\in[a,b]$. 

Our object of interest is the following linear reconstruction problem:
\begin{align}\label{eq:FPM:1}
  X_i^M=L(X_i^O)+Z_i,
\end{align}
which aims to reconstruct the unobserved missing parts $X_i^M\in\mathbb{L}^2(M)$ given the partial observation $X_i^O\in\mathbb{L}^2(O)$. Our objective is to identify the optimal linear reconstruction operator $L:\mathbb{L}^2(O)\to\mathbb{L}^2(M)$ which minimizes the mean squared error loss $\E\big[\big(X_i^M(u)-L(X_i^O)(u)\big)^2\big]$ at any $u\in M$.

The case of partially observed functional data was initially considered in the applied work of \cite{Liebl2013} and the theoretical works of \cite{Goldberg_JSPI_2014} and \cite{kraus2015}. The work of \cite{gromenko2017evaluation} is also related as it proposes an inferential framework for incomplete spatially and temporally correlated functional data. \cite{Goldberg_JSPI_2014} consider the case of finite dimensional functional data and their results have well-known counterparts in multivariate statistics. \cite{kraus2015} starts by deriving his ``optimal'' reconstruction operator as a solution to the Fr\'echet-type normal equation, where he assumes the existence of a bounded solution. The theoretical results in our paper imply that this assumption generally holds only under the very restrictive case of linear regression operators, i.e., Hilbert-Schmidt operators. For showing consistency of his empirical reconstruction operator, \cite{kraus2015} restricts his work to this case of Hilbert-Schmidt operators. We demonstrate, however, that a Hilbert-Schmidt operator generally cannot be the optimal reconstruction operator.

In order to see the latter, we need some conceptional work. Hilbert-Schmidt operators on $\mathbb{L}^2$ spaces correspond to linear regression operators,
\begin{equation}\label{eq:RF_0}
L(X_i^O)(u)=\int_{O}\beta(u,v)X_i^O(v)dv,\quad\text{with}\quad\beta\in\mathbb{L}^2(M\times O).
\end{equation}
However, such a regression operator generally does \emph{not} provide the optimal solution of the reconstruction problem in \eqref{eq:FPM:1}. For instance, let us consider the ``last observed''(=``first missing'') points, namely, the boundary points\footnote{The boundary $\partial M$ of a subset $M$ is defined as $\partial M:=\overline{M}\cap\overline{O}$, where $\overline{M}$ and $\overline{O}$ denote the closures of the subsets $M$ and $O$.} $\vartheta\in\partial M$ of $M$. For any optimal reconstruction operator $L$, it must hold that the ``first reconstructed'' value, $L(X_i^O)(\vartheta)$, connects with the ``last observed'' value, $X_i^O(\vartheta)$, i.e., that
\begin{equation*}
X_i^O(\vartheta)=L(X_i^O)(\vartheta)\quad\text{for all}\quad\vartheta\in\partial M.
\end{equation*}
There is no hope, though, of finding a slope function $\beta(\vartheta,.)\in\mathbb{L}^2(O)$ that fulfills the equation $X_i^O(\vartheta)=\int_{O}\beta(\vartheta,v)X_i^O(v)dv$ (the Dirac-$\delta$ function is not an element of $\mathbb{L}^2(O)$). It is therefore impossible to identify the optimal reconstruction operator $L$ within the class of linear regression operators defined by \eqref{eq:RF_0}.

Best possible linear reconstruction operators depend, of course, on the structure of the random function $X_i$, and possible candidates have only to be well-defined for any function in the support of $X_i$. We therefore consider the class of all linear operators $L$ with finite variance $\V(L(X_i^O)(u))<\infty$ and thus $\Prob(|L(X_i^O)(u)|<\infty)=1$ for any $u\in M$. This class of reconstruction operators is much larger than the class of regression operators and contains the latter as a special case. A theoretical characterization is given in Section \ref{sec:opt_pred}. We then show that the optimal linear reconstruction operator, minimizing $\E[(X_i^M(u)-L(X_i^O)(u))^2]$ for all $u\in M$, is given by
\begin{align}\label{eq:PF_0}
\mathcal{L}(X_i^O)(u)=
\sum_{k=1}^\infty\frac{\xi^{O}_{ik}\,\E[X_i^M(u)\xi^{O}_{ik}]}{\lambda_k^{O}}
=\sum_{k=1}^\infty\xi^{O}_{ik}\ \frac{\langle\phi_k^O,\gamma_u\rangle_2}{\lambda_k^{O}},
\end{align}
where $(\phi_k^O,\lambda_k^O)_{k\geq 1}$ denote the pairs of orthonormal eigenfunctions and nonzero eigenvalues of the covariance operator $\Gamma^O(x)(u)=\int\gamma^O(u,v)x(v)dv$ with $x\in \mathbb{L}^2(O)$, while
$\xi^{O}_{ik}:=\langle\phi_k^O,X_i^O\rangle_2$. Here  $\gamma^O(u,v)=\Cov(X_i^O(u),X_i^O(v))$ denotes the covariance function of $X_i^O$, and $\gamma_u(v)=\gamma(u,v)$ the covariance function $\gamma(u,v)=\Cov(X_i^M(u),X_i^O(v))$.

The general structure of $\mathcal{L}$ in \eqref{eq:PF_0} is similar to the structure of the operators considered in the literature on functional linear regression, which, however, additionally postulates that $\mathcal{L}$ has an (restrictive) integral-representation as in \eqref{eq:RF_0}; see, for instance, \cite{cardot2007clt}, \cite{cai2006prediction}, \cite{hall2007methodology} in the context of functional linear regression, or \cite{kraus2015} in a setup similar to ours.

There is, however, no reason to expect that the optimal reconstruction operator $\mathcal{L}$ satisfies \eqref{eq:RF_0}. To see the point note that $\mathcal{L}(X_i^O)(u)$  can be represented in the form \eqref{eq:RF_0} \emph{if and only if} the \emph{additional} square summability condition $\sum_{k=1}^\infty \langle\phi_k^O,\gamma_u\rangle_2^2/(\lambda_k^{O})^2<\infty$ is satisfied for $u\in M$. Only then the series $\sum_{k=1}^L(\langle\phi_k^O,\gamma_u\rangle_2/\lambda_k^{O})\phi_k^O(v)$, $v\in O$, converge as $L\rightarrow\infty$ and  define a function $\beta_u:=\sum_{k=1}^\infty(\langle\phi_k^O,\gamma_u\rangle_2/\lambda_k^{O})\phi_k^O(\cdot)\in \mathbb{L}^2(O)$ such that  $\int_{O}\beta_u(v)X_i^O(v)dv=\sum_{k=1}^\infty \xi^{O}_{ik} \langle\phi_k^O,\beta_u\rangle_2=\sum_{k=1}^\infty\xi^{O}_{ik}\langle\phi_k^O,\gamma_u\rangle_2/\lambda_k^{O}$.

But consider again the reconstruction at a boundary point $\vartheta\in\partial M$, where $\langle\phi_k^O,\gamma_\vartheta\rangle_2$ simplifies to $\langle\phi_k^O,\gamma_\vartheta\rangle_2=\lambda_k^O\phi_k^O(\vartheta)$, since for boundary points $\vartheta$ we have $\gamma_\vartheta=\gamma^O_\vartheta$ and $\gamma^O_\vartheta(\cdot)=\gamma^O(\vartheta,\cdot)=\sum_{k=1}^\infty\lambda_k^O\phi_k^O(\vartheta)\phi_k^O(\cdot)$. Plugging this simplification into \eqref{eq:PF_0} and using the Karhunen-Lo\'eve decomposition of $X_i^O$ implies that $\mathcal{L}(X_i^O)(\vartheta)=\sum_{k=1}^\infty \xi^{O}_{ik} \phi_k^O(\vartheta)=X_i^O(\vartheta)$. This means that our reconstruction operator $\mathcal{L}$ indeed connects the ``last observed'' value $X_i^O(\vartheta)$ with the ``first reconstructed'' value $\mathcal{L}(X_i^O)(\vartheta)$. On the other hand, the sum $\sum_{k=1}^L\langle\phi_k^O,\gamma_\vartheta\rangle_2^2/(\lambda_k^{O})^2=\sum_{k=1}^L(\phi_k^O(\vartheta))^2$ will generally tend to infinity as $L\rightarrow\infty$,  which violates the additional condition necessary for establishing \eqref{eq:RF_0}. Therefore, in general, $\mathcal{L}$ does not constitute a regression operator.\footnote{A frequently used justification of the use of regression operators relies on the Riesz representation theorem which states that any continuous linear functional $L(X_i^O)(u)$ can be represented in the form \eqref{eq:RF_0}. This argument, however, does not necessarily apply to the optimal linear functional $\mathcal{L}(X_i^O)(u)$ which may not be a continuous functional $\mathbb{L}^2(O)\rightarrow \mathbb{R}$. In particular, although being a well-defined linear functional, the point evaluation $\mathcal{L}(X_i^O)(\vartheta)=X_i^O(\vartheta)$ is {\bf not continuous}, since for two functions $f,g\in \mathbb{L}^2(O)$ an arbitrarily small $\mathbb{L}^2$-distance $\Vert f-g\Vert_2$ may go along with a very large pointwise distance $|f(\vartheta)-g(\vartheta)|$ (see the example in Appendix \ref{appendix:FE1} of the supplementary paper \cite{KL_Suppl_19}).}

The above arguments indicate that methods for estimating $\mathcal{L}$ should not be based on \eqref{eq:RF_0}. Any theoretical justification of such procedures has to rely on non-standard asymptotics avoiding the restrictive assumption that $\sum_{k=1}^\infty\langle\phi_k^O,\gamma_u\rangle_2^2/(\lambda_k^{O})^2<\infty$. This  constitutes a major aspect of our asymptotic theory given in Section \ref{sec:asymp}.

The problem of estimating $\mathcal{L}(X_i^O)$ from real data is considered in Section \ref{sec:estim}. Motivated by our application, the estimation theory allows for an autocorrelated time series of functional data and considers the practically relevant case where the function parts $X_i^O$ are only observed at $m_i$ many discretization points $(Y_{i1},U_{i1}),\dots,(Y_{im_{i}},U_{im_{i}})$ with $Y_{ij}=X_i^O(U_{ij})+\varepsilon_{ij}$, $i=1,\dots,n$, and $j=1,\dots,m_{i}$.

We basically follow the standard approach to estimate $\mathcal{L}(X_i^O)$ through approximating the infinite series \eqref{eq:PF_0} by a truncated sequence relying only on the $K$ largest eigenvalues of the covariance operator. But note that our data structure implies that we are faced with two simultaneous estimation problems. One is efficient estimation of $\mathcal{L}(X_i^O)(u)$ for $u\in M$, the other one is a best possible estimation of the function $X_i(u)$ for $u\in O$ from the observations $(Y_{i1},U_{i1}),\dots,(Y_{im_{i}},U_{im_{i}})$. We consider two different estimation strategies; both allow us to accomplish these two estimation problems.

The first consists in using a classical functional principal components based approximation of $X_i$ on $O$, which is simply given by extending the operator $\mathcal{L}$ in \eqref{eq:PF_0} by extending $\gamma_u(v)=\Cov(X_i^M(u),X_i^O(u))$ to $\gamma_u(v)=\Cov(X_i(u),X_i(u))$. This way the empirical counterpart of the truncated sum
\begin{align*}
\mathcal{L}_{K}(X_i^O)(u)&=
\sum_{k=1}^K \xi^{O}_{ik}\ \frac{\langle\phi_k^O,\gamma_u\rangle_2}{\lambda_k^{O}},\quad\text{for}\quad u\in O\cup M,
\end{align*}
will simultaneously provide estimates of the true function $X_i^O(u)$ on the observed interval $O$ and of the optimal reconstruction $\mathcal{L}(X_i^O)(u)$ on the unobserved interval $M$.


The second consists in estimating the true function $X_i^O(u)$ on the observed interval $O$ directly from the observations $(Y_{i1},U_{i1}),\dots,(Y_{im_{i}},U_{im_{i}})$ using, for instance, a local linear smoother and to estimate $\mathcal{L}(X_i^O)(u)$ for $u\in M$ through approximating the infinite series \eqref{eq:PF_0} by its truncated version. But a simple truncation would result in a jump at a boundary point $\vartheta_u$, with $\vartheta_u$ denoting the closest boundary point to the considered $u\in M$, i.e., $\vartheta_u=A_i$ if $|A_i-u|<|B_i-u|$ and $\vartheta_u=B_i$ otherwise. We know, however, that for any $u\approx\vartheta_u$ we must have $\langle\phi_k^O,\gamma_u\rangle_2\approx\lambda_k^O \phi_k^O(\vartheta)$ for all $k\geq 1$, since $\langle\phi_k^O,\gamma_\vartheta\rangle_2=\lambda_k^O\phi_k^O(\vartheta)$ for all boundary points $\vartheta_u\in\partial M$. Therefore, we explicitly incorporate boundary points and estimate $\mathcal{L}(X_i^O)$ by the empirical counterpart of the truncated sum
\begin{equation*}
\mathcal{L}_{K}^*(X_i^O)(u)=X_i^O(\vartheta_u)+
\sum_{k=1}^K\xi^{O}_{ik}\, \left(\frac{\langle\phi_k^O,\gamma_u\rangle_2}{\lambda_k^O}-\phi_k^O(\vartheta_u)\right),\quad u\in M.
\end{equation*}
The above truncation does not lead to an artificial jump at a boundary point $\vartheta_u$, since $(\langle\phi_k^O,\gamma_u\rangle_2/\lambda_k^{O}-\phi_k^O(\vartheta_u))\to 0$ continuously as $u\to\vartheta_u$ for all $k=1,\dots,K$.

For estimating the mean and covariance functions -- the basic ingredients of our reconstruction operator -- we suggest using Local Linear Kernel (LLK) estimators. These LLK estimators are commonly used in the context of sparse functional data \citep[see, e.g.,][]{Yao2005}, though, we do \emph{not} consider the case of sparse functional data. In the context of partially observed functional data, it is advisable to use LLK estimators, since these will guarantee smooth estimation results, which is not the case when using the empirical moment estimators for partially observed functions as proposed in \cite{kraus2015}.

We derive consistency as well as uniform rates of convergence under a double asymptotic
which allows us to investigate all data scenarios from almost sparse to dense functional data. This leads to different convergence rates depending on the relative order of $m$ and $n$. For data situations, as in our real data application where $m$ is considerably smaller than $n$ and all sample curves are of similar structure, we show that our functional principal components based estimator achieves almost parametric convergence rates and can provide better rates of convergence than  conventional nonparametric smoothing methods, such as, for example, local linear regression.

Our development focuses on the regular situation where (with probability tending to 1) there exist functions that are observed over the total interval $[a,b]$. Only then is it possible to consistently estimate the covariance function $\gamma(u,v)$ for all possible pairs $(u,v)\in[a,b]^2$. In our application this is not completely fulfilled, and there is no information on $\gamma(u,v)$ for very large values $|u-v|$. Consequently, for some intervals $O$ and $M$ the optimal reconstruction operator cannot be identified. This situation corresponds to the case of so-called fragmentary observations, as considered by \cite{delaigle2013classification}, \cite{Delaigle11112016}, \cite{DP2017}, and \cite{DHHK2018}. To solve this problem we suggest an iterative reconstruction algorithm. Optimal reconstruction operators are determined for a number of smaller subintervals, and a final operator for a larger interval is obtained by successively plugging in the reconstructions computed for the subintervals. We also provide some inequality bounding the accumulating reconstruction error.

The rest of this paper is structured as follows: Section \ref{sec:opt_pred} introduces our reconstruction operator and contains the optimality result. Section \ref{sec:estim} comprises our estimation procedure. The asymptotic results are presented in Section \ref{sec:asymp}. Section \ref{sec:pred_algo} describes the iterative reconstruction algorithm. Section \ref{sec:sim} contains the simulation study and Section \ref{sec:appl} the real data application. All proofs can be found in the online supplement supporting this article \citep{KL_Suppl_19}.

\section{Optimal reconstruction of partially observed functions}\label{sec:opt_pred}
Let our basic setup be as described in Section \ref{sec:intro}. Any (centered) random function $X_i^O\in\mathbb{L}^2(O)$ then adopts the well-known Karhunen-Lo\'eve (KL) representation
\begin{align}\label{eq:KLR}
  X_i^{O}(u)=\sum_{k=1}^\infty\xi^{O}_{ik}\phi^{O}_{k}(u), \quad u\in O,
\end{align}
with the principal component (pc) scores $\xi^{O}_{ik}=\langle X_i^{O},\phi^{O}_{k}\rangle_2$, where $\E(\xi^{O}_{ik})=0$ and $\E(\xi^{O}_{ik}\,\xi^{O}_{il})=\lambda^{O}_k$ for all $k=l$ and zero else and $\lambda^{O}_1> \lambda^{O}_2>\dots >0$. We want to note that all arguments in this section also apply to the more general case where the observed subdomain $O=\bigcup_{j=1}^J[A_j,B_j]$ consists of a finite number $1\leq J<\infty$ of mutually disjoint subintervals $[A_j,B_j]\subseteq[a,b]$.

By the classical eigen-equations we have that
\begin{align}\label{eq:obs.ef}
\phi^{O}_{k}(u)=
\frac{\langle\phi_k^{O},\gamma_u^O\rangle_2}{\lambda_k^O}, \quad u\in O,
\end{align}
where $\gamma_u^O(v)=\gamma^O(u,v)=\E(X_i^O(u)X_i^O(v))$. Equation \eqref{eq:obs.ef} can obviously be generalized for all $u\in M$ which leads to the following ``extrapolated'' $k$th basis function:
\begin{align}\label{eq:PEF}
\tilde{\phi}_{k}^O(u)=
\frac{\langle\phi_k^{O},\gamma_u\rangle_2}{\lambda_k^O},\quad u\in M,
\end{align}
where $\gamma_u(v)=\E(X_i^M(u)X_i^O(v))$. Equation \eqref{eq:PEF} leads to the definition of our reconstruction operator $\mathcal{L}_u$ as a generalized version of the KL representation in \eqref{eq:KLR}:
\begin{align}\label{eq:PEFexp}
\mathcal{L}(X_i^O)(u)=\sum_{k=1}^\infty\xi_{ik}^O\;\tilde{\phi}_{k}^O(u), \quad u\in M.
\end{align}

\paragraph{\normalfont Remark}
\textit{
Note that the KL representation provides the very basis of a majority of the works in functional data analysis \citep[cf.][]{RamsayfdaBook2005,horvath2012inference}. Functional Principal Component Analysis (FPCA) relies on  approximating $X_i$ by its first $K$ principal components. This is justified by the {\bf best basis property}, i.e., the property that for any $K\geq 1$
\begin{align}\label{eq:best basis}
\sum_{k=K+1}^{\infty} \lambda_k^O &=\E\left(
 \Vert X_i^{O}(u)-\sum_{k=1}^K\xi^{O}_{ik}\phi^{O}_{k}(u)\Vert_2^2\right)\nonumber \\
& =\min_{v_1,\dots,v_K\in \mathbb{L}^2(O)} \E\left(\min_{a_{i1},\dots,a_{ik}\in\mathbb{R}}
 \Vert X_i^{O}(u)-\sum_{k=1}^K a_{ik}v_{k}(u)\Vert_2^2\right).
\end{align}
}

\paragraph{\normalfont Remark}
\textit{
For later use it is important to note that the definitions of $\tilde{\phi}_k^O(u)$ and $\mathcal{L}(X_i^O)(u)$ in \eqref{eq:PEF} and \eqref{eq:PEFexp} can be extended for all $u\in O\cup M$ by setting $\gamma_u=\E(X_i(u)X_i(v))$. Then by construction $\tilde{\phi}_k^O(u)=\phi_k^O(u)$ for all $u\in O$ and, therefore, $\mathcal{L}(X_i^O)(u)=X_i^O(u)$ for all $u\in O$.
}

\subsection{A theoretical framework for reconstruction operators}
Before we consider the optimality properties of $\mathcal{L}$, we need to define a sensible class of operators against which to compare our reconstruction operator. We cannot simply choose the usual class of regression operators, since  $\mathcal{L}$ does generally not belong to this class, as pointed out in Section \ref{sec:intro}. Therefore, we introduce the following (very large) class of ``reconstruction operators'':
\begin{definition}[Reconstruction operators]\label{def:SRO}
Let the (centered) random function $X_i^O$ have a KL representation as in \eqref{eq:KLR}. We call every linear operator $L:\mathbb{L}^2(O)\to\mathbb{L}^2(M)$ a ``reconstruction operator with respect to $X_i^O$'' if $\;\V(L(X_i^O)(u))<\infty$ for all $u\in M$.
\end{definition}

It is important to note that this definition of ``reconstruction operators'' is specific to the considered process $X_i$. This should not be surprising, since a best possible linear reconstruction will of course depend on the structure of the relevant random function $X_i$. The following theorem provides a useful representation of this class of linear operators:
\begin{theorem}[Representation of reconstruction operators]\label{th:SRO}
Let $L:\mathbb{L}^2(O)\to\mathbb{L}^2(M)$ be a ``reconstruction operator with respect to $X_i^O$'' according to Definition \ref{def:SRO}. Then there exists a unique (deterministic) parameter function $\alpha_u\in H$ such that almost surely
\begin{equation*}
L(X_i^O)(u)=\langle\alpha_u,X_i^O\rangle_H,\quad u\in M,
\end{equation*}
where $H:=\left\{f\in \mathbb{L}^2(O):\;||f||^2_H<\infty\right\}$ is a Hilbert space with inner product $\langle f,g\rangle_H:=\sum_{k=1}^\infty\langle f,\phi_k^O\rangle_2\;\langle g,\phi_k^O\rangle_2/\lambda_k^O$  for all $f,g\in \mathbb{L}^2(O)$ and induced norm $||f||_H=\sqrt{\langle f,f\rangle_H}$.
\end{theorem}
The space $H$ is the Reproducing Kernel Hilbert Space (RKHS) that takes the covariance kernel $\gamma^O(u,v)=\sum_{k=1}^\infty \lambda_k^O\phi_k^O(u)\phi_k^O(v)$ as its reproducing kernel. By construction, we obtain that the variance of $L(X_i^O)(u)$ equals the $H$-norm of the parameter function $\alpha_u$, i.e., $\V(L(X_i^O)(u))=||\alpha_u||^2_H$.

\smallskip
Let us consider two examples of possible reconstruction operators. While the first example does not belong the class of regression operators, the second example is a regression operator demonstrating the more restrictive model assumptions.
\smallskip

\noindent{\bf Example 1 - Point of impact:}
Consider $L(X_i^O)(u) = X_i^O(\tau)$, i.e., a model with one ``impact point'' $\tau\in O$ for all missing points $u\in M$.  With $\gamma_\tau(v):=\gamma(\tau,v)=\sum_{k=1}^\infty \lambda_k^O \phi_k^O(\tau) \phi_k^O(v)$ we have $\lambda_k^O \phi_k^O(\tau)=\langle \gamma_\tau,\phi_k^O \rangle_2$, and hence
\begin{align}
&L(X_i^O)(u)= X_i^O(\tau) = \sum_{k=1}^\infty\xi^{O}_{ik}\phi^{O}_{k}(\tau)= \sum_{k=1}^\infty\frac{\langle X_i^O,\phi_k^O \rangle_2 \lambda_k^O \phi_k^O(\tau)}{\lambda_k^O}=\nonumber\\
&=\sum_{k=1}^\infty\frac{ \langle X_i^O,\phi_k^O \rangle_2 \langle \gamma_\tau,\phi_k^O \rangle_2}{\lambda_k^O}=\langle\gamma_\tau,X_i^O\rangle_H,
\label{rcoex1}
\end{align}
where $\gamma_\tau(\cdot):=\sum_{k=1}^\infty \lambda_k^O \phi_k^O(\tau) \phi_k^O(\cdot)\in H$ with
$||\gamma_\tau||^2_H=\sum_{k=1}^\infty \frac{(\lambda_k^O)^2\phi_k^O(\tau)^2}{\lambda_k^O} =\sum_{k=1}^\infty \lambda_k^O \phi_k^O(\tau)^2=\V(X_i(\tau))<\infty$.

\smallskip

\noindent{\bf Example 2 - Regression operator:}
Let $L$ be a regression operator (see \eqref{eq:RF_0}). Then there exists a $\beta_u\in\mathbb{L}^2(O)$ such that $L(X_i^O)(u)=\langle\beta_u,X_i^O\rangle_2$. Since eigenfunctions can be completed to an orthonormal basis of $\mathbb{L}^2(O)$, we necessarily have that $\sum_{k=1}^\infty \beta_{u,k}^2<\infty$ for $\beta_{u,k}:=\langle\beta_u, \phi_k^O\rangle_2$. Then
\begin{align}
&L(X_i^O)(u)=\langle\beta_u,X_i^O\rangle_2=\sum_{k=1}^\infty \xi^{O}_{ik} \beta_{u,k}
=\sum_{k=1}^\infty\frac{ \langle X_i^O,\phi_k^O \rangle \lambda_k^O \beta_{u,k} }{\lambda_k^O}\nonumber \\
&=\sum_{k=1}^\infty\frac{ \langle X_i^O,\phi_k^O \rangle \langle\alpha_u,\phi_k^O \rangle}{\lambda_k^O}
=\langle\alpha_u,X_i^O\rangle_H,  \label{rcoex2}
\end{align}
where
$\alpha_u(\cdot):=\sum_{j=1}^\infty \lambda_k^O \beta_{u,k} \phi_k^O(\cdot)\in H$
with
$||\alpha||^2_H=\sum_{k=1}^\infty \frac{(\lambda_k^O)^2\beta_{u,k}^2}{\lambda_k^O}
=\sum_{k=1}^\infty \lambda_k^O \beta_{u,k}^2<\infty$. Also note that for any $k$ we have $\langle\alpha_u,\phi_k^O\rangle_2=\lambda_k^O \beta_{u,k}$. This means that for $\alpha_u\in H$ the operator
$\langle\alpha_u,X_i^O\rangle_H$ constitutes a regression operator if and only if in addition to $\Vert\alpha_u\Vert_H^2
=\sum_{k=1}^\infty \langle\alpha_u,\phi_k^O\rangle_2^2/\lambda_k^O<\infty$ we also have that
$\sum_{k=1}^\infty \langle\alpha_u,\phi_k^O\rangle_2^2/(\lambda_k^O)^2<\infty$ (the latter is not satisfied in Example 1).

\smallskip

These examples show that Definition \ref{def:SRO} leads to a very large class of linear operators which contains the usually considered class of regression operators as a special case. Of course, the class of reconstruction operators as defined by Definition \ref{def:SRO} also contains much more complex operators than those illustrated in the examples.

Using Theorem \ref{th:SRO}, our reconstruction problem in \eqref{eq:PF_0} of finding a ``best linear'' reconstruction operator minimizing the squared error loss can now be restated in a theoretically precise manner: Find the linear operator $L:\mathbb{L}^2(O)\to\mathbb{L}^2(M)$ which for all $u\in M$ minimizes
\begin{equation*}
\E\Big[\big(X_i^M(u)-L(X_i^O)(u)\big)^2\Big]
\end{equation*}
with respect to all reconstruction operators $L$ satisfying $L(X_i^O)(u)=\langle\alpha_u,X_i^O\rangle_H$ for some $\alpha_u\in H$. In the next subsection we show that the solution is given by the operator $\mathcal{L}$ defined in \eqref{eq:PEFexp} which can now be rewritten in the form
\begin{align}\label{eq:PF_S}
\mathcal{L}(X_i^O)(u)=\langle\gamma_u,X_i^O\rangle_H, \quad u\in M,
\end{align}
where $\gamma_u(v)=\gamma(u,v)$ for $v\in O$ and $u\in M$. In particular, Theorem \ref{WellDef} below shows that $\V(\mathcal{L}(X_i^O)(u))=||\gamma_u||^2_H<\infty$ for any $u\in M$, i.e., that $\mathcal{L}$ is indeed a reconstruction operator according to Definition \ref{def:SRO}.

\paragraph{\normalfont Remark}
\textit{
In the context of reconstructing functions, problems with the use of regression operators are clearly visible. But the above arguments remain valid for standard functional linear regression, where for some real-valued (centered) response variable $Y_i$ with $\V(Y_i)<\infty$ one aims to determine the best linear functional $\tilde{L}:\mathbb{L}^2(O)\to\mathbb{R}$ according to the model $Y_i=\tilde{L}(X_i^O)+\varepsilon_i$. Straightforward generalizations of Theorems \ref{WellDef} and \ref{OptimalPrediction} below then show that the optimal functional $\tilde{\mathcal{L}}(X_i^O)$ is given by
\begin{equation*}
\tilde{\mathcal{L}}(X_i^O)=\langle\sigma,X_i^O\rangle_H,
\end{equation*}
where $\sigma(u):=\E(Y_iX_i^O(u))$ for $u\in O$. Following the arguments of Example 2 it is immediately seen that it constitutes a restrictive, additional condition, to assume that $\tilde{\mathcal{L}}(X_i^O)$ can be rewritten in the form $L(X_i^O)(u)=\langle\beta,X_i^O\rangle_2$ for some $\beta_u\in\mathbb{L}^2(O)$.
}

\subsection{Theoretical properties}
Result \textit{(a)} of the following theorem assures that $\mathcal{L}$ is a reconstruction operator according to Definition \ref{def:SRO}, and result \textit{(b)} assures unbiasedness.
\begin{theorem}\label{WellDef}
Let the (centered) random function $X_i^O$ have a KL representation as in \eqref{eq:KLR}.
\begin{itemize}
\item[(a)]$\mathcal{L}(X_i^O)(u)$ in \eqref{eq:PEFexp} has a continuous and finite variance function, i.e., $\V(\mathcal{L}(X_i^O)(u))<\infty$ for all $u\in M$.
\item[(b)]$\mathcal{L}(X_i^O)(u)$ is unbiased in the sense that $\E(\mathcal{L}(X_i^O)(u))=0$ for all $u\in M$.
\end{itemize}
\end{theorem}

The following theorem describes the fundamental properties of the reconstruction error
\begin{align*}
  \mathcal{Z}_i:=X_i^M-\mathcal{L}(X_i^O),\quad\mathcal{Z}_i\in\mathbb{L}^2(M),
\end{align*}
and contains the optimality result for our reconstruction operator $\mathcal{L}$. Result \textit{(a)} shows that the reconstruction error $\mathcal{Z}_i$ is orthogonal to $X^{O}_i$. This result serves as an auxiliary result for result \textit{(b)} which shows that $\mathcal{L}(X_i^O)$ is the optimal linear reconstruction of $X_i^M$. Finally, result \textit{(c)} allows us to identify cases where $X_i^M$ can be reconstructed without any reconstruction error.

\begin{theorem}[Optimal linear reconstruction]\label{OptimalPrediction}
Under our setup it holds that:
\begin{itemize}
\item[(a)] For every $v\in O$ and $u\in M$,
\begin{align}
  &\E\left(X^{O}_i(v)\mathcal{Z}_i(u)\right)=0 \quad\text{and}\label{Zcorr}\\
  &\V(\mathcal{Z}_i(u))=\E\left((\mathcal{Z}_i(u))^2\right)=\gamma(u,u)-\sum_{k=1}^\infty\lambda_k^{O}\big(\tilde{\phi}_{k}^O(u)\big)^2.\label{ZVar}
  \end{align}
  \item[(b)] For any linear operator
    $L:\mathbb{L}^2(O)\rightarrow\mathbb{L}^2(M)$ that is a reconstruction operator with respect to $X_i^O$, according to Definition \ref{def:SRO},
    \begin{align*}
      \E\Big(\big(X_i^M(u)-L(X_i^{O})(u)\big)^2 \Big)\geq\V(\mathcal{Z}_i(u)),\quad\text{for all}\quad u\in M.
    \end{align*}
  \item[(c)] Assume that the underlying process $X_i$ is Gaussian, and let $X_{i,1}$ and $X_{i,2}$ be two independent copies of the random variable $X_i$.
  Then for all $u\in M$ the variance of the reconstruction error can be written as
    \begin{align}\label{Zgauss}
      \V(\mathcal{Z}_i(u))=\frac{1}{2}
      \E\Big(\E\Big(\big(X_{i,1}(u)-X_{i,2}(u)\big)^2 \bigl|
          X_{i,1}^{O}=X_{i,2}^{O}\Big)\Big)
    \end{align}
     where $X_{i,1}^{O}= X_{i,2}^{O}$ means that $X_{i,1}(v)=X_{i,2}(v)$ for all $v\in O$.
  \end{itemize}
\end{theorem}

Whether or not a sensible reconstruction of partially observed functions is possible, of course, depends on the character of the underlying process. For very rough and unstructured processes no satisfactory results can be expected. An example is the standard Brownian motion on $[0,1]$  which is a pure random process with independent increments. If Brownian motions $X_i$ are only observed on an interval $O:=[0,\vartheta]$, it is well known that the ``best'' (and only unbiased) prediction of $X_i(u)$ for $u\in M:=(\vartheta,1]$ is the last observed value $X_i^O(\vartheta)$. This result is consistent with our definition of an ``optimal'' operator $\mathcal{L}$: The covariance function of the Brownian motion is given by $\gamma_u(v)=\min(u,v)$, and hence for all $v\in [0,\vartheta]$ and $u\geq \vartheta$ one obtains $\gamma_u(v)=\min(u,v)=\min(\vartheta,v)=\gamma_\vartheta(v)=v$. Therefore, by \eqref{eq:PF_S} and \eqref{rcoex1} we have $\mathcal{L}(X_i^O)(u)=\langle\gamma_u,X_i^O\rangle_H=\langle\gamma_\vartheta,X_i^O\rangle_H=X_i^O(\vartheta)$ for all $u\in[\vartheta,1]$. Although in this paper we focus on processes that lead to smooth, regularly shaped sample curves, the Brownian motion is of some theoretical interest since it defines a reconstruction operator which obviously does not constitute a regression operator. Also note that $\mathcal{L}(X_i^O)(u)=X_i^O(\vartheta)$
will provide  perfect reconstructions if a.s. sample functions $X_i(u)$ are constant for all $u\in M$.


Result \textit{(c)} of Theorem \ref{OptimalPrediction} may be useful to identify cases that allow for a perfect reconstruction. By \eqref{Zgauss} there is no reconstruction error, i.e., $\V(\mathcal{Z}_i(u))=0$ for $u\in M$ if the event $X_i^{O}=X_j^{O}$ implies that also $X_i^M=X_j^M$. This might be fulfilled for very simply structured processes. It is necessarily satisfied for finite dimensional random functions $X_i^K(u)=\sum_{k=1}^K\xi_{ik}\phi_k(u)$, $\lambda_{K+1}=\lambda_{K+2}=\dots=0$, as long as the basis functions $\phi_1,\dots,\phi_K$ are linearly independent over $O$.

\subsection[]{A deeper look at the structure of $\mathcal{L}$}\label{sec:CritCons}
Remember that the definition of $\mathcal{L}$ can be extended to an operator $\mathcal{L}:\mathbb{L}^2(O)\to\mathbb{L}^2(O\cup M)$. For elements $u\in O$ of the observed part $O$ the best ``reconstruction'' of $X_i(u)$ is obviously the observed value $X_i(u)$ itself, and indeed for any $u\in O$ \eqref{eq:PF_S} yields $\mathcal{L}(X_i^O)(u)=\langle\gamma_u,X_i^O\rangle_H=X_i(u)$. Equation \eqref{eq:PEFexp} then holds with
\begin{align*}
\tilde{\phi}_{k}^O(u):=
\frac{\langle\phi_k^{O},\gamma_u\rangle_2}{\lambda_k^O}=\phi_k^{O}(u),\quad u\in O.
\end{align*}

Since $\gamma_u(v)=\gamma(u,v)=\E(X_i(u)X_i(v))$ is a continuous function on $O\cup M$ it follows that the resulting ``reconstructed'' function $[\mathcal{L}(X_i^O)]$ is continuous on $O\cup M$. In particular, $[\mathcal{L}(X_i^O)]$ is continuous at any boundary point $\vartheta_u\in\partial M$, and
\begin{align*}
&\lim_{u\in M, u\rightarrow \vartheta_u}\mathcal{L}(X_i^O)(u)=X_i(\vartheta_u), \text{ as well as }\\
&\lim_{u\in M, u\rightarrow \vartheta_u}\tilde{\phi}_{k}^O(u)=\phi_{k}^O(\vartheta_u), \ k=1,2,\dots
\end{align*}

Equation \eqref{eq:PEFexp} together with our definition of $\mathcal{Z}_i$ imply that the complete function $X_i$ on $O\cup M$ can be represented in the form
\begin{align}\label{completerepresent}
X_i(v)=\sum_{k=1}^\infty \xi^{O}_{ik}\phi^{O}_{k}(v), v\in O,\text{ and }
X_i(u)=\sum_{k=1}^\infty \xi^{O}_{ik}\tilde{\phi}^{O}_{k}(u)+\mathcal{Z}_i(u), u\in M.
\end{align}
This sheds some additional light on result \eqref{Zgauss}. We will  have $\mathcal{Z}_i(u)\approx 0$ and $X_i(u)\approx\sum_{k=1}^\infty \xi^{O}_{ik}\tilde{\phi}^{O}_{k}(u)$ if on the segment $M$ the process is {\em essentially} driven by the same random components $\xi^{O}_{ik}$ as those determining its structure on $O$. Additional random components $\mathcal{Z}_i(u)$, not present on $O$, and uncorrelated with $\xi^{O}_{ik}$, then have to be of minor importance. If the observed interval is sufficiently long, then this may be approximately true for processes with smooth, similarly shaped trajectories. Note that even if $X_i(u)=\sum_{k=1}^\infty \xi^{O}_{ik}\tilde{\phi}^{O}_{k}(u)$ for $u\in M$, the eigenfunctions of $X_i^M$ will usually not coincide with $\tilde{\phi}^{O}_{k}$  for $u\in M$, since there is no reason to expect that these functions are mutually orthogonal.

\section{Estimation}\label{sec:estim}
\begin{figure}[htbp]
\centering
\includegraphics[width=.55\textwidth]{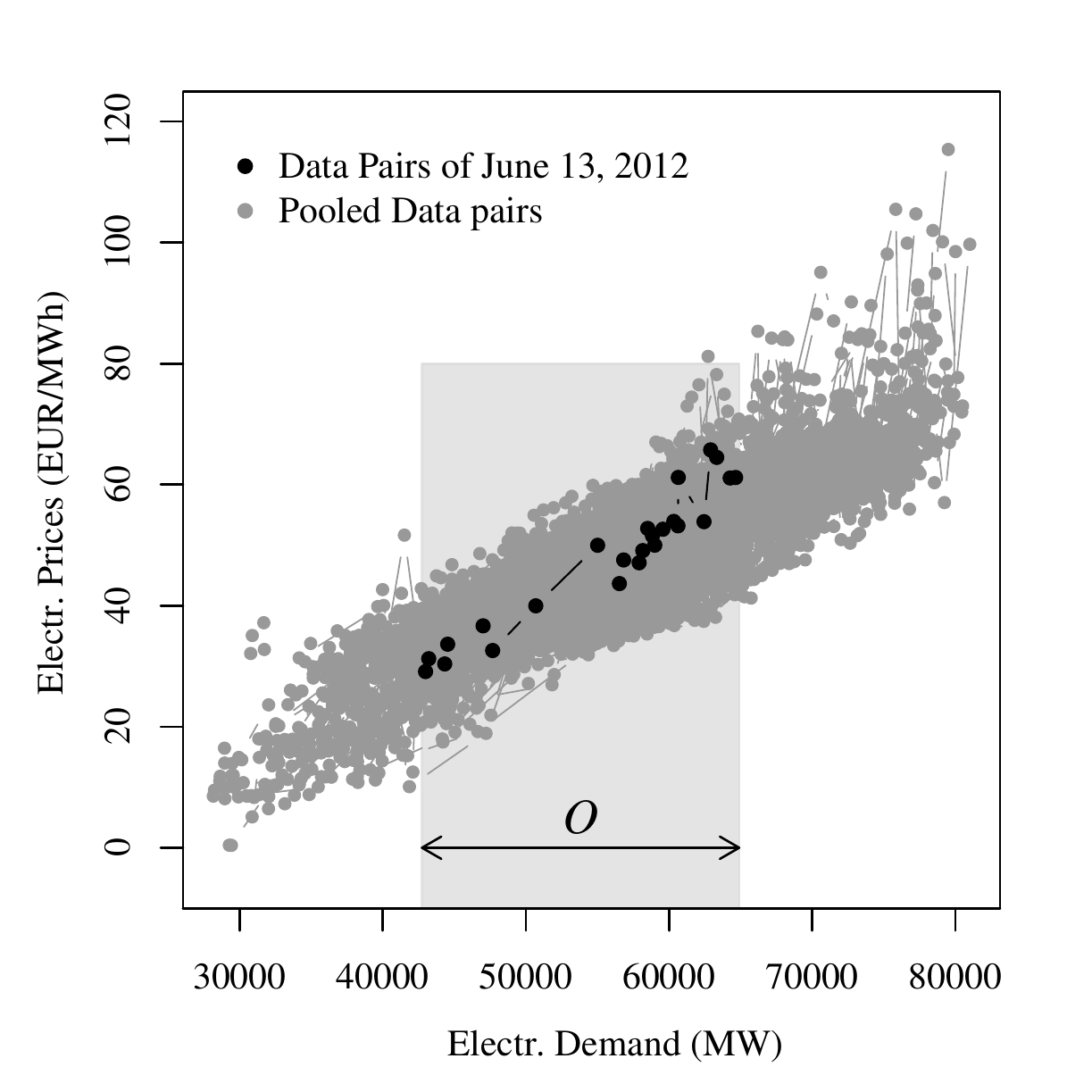}
\caption[]{Scatter plot of the observed data pairs $(Y_{ij},U_{ij})$.}
\label{fig:scatter}
\end{figure}
We typically do not observe a functional trajectory directly, but only its discretization with or without measurement errors. For instance, Figure \ref{fig:presm_scatter} shows the \emph{pre-smoothed} functions; however, the actual raw data is shown in Figure \ref{fig:scatter}. Let $\mathbb{X}_i^O:=((Y_{i1},U_{i1}),\dots,(Y_{im_{i}},U_{im_{i}}))$ denote the observable data pairs of a function $X_i^O$, where
\begin{align}\label{DGP}
Y_{ij}=X_i^O(U_{ij})+\varepsilon_{ij},\quad i\in\{1,\dots,n\},\quad j\in\{1,\dots,m_{i}\},
\end{align}
and $U_{ij}\in O_i$.

For the rest of the paper, we focus on the case where $O_i=[A_i,B_i]$ as in our real data application. However, we give detailed descriptions on how to use our methods in the more general cases where $O_i$ consists of several mutually disjoint subintervals. We consider the case where $U_{i1},\dots,U_{im_{i}}$ are iid random variables with strictly positive density over the random subinterval $[A_i,B_i]$, which in practice can be approximated by $A_i\approx\min_{1\leq j\leq m_{i}}(U_{ij})$ and $B_i\approx\max_{1\leq j\leq m_{i}}(U_{ij})$. Let the error term $\varepsilon_{ij}$ be a real iid random variable that is independent from all other stochastic model components and has mean zero and finite (possibly zero) variance $\mathbb{V}(\varepsilon_{ij})=\sigma^2$ with $0\leq\sigma^2<\infty$. Motivated by our real data application we will concentrate on the case that $n$ is considerably larger than $m_{i}$, which also holds in many other important applications.

So far, we have considered centered random functions $X_i^O$. Henceforth, we consider \textbf{non-centered} functions and will make the empirical centering explicit in all estimators. As already outlined in Section \ref{sec:intro}, we propose to estimate $\mathcal{L}(X_i^O)(u)$ by the empirical counterpart of the truncated sum
\begin{align}\label{Ltrunc}
\mathcal{L}_{K}(X_i^O)(u)&=\mu(u)+\sum_{k=1}^K \xi_{ik}^O\tilde{\phi}^O_k(u)=
\mu(u)+
\sum_{k=1}^K \xi_{ik}^O\, \frac{\langle\phi_k^O,\gamma_u\rangle_2}{\lambda_k^{O}},
\end{align}
where the unknown true values of $\xi_{ik}^O$ and $\tilde{\phi}^O_k(u)$ are replaced by suitable estimates defined below.

Remember, however, that our data structure in \eqref{DGP} implies that we are faced with two simultaneous estimation problems. One is efficient estimation of $\mathcal{L}(X_i^O)(u)$ for $u\in M$, the other one is the estimation of the underlying function $X_i(u)$ for $u\in O$. There are two possible strategies which can be employed.

The first is motivated by the best basis property \eqref{eq:best basis} and simply consists in using an FPCA-approximation of $X_i$ on $O$. Recall that $\mathcal{L}(X_i^O)(u)$ can be extended to an operator on $O\cup M$. For $u\in O$ we then obtain $\tilde{\phi}_k^O(u)=\langle\phi_k^O,\gamma_u\rangle_2/\lambda_k^{O}=\phi_k^O(u)$, and thus $\mathcal{L}(X_i^O)(u)=X_i(u)$. That is, estimates $\widehat{\mathcal{L}}_{K}(\mathbb{X}_i^O)(u)$ of $\mathcal{L}_{K}(X_i^O)(u)$ for $u\in O\cup M$ will simultaneously provide estimates of the true function $X_i(u)$ on the observed interval $O$ and of the optimal reconstruction $\mathcal{L}(X_i^O)(u)$ on the unobserved interval $M$.

The second approach is to rely on nonparametric curve estimation, e.g., local linear smoothers, to approximate $X_i^O$ on $O$, while \eqref{Ltrunc} is only used for reconstructing the unobserved part $M$. We then, however, run into the boundary problem already mentioned in the introduction. Let $\vartheta_u$ be the boundary point closest to the considered $u\in M$, i.e., $\vartheta_u=A_i$ if $|A_i-u|<|B_i-u|$ and $\vartheta_u=B_i$ else. Usually nonparametric estimates of $X_i^O$ and reconstruction estimates based on \eqref{Ltrunc} will not coincide for $u=A_i$ or $u=B_i$. A correction, leading to continuous function estimates on $O\cup M$ may then be based on the identity $\mathcal{L}(X_i^O)(u)=X_i^O(\vartheta_u)+\mathcal{L}(X_i^O)(u)-\mathcal{L}(X_i^O)(\vartheta_u)$ and its truncated version
\begin{align}\label{Ltruncast}
\mathcal{L}_{K}^*(X_i^O)(u)&=X_i^O(\vartheta_u)+\mathcal{L}_{K}(X_i^O)(u)-\mathcal{L}_{K}(X_i^O)(\vartheta_u)
\nonumber \\
&=
\mu(u)-\mu(\vartheta_u)+
\sum_{k=1}^K \xi_{ik}^O\Big(\tilde{\phi}^O_k(u)-\tilde{\phi}^O_k(\vartheta_u)\Big), \quad\text{for } u\in M
\end{align}

In this paper we propose to use the following empirical counterparts of $\mathcal{L}_K(X_i^O)(u)$ and $\mathcal{L}_{K}^*(X_i^O)(u)$:
\begin{align}
&\widehat{\mathcal{L}}_{K}(\mathbb{X}_i^O)(u):=\hat{\mu}(u;h_\mu)+
\sum_{k=1}^K \hat{\xi}_{ik}^O\hat{\tilde{\phi}}^O_k(u) \quad\text{for } u\in O\cup M, \label{Estim_L} \\
&\quad \text{with } \hat{\tilde{\phi}}^O_k(u):= \frac{\langle\hat{\phi}_k^O,\hat{\gamma}_u\rangle_2}{\hat{\lambda}_k^{O}}, \quad k=1,\dots,K,
\nonumber
\end{align}
\begin{align}
\widehat{\mathcal{L}}_{K}^*(\mathbb{X}_i^O)(u):&=\widehat{X}_i^{O}(\vartheta_u;h_X)
+\widehat{\mathcal{L}}_{K}(\mathbb{X}_i^O)(u)-\widehat{\mathcal{L}}_{K}(\mathbb{X}_i^O)(\vartheta_u)
\nonumber \\
&\hspace*{-.75cm}=\widehat{X}_i^{O}(\vartheta_u;h_X)+
\hat{\mu}(u;h_\mu)-\hat{\mu}(\vartheta_u;h_\mu)+
\sum_{k=1}^K \hat{\xi}_{ik}^O\,\left( \hat{\tilde{\phi}}^O_k(u)-\hat{\phi}_k^O(\vartheta_u)\right),
\label{Estim_D}
\end{align}
where
$\widehat{X}_i^{O}$ denotes the LLK estimator of $X_i^{O}$ (see \eqref{X_cent_Estim}),
$\hat{\mu}$ denotes the LLK estimator of $\mu$ (see \eqref{Estimator_mu_paper}),
$\hat{\gamma}_u$ denotes the LLK estimator of the covariance function (see \eqref{Estimator_gamma_paper}),
$\hat{\phi}^O_k$ and $\hat{\lambda}_k^{O}$ denote the estimators of the eigenfunctions and eigenvalues (see \eqref{small_eigen}), and
$\hat{\xi}_{ik}$ denote the estimators of the pc-scores (see \eqref{PC-Score-Estim}).

\paragraph{\normalfont Remark}
\textit{
Estimator \eqref{Estim_L} can be directly applied in the general case, where $O=\bigcup_{j=1}^J[A_j,B_j]$ consists of a union of finitely many mutually disjoint subintervals $[A_j,B_j]\subseteq[a,b]$. Estimator \eqref{Estim_D}, however, must be adjusted for this general case as follows.
First, consider a point $u\in M$ located between the observed intervals $[A_j,B_j]$ and $[A_{j+1},B_{j+1}]$ for any $j=1,\dots,J-1$. In this case the quantities $\widehat{X}_i^{O}(\vartheta_u;h_X)$ and  $\hat{\phi}_k^O(\vartheta_u)$ in \eqref{Estim_D} have to be replaced by the linear interpolations
$(1-w_u) \widehat{X}_i^{O}(B_j;h_X)+
    w_u  \widehat{X}_i^{O}(A_{j+1};h_X)$
and
$(1-w_u)\hat{\phi}_k^O(B_j) +
    w_u \hat{\phi}_k^O(A_{j+1})$
with $w_u=(u-B_j)/(A_{j+1}-B_j)$.
Second, for $0\leq u<A_1$ replace
$\widehat{X}_i^{O}(\vartheta_u;h_X)$ and
$\hat{\phi}_k^O(\vartheta_u)$ by
$\widehat{X}_i^{O}(A_1;h_X)$ and
$\hat{\phi}_k^O(A_1)$.
Third, for $B_J< u\leq 1$ replace
$\widehat{X}_i^{O}(\vartheta_u;h_X)$ and
$\hat{\phi}_k^O(\vartheta_u)$ by
$\widehat{X}_i^{O}(B_J;h_X)$ and
$\hat{\phi}_k^O(B_J)$.
}

\textit{In our asymptotic analysis (Section \ref{sec:asymp}) we focus on the case of single subintervals $O_i=[A_i,B_i]$ which leads to comprehensible theorems and proofs.}


\smallskip

For $u\in O$ the LLK estimator $\widehat{X}_i^{O}(u;h_X)$ is defined by $\widehat{X}_i^{O}(u;h_X)=\hat\beta_0$, where
\begin{align}\label{X_cent_Estim}
  (\hat\beta_0,\hat\beta_1)=\underset{\beta_0,\beta_1}{\arg\min}\sum_{j=1}^{m_{i}}[Y_{ij}-\beta_0-\beta_1(U_{ij}-u)]^2 K_{h_X}(U_{ij}-u)
\end{align}
for $K_{h}(.)=\kappa(./h)/h$. The kernel function $\kappa$ is assumed to be a univariate symmetric pdf with compact support $\supp(\kappa)=[-1,1]$ such as, e.g., the Epanechnikov kernel (see Assumption A5). The usual kernel constants are given by $\nu_{2}(\kappa):=\int v^2\kappa(v)dv$, and $R(\kappa):=\int\kappa(v)^2dv$.

The LLK mean estimator $\hat{\mu}(u;h_{\mu})$ is defined by $\hat{\mu}(u;h_{\mu})=\hat\beta_0$, where
\begin{align}\label{Estimator_mu_paper}
  (\hat\beta_0,\hat\beta_1)=\underset{\beta_0,\beta_1}{\arg\min}\sum_{i=1}^n\sum_{j=1}^{m_{i}}[Y_{ij}-\beta_0-\beta_1(U_{ij}-u)]^2 K_{h_\mu}(U_{ij}-u).
\end{align}

The LLK estimator $\hat{\gamma}_u(v)=\hat\gamma(u,v;h_\gamma)$ is defined as $\hat{\gamma}(u,v;h_{\gamma})=\hat\beta_0$, where
\begin{align}\label{Estimator_gamma_paper}
  (\hat{\beta}_0,\hat\beta_1,\hat\beta_2)=&\underset{\beta_0,\beta_1,\beta_2}{\arg\min}\sum_{i=1}^n\sum_{1\leq j,l\leq m_{i}}[\widehat{C}_{ijl}-\beta_0-\beta_1(U_{ij}-u)-\beta_2(U_{il}-v)]^2\\
  &\times K_{h_\gamma}(U_{ij}-u)K_{h_\gamma}(U_{il}-v),\notag
\end{align}
with raw-covariance points
$\widehat{C}_{ijl}=(Y_{ij}-\hat{\mu}(U_{ij}))(Y_{il}-\hat{\mu}(U_{il}))$. Like \cite{Yao2005}, we do not include the diagonal raw-covariances $\widehat{C}_{ijj}$ for which $U_{ij}=U_{ij}$ as these would introduce an estimation bias through taking squares of the error term $\varepsilon_{ij}$ contained in $Y_{ij}$.

Estimates of the eigenvalues $\lambda_k^{O}$ and the eigenfunctions $\phi_k^{O}$ are defined by the corresponding solutions of the empirical eigen-equations
\begin{align}
\int_{O}\hat{\gamma}(u,v;h_\gamma)\hat{\phi}^{O}_k(v)\,dv&=\hat{\lambda}^{O}_k\,\hat{\phi}^{O}_k(u),\quad u\in O.\label{small_eigen}
\end{align}

\paragraph{\normalfont Remark}
\textit{The implementation of \eqref{small_eigen} can be done as usually by discretizing the smoothed covariance $\hat{\gamma}(u_r,v_s)$ using regular grid points $(u_r,v_s)\in[a,b]^2$, $r,s\in\{1,\dots,L\}$ \citep[see, for instance,][]{rice1991estimating}. For approximating the eigenvalues and eigenfunctions of $\Gamma^O(x)(u)=\int\gamma^O(u,v)x(v)dv$ one needs to construct the matrix $(\hat\gamma^O(u_r,v_s))_{r,s}$ from the grid points falling into $[A,B]\times[A,B]$. In the case of several disjoint intervals the matrix must be assembled from the grid points falling into the intervals $[A_j,B_j]\times [A_{j'},B_{j'}]$, $j,{j'}\in\{1,\dots,J\}$.}

Finally, the empirical pc-score $\hat{\xi}^{O}_{ik}$ is defined by the following integral approximation of $\xi_{ik}^O$:
\begin{align}\label{PC-Score-Estim}
\hat{\xi}^O_{ik}=\sum_{j=2}^{m_i}\hat{\phi}_k^O(U_{i(j)})
(Y_{i(j)}-\hat{\mu}(U_{i(j)};h_\mu))(U_{i(j)}-U_{i,(j-1)}),
\end{align}
where $(Y_{i(j)},U_{i(j)})$ are ordered data pairs for which the ordering is determined through the order sample $U_{i(1)}\leq\dots\leq U_{i(m_{i})}$

In our theoretical analysis we consider $K\equiv K_{nm}\to\infty$ as the sample size $nm\to\infty$, where $m\leq m_i$ for all $i=1,\dots,n$. In practice, the truncation parameter $K$ can be chosen by one of the usual procedures such as, for instance, Cross Validation or the Fraction of Variance Explained (FVE) criterion.

Alternatively, one can use an adapted version of the GCV criterion in \cite{kraus2015} in order to define an $M$-specific GCV criterion. For this let $\mathcal{C}$ denote the index set of the completely observed functions $\mathbb{X}_l$, $l\in\mathcal{C}$, with $[a,b]\approx[\min_{1\leq j\leq m_i}(U_{lj}),\max_{1\leq j\leq m}(U_{lj})]$, for instance, with $\min_{1\leq j\leq m_i}(U_{lj})\in[a,a+(b-a)/10]$ and $\min_{1\leq j\leq m_i}(U_{lj})\in[b-(b-a)/10,b]$ and define the following vectors by partitioning the complete data-vectors into pseudo-missing and pseudo-observed parts:
\begin{align*}
  \mathbf{Y}_l^{M} &=\big(Y_{lj}:j=1,\dots,m_l;\; U_{lj}\in M\big)^\top,\\
  \mathbb{X}_l^{O} &=\big((Y_{lj},U_{lj}):j=1,\dots,m_l;\; U_{lj}\in O\big)^\top,\text{ and}\\
  \widehat{\mathbf{Y}}_{lK}^{M} &=\big(\widehat{\mathcal{L}}_{K}(\mathbb{X}_l^{O})(U_{lj}):j=1,\dots,m_l;\; U_{lj}\in M\big)^\top.
\end{align*}
This allows us to compute the weighted sum of the residual sum of squares $||\mathbf{Y}_l^M - \widehat{\mathbf{Y}}_{lK}^M||^2$ for reconstructions over $M$
\begin{align*}
\operatorname{RSS}_M(K)=\sum_{l\in C}||\mathbf{Y}_l^M - \widehat{\mathbf{Y}}_{lK}^M||^2/|\mathbf{Y}_l^M|,
\end{align*}
where $|\mathbf{Y}_l^M|$ is the number of elements in $\mathbf{Y}_l^M$. The GCV criterion for reconstructing functions over $M$ is
\begin{align}\label{GCV}
  \operatorname{GCV}_M(K)=\frac{\operatorname{RSS}_M(K)}{\big(1-K/|\mathcal{C}|\big)^2},
\end{align}
where $|\mathcal{C}|$ is the number of elements in $\mathcal{C}$, i.e., the number of complete functions.

\section{Asymptotic results}\label{sec:asymp}
Our theoretical analysis analyzes the reconstruction of an arbitrary sample function $X_i$ satisfying $O\subseteq O_i=[A_i,B_i]$.

Our asymptotic results on the convergence of our nonparametric estimators are developed under the following assumptions which are generally close to those in \cite{yao2005functional} and \cite{hall2006properties}. We additionally allow for weakly dependent time series of random functions $(X_i)_i$, and we consider a different asymptotic setup excluding the case of sparse functional data. Only second-order kernels are employed.\\[-1.5ex]

\noindent\textbf{A1} (Stochastic)
For some $d_{\min}>0$ the conditional random variables\\ $U_{i1}|O_i$,\dots,$U_{im}|O_i$ are iid with pdf $f_{U|O_i}(u)\geq d_{\min}$ for all $u\in O_i=[A_i,B_i]$ and zero else. For the marginal pdf $f_U$ it is assumed that $f_{U}(u)>0$ for all $u\in[a,b]$ and zero else. The time series $(A_i)_{i=1,\dots,n}$, $(B_i)_{i=1,\dots,n}$, and $(X_i)_{i=1,\dots,n}$ are strictly stationary ergodic (functional) time series with finite fourth moments (i.e., $\E(||X_i||_2^4)<\infty$ in the functional case) and autocovariance functions with geometric decay. I.e., there are constants $C_A,C_B,C,\dot{C},\iota_A,\iota_B,\iota,\dot{\iota}$ with $0<C_A,C_B,C,\dot{C}<\infty$ and $0<\iota_A,\iota_B,\iota,\dot{\iota}<1$, such that
$|\Cov(A_i,B_{i+h})|\leq C_A \iota_A^{h}$,
$|\Cov(B_i,B_{i+h})|\leq C_B \iota_B^{h}$, \\
$\sup_{(u,v)\in [a,b]^2}|\gamma_h(u,v)|\leq C \iota^{h}$, and \\
$\sup_{(u_1,v_1,u_2,v_2)\in[a,b]^4}|\dot\gamma_h((u_1,v_1),(u_2,v_2))|\leq \dot{C}\dot{\iota}^{h}$ for all $h\geq 0$, where\\
$\gamma_h(u,v):=\Cov(X_{i+h}(u),X_i(v))$ and\\
$\dot\gamma_h((u_1,v_1),(u_2,v_2)):=\Cov(X_{i+h}(u_1)X_{i+h}(v_1),X_i(u_2)X_i(v_2))$.\\
The error term $\varepsilon_{ij}$ is assumed to be independent from all other random variables. The random variables $U_{ij}$ and $O_i$ are assumed to be independent from $(X_i)_{i=1,\dots,n}$, which leads to the so-called ``missing completely at random'' assumption. The event $O_i\times O_i=[a,b]^2$ has a strictly positive probability and $B_i>A_i$ almost surely.
\\[1.5ex]
\noindent\textbf{A2} (Asymptotic scenario) $nm\to\infty$ with $m\leq m_i$ for all $i=1,\dots,n$, where $n\to\infty$ and $m=m(n)\asymp n^\theta$ with $0<\theta<\infty$. Here, $a(n)\asymp b(n)$ is used to denote that $(a(n)/b(n)) \to c$ as $n\to\infty$, where $c$ is some constant $0<c<\infty$.\\[1.5ex]
\noindent\textbf{A3} (Smoothness)
For $\hat{\mu}$: All second order derivatives of
$\mu(u)$ on $[a,b]$,
$f_{U}(u)$ on $[a,b]$,
$\gamma(u,v)$ on $[a,b]^2$, and of
$f_{YU}(y,u)$ on $\mathbb{R}\times[a,b]$
are uniformly continuous and bounded, where $f_{YU}$ is the joint pdf of $(Y_{ij},U_{ij})$.
For $\hat{\gamma}$: All second order derivatives of $\gamma(u,v)$ on
$[a,b]^2$, $f_{UU}(u,v)$ on $[a,b]^2$,
$\dot\gamma((u_1,v_1),(u_2,v_2))$ on $[a,b]^4$, and of $f_{CUU}(c,u,v)$ on
$\mathbb{R}\times[a,b]^2$ are uniformly continuous and bounded, where
$f_{CUU}$ is the joint pdf of $(C_{ijl},U_{ij},U_{il})$.
Finally, $f_{U|O_i}(u)$ is a.s.~continuously differentiable, and
$\E\big(|f'_{U|O_i}(u)|/f_{U|O_i}(u)^2\big)<\infty$, and $X_i$ is a.s.~twice continuously differentiable. \\[1.5ex]
\noindent\textbf{A4} (Bandwidths)
For estimating $X_i^O$:
$h_{X}\to 0$ and
$(m\,h_{X})\to \infty$ as $m\to\infty$.
For estimating $\mu$:
$h_{\mu}\to 0$ and
$(nm\,h_{\mu})\to \infty$ as $nm\to\infty$.
For estimating $\gamma$:
$h_{\gamma}\to 0$ and
$(n\mathcal{M}\,h_{\gamma})\to \infty$ as $n\mathcal{M}\to\infty$,
where $\mathcal{M}=m^2-m$.\\[1.5ex]
\noindent\textbf{A5} (Kernel function)
$\kappa$ is a second-order kernel with compact support $\supp(\kappa)=[-1,1]$.
\smallskip


In Assumption A2, we follow \cite{ZC2007} and consider a deterministic sample size $m\to\infty$, where $m\leq m_i$ for all $i=1,\dots,n$. As \cite{hall2006properties}, \cite{ZC2007} and \cite{zhang2016sparse} we do not consider random numbers $m_i$, but if $m_i$ are random, our theory can be considered as conditional on $m_i$.

While A1-A5 suffice to determine rates of convergence of mean and covariance estimators, it is well-known from the literature that rates of convergence of estimated eigenfunctions will depend on the rate of decay characterizing the convergence of $\lambda_k^O$ to zero as $k\rightarrow\infty$.

We want to note that for a subinterval $O\subset [a,b]$ the decay of eigenvalues $\lambda^O_1,\lambda_2^O,\dots$ will usually be faster than the rate of decay of the eigenvalues $\lambda^C_1,\lambda^C_2,\dots$ of the complete covariance operator defined on $[a,b]^2\supset O^2$. This is easily seen. Let $\gamma_1^C,\gamma_2^C,\dots$ denote the corresponding eigenfunctions on $[a,b]$, and define $\gamma_k^{C|O}\in \mathbb{L}^2(O)$ by $\gamma_k^{C|O}(u)=\gamma_k(u)$ for $u\in O$ and $k=1,2,\dots$. For the special case $v_k=\gamma_k^{C|O}$, $k=1,\dots,K$, inequality \eqref{eq:best basis} then implies that for all $K\geq 1$ we have  $\sum_{k=K+1}^\infty\lambda^O_k\leq \sum_{k=K+1}^\infty\lambda^C_k \int_O \gamma_k^{C|O}(u)^2du\leq \sum_{k=K+1}^\infty\lambda^C_k$, since $\int_O \gamma_k^{C|O}(u)^2du\leq \int_a^b \gamma_k^{C}(u)^2du=1$ for all $k=1,2,\dots$.

To complete our asymptotic setup, we consider the reconstruction of arbitrary sample functions $X_i$ observed over an interval $O_i=[A_i,B_i]$ with  length $B_i-A_i\geq \ell_{\min}$, where $0<\ell_{\min}<b-a$ is an (arbitrary) constant. We then impose the following additional assumptions.

\smallskip

\noindent\textbf{A6} (Eigenvalues)
For any subinterval $O=[A,B]\subset [a,b]$ with $B-A\geq \ell_{\min}$ the ordered eigenvalues $\lambda_1^O>\lambda_2^O>\dots>0$ have all multiplicity one. Furthermore, there exist some $a_O>1$ and some $0<c_O<\infty$, possibly depending on $O$, such that $\lambda^O_k-\lambda^O_{k+1}\geq c_O k^{-a_O-1}$ with $0<c_O<\infty$, and $\lambda_k^O=\mathcal{O}(k^{-a_O})$ as well as $1/\lambda_k^O =O(k^{a_O})$ as $k\rightarrow\infty$.\\[1.5ex]
\noindent\textbf{A7} (Eigenfunctions)
For any subinterval $O=[A,B]\subset [a,b]$ with $B-A\geq \ell_{\min}$ there exists a constant $0<D_O<\infty$ such that $\sup_{u\in [a,b]} \sup_{k\ge 1} |\tilde{\phi}^O_k(u)|\leq D_O $ (recall that $\tilde{\phi}^O_k(U)= \phi^O_k(u)$ for $u\in O$).
\smallbreak

Assumption A6 requires a polynomial decay of the sequence of eigenvalues. It cannot be tested, but it corresponds to the usual assumption characterizing a majority of work concerning eigenanalysis of functional data, although some authors also consider exponential decays. There exist various types of functional data, but this paper focuses on applications where the true sample functions are  smooth and all possess a similar functional structure. This is quite frequent in practice, and in applied papers it is then often found that few functional principal components suffice to approximate sample functions with high accuracy. In view of the best basis property \eqref{eq:best basis} one may then tend to assume that A6 holds for some very large $a_O\gg 1$. Indeed, for increasing $k$ eigenfunctions $\phi^O_k$ will become less and less ``smooth'' since the number of sign changes will necessarily tend to infinity. If observed trajectories are  smooth, then the influence of such high-frequency components must be very small, indicating a very small eigenvalue $\lambda_{k}^O=\E(\xi_k^O)$ for large $k$. This is of substantial interest, since the theorems below show that rates of convergence of our final estimators are better the larger $a_O$.

Assumption A7 imposes a (typical) regularity condition on the structure of the eigenfunctions $\phi^O_k(u)$, since $\tilde{\phi}^O_k(u)=\phi^O_k(u)$ for $u\in O$. For $u\in M=[a,b]\setminus O$ condition $|\tilde{\phi}^O_k(u)|\leq D_O$ is much weaker than the standard assumption of a regression operator which would go along with the requirement $\sum_{k=1}^\infty\tilde{\phi}^O_k(u)^2<\infty$. But, for $u\in M$, theory only ensures that $\sum_{k=1}^\infty\lambda_k(\tilde{\phi}^{O}_{k}(u))^2<\infty$ (see Theorem \ref{OptimalPrediction} (a)) and A7 is restrictive in so far as it excludes the possible case that for $u\in M$ we have $|\tilde{\phi}^{O}_{k}(u)|\rightarrow \infty$ as $k\rightarrow \infty$. We are not sure whether the latter excluded case constitutes a realistic scenario in practical applications, since by \eqref{completerepresent} it would correspond to the fairly odd situation that for large $k$ the high-frequency components $\xi^{O}_{ik}$ possess much larger influence on $M$ than on $O$. Nevertheless, we want to emphasize that the arguments used in the proof of our theorems may easily be generalized to prove consistency of our estimators even in this excluded case; however, rates of convergence deteriorate and asymptotic expressions become much more involved.


\begin{theorem}[Preliminary consistency results]\label{Theorem_UR}\ \\
Under Assumptions A1-A5 we have that:
\begin{itemize}
\item[(a)]$\sup_{u\in[a,b]}|\hat{\mu}(u;h_\mu)-\mu(u)|=\mathcal{O}_p\left(r_\mu\right)$
\item[(\~a)]Conditional on $X_i^O$:\;
$\sup_{u\in O}|\widehat{X}_i^O(u;h_\mu,h_X)-X_i^O(u)|=\mathcal{O}_p\left(r_X\right)$
\item[(b)]$\sup_{(u,v)\in [a,b]^2}|\hat{\gamma}(u,v;h_\gamma)-\gamma(u,v)|=\mathcal{O}_p\left(r_\mu+r_\gamma\right)$, where
\end{itemize}
\begin{equation*}
\hspace*{-5.5ex}\begin{array}{rcl}
r_\mu&\equiv& r_\mu(h_\mu,n,m):=h_\mu^2+1/\sqrt{nm\,h_\mu}+1/\sqrt{n}\\
r_X&\equiv& r_X(h_X,m):=h_X^2+1/\sqrt{m\,h_X}\\
r_\gamma&\equiv& r_\gamma(h_\gamma,n,\mathcal{M}):=h_\gamma^2+1/\sqrt{n\mathcal{M}\,h^2_\gamma}+1/\sqrt{n},
\end{array}
\end{equation*}
and where $\mathcal{M}=m^2-m$ and $m\leq m_i$ for all $i=1,\dots,n$ (see A2 and A4).\\

\noindent If additionally Assumption A6 and A7 hold, we obtain for every subinterval  $O=[A,B]\subset [a,b]$ with $B-A\geq \ell_{\min}$:
\begin{itemize}
\item[(c)]$\sup_{k\geq 1}|\hat{\lambda}_k^O-\lambda_k^O|=\mathcal{O}_p\left(r_\mu+r_\gamma\right)$\quad for all\quad $k\geq 1$
\item[(d)]
$\sup_{1\leq k\leq K}\delta_k^O \Vert \hat{c}_k\hat{\phi}_k^O-\phi_k^O\Vert_2=\mathcal{O}_p\left(r_\mu+r_\gamma\right)$\\
\end{itemize}
where $\hat{c}_k:=\operatorname{sgn}(\langle\hat{\phi}_k^{O},\phi_k^{O}\rangle_2)$ and
$\delta_k^O:=\min_{j\neq k}\{\lambda_j^O-\lambda_{k}^O\}$.
\end{theorem}
Related results can be found in \cite{Yao2005}, \cite{li2010uniform}, and \cite{zhang2016sparse}. Our proof of results (a)-(b) follows that of \cite{Yao2005}, but is more restrictive as we allow only for compact second order kernels. Results (c) and (d) follow from standard arguments as used in \cite{bosq2000linear}.

\begin{theorem}[Consistency results for $\widehat{\mathcal{L}}_{K}(\mathbb{X}_i^O)$]\label{Theorem_Main_Estim}
Consider an arbitrary $i\in {1,\dots,n}$ and assume that $O=[A,B]\subseteq[a,b]$ satisfies $B-A\geq \ell_{\min}>0$. For some $0<C<\infty$ let $\bar{K}_{mn}=C\cdot  (\min\{n^{1/2},(n\mathcal{M})^{1/3}\})^{1/(a_O+3/2)}$. The following results hold then under Assumptions A1-A7, for $1\leq K\leq \bar{K}_{mn}$, $h_X\asymp m^{-1/5}$,  $h_\mu\asymp (nm)^{-1/5}$ and $h_\gamma\asymp (n\mathcal{M})^{-1/6}$, as $n\to\infty$ and $m\to\infty$ with $m\asymp n^{\theta}$, $0<\theta<\infty$. For any $u\in [a,b]$:

\smallskip
\begin{align}
&\widehat{\mathcal{L}}_{K}(\mathbb{X}_i^O)(u)=\mathcal{L}_{K}(X_i^O)(u)+
\mathcal{O}_p\left(K\left(\frac{1}{m^{1/2}}+\frac{K^{a_O/2+3/2}}{\min\{n^{1/2},(n\mathcal{M})^{1/3}\}}\right)\right)
\nonumber \\
&\mathcal{L}(X_i^O)(u)-\mathcal{L}_{K}(X_i^O)(u)=\mathcal{O}\left(\left(\sum_{k=K+1}^\infty \lambda_k^O\right)^{1/2}\right)=\mathcal{O}\left(K^{-(a_O-1)/2} \right)
\label{AME1}
\end{align}
\smallskip
Furthermore, for all $u\in M:=[a,b]\setminus O$
\begin{align}
&\widehat{\mathcal{L}}^*_{K}(\mathbb{X}_i^O)(u)=\mathcal{L}^*_{K}(X_i^O)(u)+
\mathcal{O}_p\left(m^{-2/5}+K\left(\frac{1}{m^{1/2}}+\frac{K^{a_O/2+3/2}}{\min\{n^{1/2},(n\mathcal{M})^{1/3}\}}\right)\right)
\nonumber \\
&\mathcal{L}(X_i^O)(u)-\mathcal{L}^*_{K}(X_i^O)(u)=\mathcal{O}\left(\left(\sum_{k=K+1}^\infty \lambda_k^O\right)^{1/2}\right)=\mathcal{O}\left(K^{-(a_O-1)/2} \right)
\label{AME2}
\end{align}
\end{theorem}
The theorem tells us that for any $u\in [a,b]$ the estimator $\widehat{\mathcal{L}}_{K}(\mathbb{X}_i^O)(u)$ achieves the same rate of convergence. But recall that for $u\in O=O_i$ we have $\mathcal{L}(X_i^O)(u)=X_i(u)$, and thus $\widehat{\mathcal{L}}_{K}(\mathbb{X}_i^O)(u)$ can be seen as a nonparametric estimator of $X_i$. In contrast, for $u\in M$ we have  $\mathcal{L}(X_i^O)(u)=X_i(u)+\mathcal{Z}_i(u)$, and therefore the distance between $X_i(u)$ and $\widehat{\mathcal{L}}_{K}(\mathbb{X}_i^O)(u)$ will additionally depend on the reconstruction error $\mathcal{Z}_i(u)$.

Note that by the above result the rates of convergence depend on $m$ and $n$, and the optimal $K$ depends on these quantities in a complex way. However, the situation simplifies if $m$ is considerably smaller than $n$ such that $m=m_n\asymp n^\theta$ for $\theta\leq 1/2$. The following corollary then is a direct consequence of \eqref{AME1}.

\begin{corollary}\label{cor:Estim_L}
Under the  conditions of Theorem \ref{Theorem_Main_Estim} additionally assume that $\theta\leq 1/2$. With
$K\equiv K_{m} \asymp m^{1/(a_O+2)}$ we obtain for all $u\in[a,b]$
\begin{align}
|\mathcal{L}(X_i^O)(u))-\widehat{\mathcal{L}}_{K}(\mathbb{X}_i^O)(u)|=
\mathcal{O}_p\left(m^{-\frac{a_O-1}{2(a_O+2)}} \right).
\label{AME3}
\end{align}
\end{corollary}

Let us consider the simple case where $m_i=m$ for all $i=1,\dots,n$, and recall that the main difference between $\widehat{\mathcal{L}}_{K}$  and $\widehat{\mathcal{L}}^*_{K}$ consists in the way of estimating $X_i$ on the observed interval $O:=O_i$. $\widehat{\mathcal{L}}^*_{K}$ is based on local linear smoothing of the individual data $(Y_{ij},U_{ij})$, $j=1,\dots,m$, and the associated  estimation error of order $m^{-2/5}$ appears in result \eqref{AME2}. Twice continuously differentiable functions are assumed, and using only individual data it is well-known that $m^{-2/5}$ constitutes the optimal rate of convergence of nonparametric function estimators with respect to this smoothness class.

In contrast, $\widehat{\mathcal{L}}_{K}(X_i^O)(u))$ combines information from all $n$ sample curves in order to estimate $X_i(u)$ for $u\in O$. If all samples curves are structurally similar in the sense that A6 holds for a very large $a_O\gg 1$, then \eqref{AME3} implies that the rate of convergence of  $\widehat{\mathcal{L}}_{K}(\mathbb{X}_i^O)(u)$ is very close to the parametric rate $m^{-1/2}$. That is, under the conditions of Corollary \ref{cor:Estim_L} ($m$  smaller than $\sqrt{n}$ and $a_O\gg 1$) it becomes advantageous to use $\widehat{\mathcal{L}}_{K}(\mathbb{X}_i^O)$ instead of $\widehat{\mathcal{L}}^*_{K}(\mathbb{X}_i^O)$ for estimating $X_i$ on the observed interval, since $\widehat{\mathcal{L}}_{K}(\mathbb{X}_i^O)$ may provide faster rates of convergence than the rate $m^{-2/5}$ achieved by nonparametric  smoothing of individual data.. We believe that this is an interesting result in its own right, which to our knowledge has not yet been established in the literature.

\section{Iterative reconstruction algorithm}\label{sec:pred_algo}
So far we have focused on the regular situation where the covariance function $\gamma(u,v)$ is estimable for all points $(u,v)\in[a,b]^2$. Under this situation we can reconstruct the entire missing parts of the functions, such that the reconstructed functions $\tilde{X}_i$ with
\begin{align}\label{eq:RCF}
\tilde{X}_i(u)=
\left\{
\begin{array}{ll}
\mathcal{L}(X_i^O)(u)&\text{if}\quad u\in M\\
X_i^O(u)&\text{if}\quad u\in O
\end{array}
\right.
\end{align}
are identifiable for \emph{all} $u\in[a,b]$.

In our application, however, we face the more restrictive situation where the mean function $\mu(u)$ can still be estimated for all $u\in[a,b]$, but where there is no information on $\gamma(u,v)$ for large values $|u-v|$; see Figure \ref{fig:emp_res}. This makes it impossible to reconstruct the entire missing part of a function, such that $\tilde{X}_i(u)$ cannot be identified for all $u\in[a,b]$.

In order to reconstruct functions $\tilde{X}_i$ that cover the total interval $[a,b]$, or at least a very large part of it, we propose successively plugging in the optimal reconstructions computed for subintervals. In the following we describe our iterative reconstruction algorithm:
\begin{algorithm}[Iterative reconstruction algorithm]\label{algo}\
\begin{description}
\item[$1$st Step] Denote the originally observed interval $O$ as $O_1$ and compute
\begin{align*}
\tilde{X}_{i,1}(u)=\left\{
\begin{array}{ll}
\mathcal{L}(X_i^{O_1})(u)&\text{if}\quad u\in M_1\\
X_i^{O_1}(u)               &\text{if}\quad u\in O_1
\end{array}
\right.
\end{align*}
\item[$r$th Step ($r\geq 2$)] Choose a new ``observed'' interval $O_r\subset O_{r-1}\cup M_{r-1}$ and use $\tilde{X}_i^{O_r}(u):=\tilde{X}_{i,r-1}(u)$ with $u\in O_r$ as the new ``observed'' fragment. Compute
\begin{align*}
\tilde{X}_{i,r}(u)=\left\{
\begin{array}{ll}
\mathcal{L}(\tilde{X}_i^{O_{r}})(u)&\text{if}\quad u\in M_r\\
\tilde{X}_i^{O_{r}}(u)               &\text{if}\quad u\in O_r.
\end{array}
\right.
\end{align*}
Join the reconstructed fragments $\tilde{X}_{i,1},\dots,\tilde{X}_{i,r}$ to form the new ``observed'' fragment $\tilde{X}_{i,r-1}$ on $O_{r-1}\cup M_{r-1}$ and repeat the $r$th step.
\item[Stopping] Stop if $\bigcup_{l=1}^{r}O_l\cup M_l=[a,b]$ or if $r=r_{\max}$.
\end{description}
\end{algorithm}

\begin{figure}[htbp]
  \centering
  \textrm{1st Step of the Reconstruction Algorithm}\par
  \begin{minipage}{0.48\textwidth}
    \centering
    \includegraphics[width=1.05\textwidth]{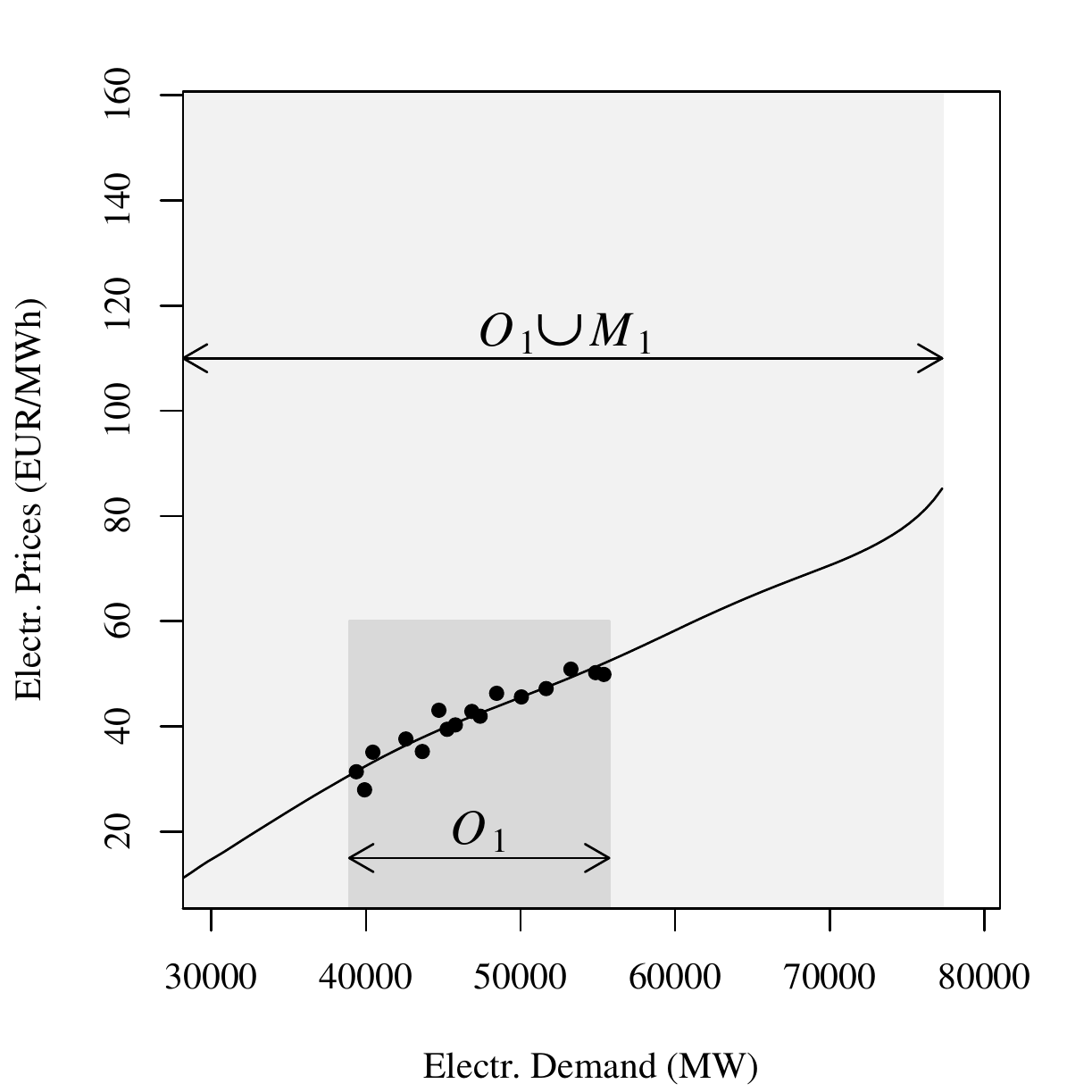}
  \end{minipage}
  \hspace*{2ex}
  \begin{minipage}{0.48\textwidth}
    \centering
    \includegraphics[width=1.05\textwidth]{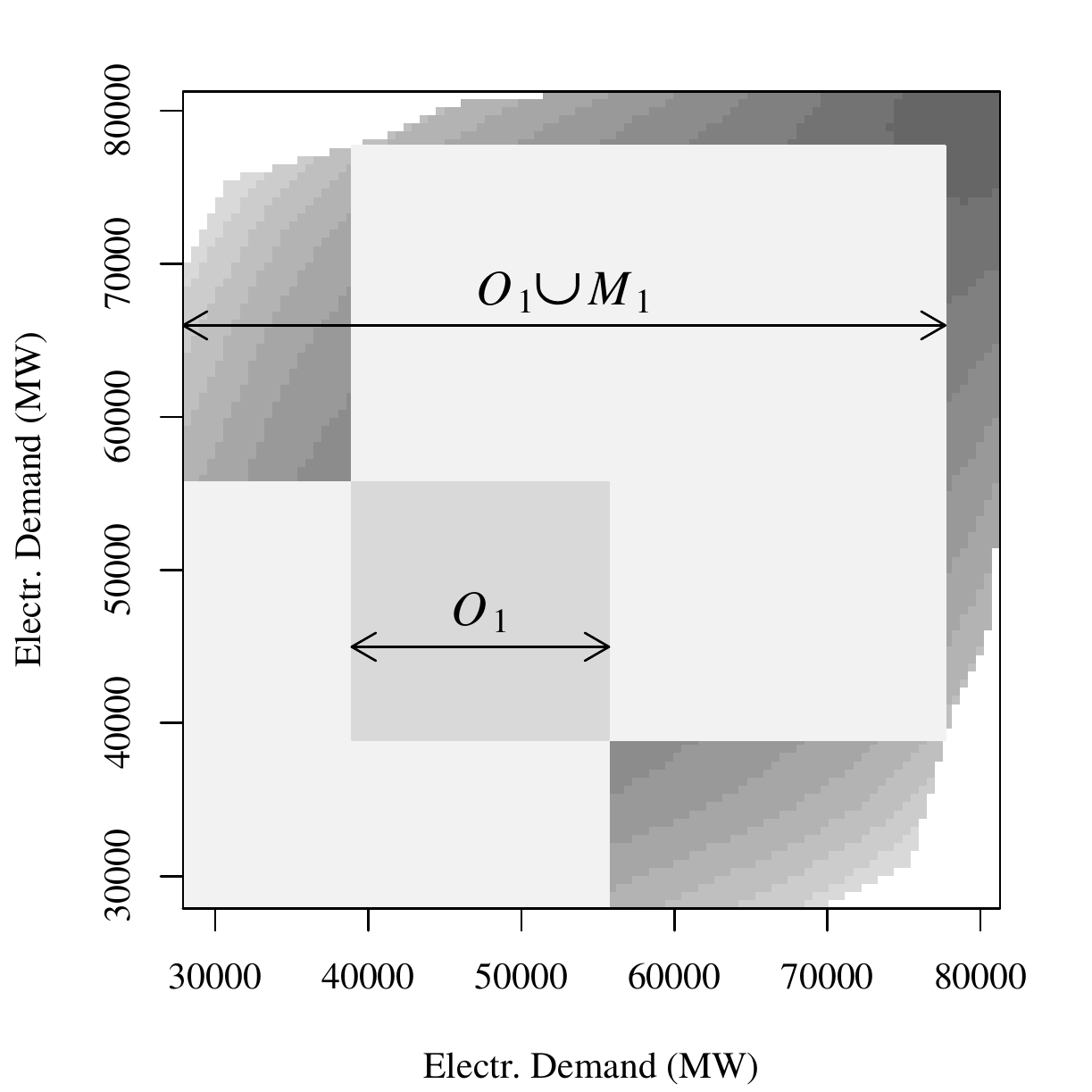}
  \end{minipage}
  \caption[]{Explanatory plots for the first run of the reconstruction algorithm.}
  \label{fig:pred_algo}
\end{figure}

This algorithm has to be applied to every fragment $X_i^O$. An exemplary first step of the reconstruction algorithm is shown in Figure \ref{fig:pred_algo}. The subinterval $O_1\cup M_1$ is determined by the original interval $O_1$ and the extend to which $\gamma$ can be estimated (see right panel). The function $\tilde{X}_{i,1}$ shown in the left panel still lacks the upper fragment for values $u\in[77362\text{ (MW)}, 82282\text{(MW)}]$ such that a second step of the reconstruction algorithm is necessary.

\begin{figure}[htbp]
  \centering
  \textrm{2nd Step of the Reconstruction Algorithm}\par
  \begin{minipage}{0.48\textwidth}
    \centering
    \includegraphics[width=1.05\textwidth]{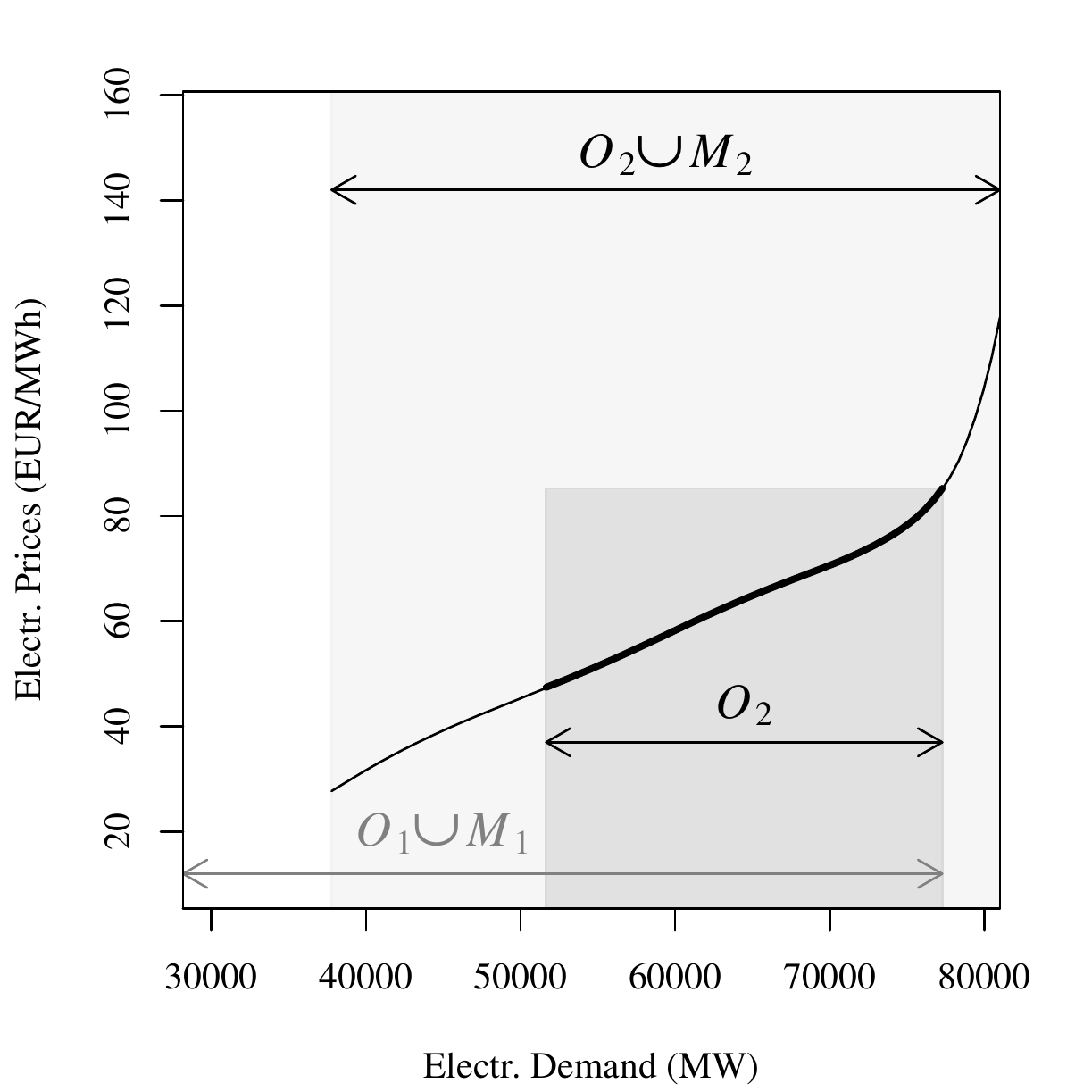}
  \end{minipage}
  \hspace*{2ex}
  \begin{minipage}{0.48\textwidth}
    \centering
    \includegraphics[width=1.05\textwidth]{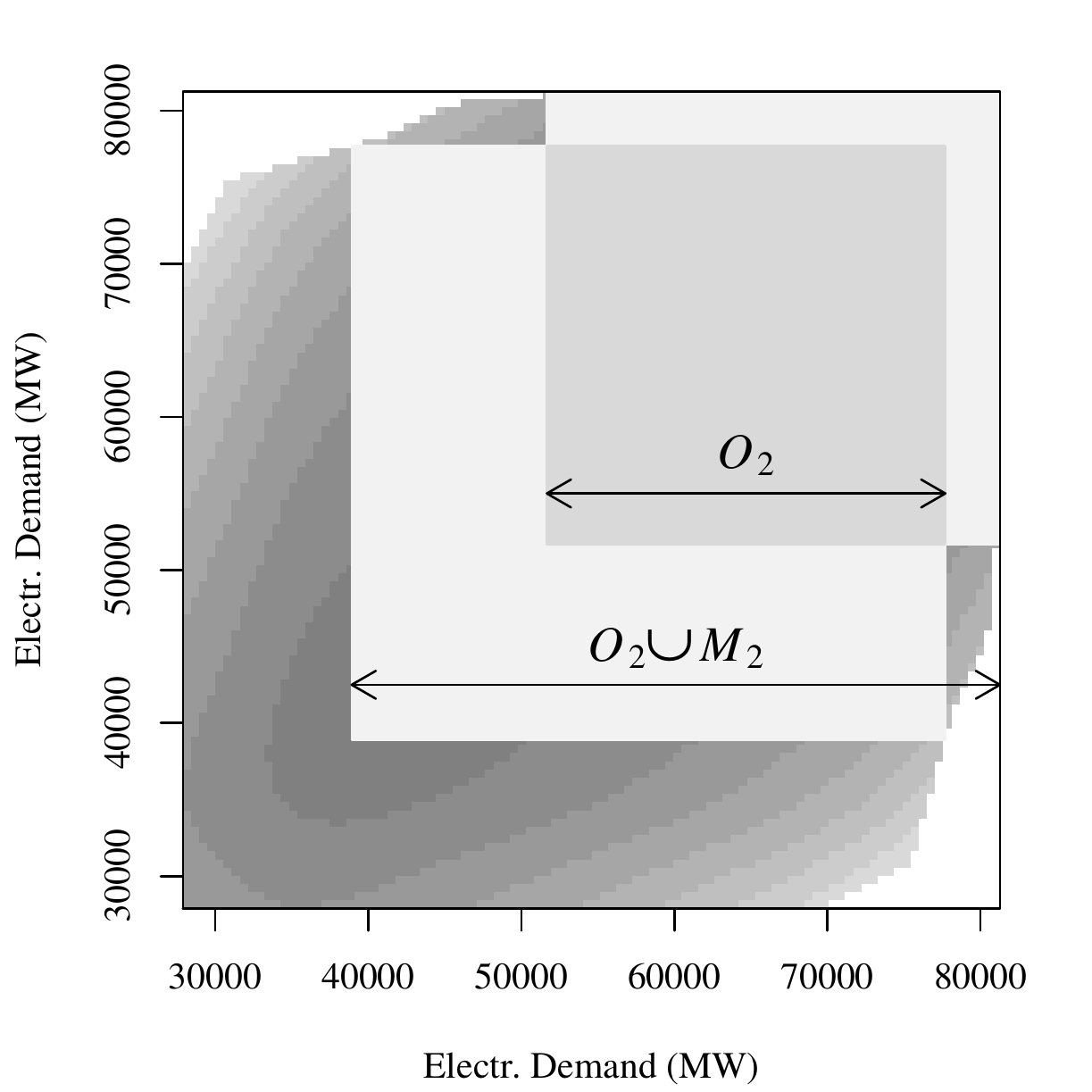}
  \end{minipage}
  \caption[]{Explanatory plots for the second run of the reconstruction algorithm.}
  \label{fig:pred_algo_2}
\end{figure}

This second step is shown in Figure \ref{fig:pred_algo_2}. There the new interval $O_2\subseteq O_1\cup M_1$ is chosen such that the still missing upper fragment becomes reconstructible. The new large interval $O_2\cup M_2$ contains the missing upper fragments, such that we can stop the  algorithm.


The choice of the subset $O_r$ in the $r$th step is crucial. On the one hand, $O_r$ should be chosen as large as possible to contain as much information as possible. On the other hand, $O_r$ must be chosen such that $M_r$ contains a still missing fragment which is -- in tendency -- met by smaller intervals $O_r$. That is, any efficient implementation of the algorithm and the choice of $r_{\max}$ depends on the extend to which $\gamma$ can be estimated. A simple practical implementation is described in our application in Section \ref{sec:appl}.

In each iteration of the reconstruction algorithm we accumulate reconstruction errors. The following proposition provides a theoretical description of this accumulation of reconstruction errors:
\begin{proposition}[Accumulated reconstruction error]\label{pro:IPA}
For simplicity, let $\E(X_i(u))=0$ for all $u\in[a,b]$ and consider the second step of the reconstruction algorithm.
Let $X_i^{M_2}(u)$ denote a missing value that we aim to reconstruct by $\mathcal{L}(\tilde{X}_i^{O_2})(u)$ using $\tilde{X}_i^{O_2}$ which is taken from the reconstruction of the $1$st Step. The mean squared reconstruction error can then be approximated as following:
\begin{align*}
\E\left(\left(X_i^{M_2}(u)-\mathcal{L}(\tilde{X}_i^{O_2})(u)\right)^2\right)\leq&
\E\left(\left(X_i^{M_2}(u)-\mathcal{L}(X_i^{O_2})(u)\right)^2\right)\\
+&\E\left(\left(X_i^{M_2}(u)-\mathcal{L}(X_i^{O_1})(u)\right)^2\right),
\end{align*}
where $\mathcal{L}(X_i^{O_1})$ and $\mathcal{L}(X_i^{O_2})$ are the hypothetical reconstruction operators if $\gamma$ were fully observed over $[a,b]^2$, and $X_i^{O_2}$ were observable.
\end{proposition}

That is, the mean squared reconstruction error in the second run of the iterative algorithm is bounded from above by the two hypothetical mean squared reconstruction errors of $\mathcal{L}(X_i^{O_1})(u)$ and $\mathcal{L}(X_i^{O_2})(u)$.


\section{Simulation study}\label{sec:sim}\
We compare the finite sample performance of our reconstruction operators \eqref{Estim_L} and \eqref{Estim_D} with that of the PACE method proposed by \cite{Yao2005} and the functional linear ridge regression model proposed by \cite{kraus2015}. A further interesting comparison method might be the functional linear regression model for sparse functional data as considered by \cite{yao2005functional}. Note, however, that this regression model becomes equivalent to the PACE method of \cite{Yao2005}, when used to predict the trajectory of $X_i$ given its own sparse, i.e., irregular and noise contaminated measurements (see Appendix \ref{appendix:FE2} in the supplementary paper \cite{KL_Suppl_19} for more detailed explanations regarding this equivalence). 

The following acronyms are used to refer to the different reconstruction methods considered in this simulation study:
\begin{description}
\item[ANo] $\widehat{\mathcal{L}}_{\hat{K}_{iM}}(\mathbb{X}_i^O)$ in \eqref{Estim_L} is denoted as ANo to indicate that this method involves No Alignment of the estimate of $X_i^O$ and the reconstruction of $X_i^M$.
\item[ANoCE] Equivalent to ANo, but with replacing the integral scores \eqref{PC-Score-Estim} using the following Conditional Exactions (CE) scores adapted from \cite{Yao2005}
\begin{align}\label{CEscores_small}
\hat{\xi}_{ik,\operatorname{CE}}^O=\hat{\lambda}_k^O\hat{\boldsymbol{\phi}}_{ik}^{O\top}\widehat{\boldsymbol{\Sigma}}^{-1}_{\mathbf{Y}_i}(\mathbf{Y}_i-\boldsymbol{\mu}_i),
\end{align}
where
$\mathbf{Y}_i=(Y_{i1},\dots,Y_{im_i})^\top$,
$\hat{\boldsymbol{\phi}}^O_{ik}=(\hat{\phi}_k^O(U_{i1}),\dots,\hat{\phi}_k^O(U_{im_i}))^\top$, \\
$[\widehat{\boldsymbol{\Sigma}}_{\mathbf{Y}_i}]_{1\leq j,k\leq m_i}=\hat{\gamma}(U_{ij},U_{ik})+\hat{\sigma}^2\delta_{jk}$, with $\delta_{jk}=1$ if $j=k$ and zero else, and with $\hat{\lambda}_k^O$ and $\hat{\phi}_k^O$ as defined in \eqref{small_eigen}. The estimate of the error variance, $\hat{\sigma}^2$, is computed using LLK estimators as described in equation (2) of \cite{Yao2005}.
\item[AYes] $\widehat{\mathcal{L}}_{\hat{K}_{iM}}^*(\mathbb{X}_i^O)$ in \eqref{Estim_D} is denoted as AYes to indicate that this method involves an alignment of the estimate of $X_i^O$ and the reconstruction of $X_i^M$.
\item[AYesCE] Equivalent to AYes, but with replacing the integral scores \eqref{PC-Score-Estim} by the conditional exaction scores of \eqref{CEscores_small}.
\item[PACE] The method of \cite{Yao2005}, who approximate the missing $X_i^{M}$ and observed $X_i^O$ parts jointly using the truncated Karhunen-Lo\`eve decomposition $\widehat{X}_i(t)=\hat{\mu}(t)+\sum_{k=1}^{\hat{K}_{iM}}\hat{\xi}_{ik}^{\operatorname{PACE}}\hat{\phi}_k(t)$ with conditional expectation scores
\begin{align}\label{CEscores_large}
\hat{\xi}_{ik}^{\operatorname{PACE}}=\hat{\lambda}_k\hat{\boldsymbol{\phi}}_{ik}^{\top}\widehat{\boldsymbol{\Sigma}}^{-1}_{\mathbf{Y}_i}(\mathbf{Y}_i-\boldsymbol{\mu}_i),
\end{align}
where $\hat{\lambda}_k$ and $\hat{\phi}_k$ are as defined in \eqref{small_eigen}, but with $O=[a,b]$.
\item[KRAUS] The functional linear ridge regression model of \cite{kraus2015}.
\end{description}

The idea of using the conditional expectation scores \eqref{CEscores_small} in ANoCE and AYesCE as an alternative to the integral scores \eqref{PC-Score-Estim} in ANo and AYes is inspired by a comment of one of the anonymous referees, who correctly pointed out that the integral scores \eqref{PC-Score-Estim} might be instable for irregular and noisy data. PACE also uses condition expectation scores, but is fundamentally different from ANoCE and AYesCE. While PACE uses approximations of the classical eigenfunctions $\phi_k$, the classical eigenvalues $\lambda_k$, and the classical scores $\xi_{ik}$, ANoCE and AYesCE use approximations of the reconstructive eigenfunctions $\tilde{\phi}_k^O$, the eigenvalues $\lambda_k^O$, and the scores $\xi^O_{ik}$ with respect to the partially observed domain $O$.

The truncation parameters $\hat{K}_{iM}$ for ANo, ANoCE, AYes, AYesCE, and PACE are selected by minimizing the GCV criterion in \eqref{GCV}. For PACE, we do not use the AIC-type criterion as proposed by \cite{Yao2005}, since this criterion determines a ``global'' truncation parameter $\hat{K}$, which performs worse than our local, i.e., $M$-specific truncation parameter $\hat{K}_{iM}$. The ridge regularization parameter for KRAUS is determined using the GCV criterion as described in \cite{kraus2015}.

We consider four different Data Generating Processes (DGPs). DGP1 and DGP2 comprise irregular evaluation points and measurement errors which facilitates the comparison of ANo, ANoCE, AYes, AYesCE and the PACE method. DGP3 and DGP4 comprise regular evaluation points and no measurements errors which facilitates the comparison of ANo, AYes, PACE and the KRAUS method. For all simulations we set $[a,b]=[0,1]$.
\smallskip\newline\noindent\textbf{DGP1} 
The data points
$(Y_{ij},U_{ij})$ are generated according to
$Y_{ij}=X_i(U_{ij})+\varepsilon_{ij}$
with error term
$\varepsilon_{ij}\sim N(0,0.0125)$
and random function
$X_i(u)=\mu(u)+\sum_{k=1}^{50}\xi_{ik,1}\cos(k\pi u)/\sqrt{5}+\xi_{ik,2}\sin(k\pi u)/\sqrt{5}$,
where
$\mu(u)=u + \sin(2\pi u)$,
$\xi_{ik,1}=50\sqrt{\exp(-(k-1)^2/5)} Z_{i,1}$, and
$\xi_{ik,2}=50\sqrt{\exp(- k^2   /5)} Z_{i,2}$ with
$Z_{i,1},Z_{i,2}\sim N(0,1)$.
The evaluation points are generated as $U_{ij}\sim \operatorname{Unif}[A_i,B_i]$, where with probability $1/2$, $A_i\sim\operatorname{Unif}[0,0.45]$ and $B_i\sim\operatorname{Unif}[0.55,1]$ and with probability $1/2$, $[A_i,B_i]=[0,1]$. That is, about one half of the sample consists of partially observed functions with mean interval-width $0.55$.
\smallskip\newline\noindent\textbf{DGP2} 
Equivalent to DGP1, except for a larger noise component with $\varepsilon_{ij}\sim N(0,0.125)$.
\smallskip\newline\noindent\textbf{DGP3} 
The data points
$(Y_{ij},U_{ij})$ are generated according to
$Y_{ij}=X_i(U_{ij})$
with random function
$X_i(u)=\mu(u)+\sum_{k=1}^{50}\xi_{ik,1}\cos(k\pi u)+\xi_{ik,2}\sin(k\pi u)$,
where
$\mu(u)=u^2 + \sin(2\pi u)$,
$\xi_{ik,1}=50\sqrt{\exp(-(k-1)^2)}  Z_{i,1}$, and
$\xi_{ik,2}=50\sqrt{\exp(- k^2   )}  Z_{i,2}$ with
$Z_{i,1},Z_{i,2}\sim N(0,1)$.
The evaluation points are equidistant grid points $U_{ij}=j/51$, with $j=1,\dots,51$, where all points $U_{ij}\not\in[A_i,B_i]$ are set to \textsf{NA}. With probability $3/4$, $A_i\sim\operatorname{Unif}[0,1/3]$ and $B_i=A_i+1/2$ and with probability $1/4$, $[A_i,B_i]=[0,1]$.
\smallskip\newline\noindent\textbf{DGP4} 
Equivalent to DGP3, but with $A_i\sim\operatorname{Unif}[0,2/3]$ and $B_i=A_i+1/3$. That is, DGP4 has smaller and therefore more challenging fragments than DGP3.

\smallskip

For each DGP, we generate 50 different targets $X_{\ell}$, $\ell=1,\dots,50$, where each target is partitioned into a (non-degenerated) missing part $X_{\ell}^{M}$ and an observed part $X_{\ell}^O$. Each of these targets $X_{\ell}$ are reconstructed in each of the $b=1,\dots,100$ simulation runs with sample sizes $n\in\{50,100\}$ for DGP1-DGP4 and $m\in\{15,30\}$ for DGP1 and DGP2.

Let $\widehat{X}_{\ell,b}$ denote the reconstructed function in simulation run $b$ using one of the reconstruction methods ANo, ANoCE, AYes, AYesCE, PCAE, or KRAUS. For each target $X_{\ell}$, we compute the integrated mean squared error, the integrated squared bias, and the integrated variance,
\begin{align*}
\operatorname{MSE}_\ell&=\operatorname{Var}_\ell+\operatorname{Bias}^2_\ell,\quad
\operatorname{Bias}^2_\ell=\int_0^1\big(\bar{X}_{\ell}(u) - X_{\ell}(u)\big)^2dt,\\
\text{and}\quad\operatorname{Var}_\ell&=\int_0^1 100^{-1}\sum_{b=1}^{100}\big(\widehat{X}_{\ell,b}(u)-\bar{X}_{\ell}(u)\big)^2dt,
\end{align*}
where
$\bar{X}_{\ell}(u)=100^{-1}\sum_{r=1}^{100}\widehat{X}_{\ell,b}(u)$. The finite sample performance is evaluated using the averages over all 50 targets,
$$\operatorname{Var}=\frac{1}{50}\sum_{\ell=1}^{50}\operatorname{Var}_\ell,\quad\operatorname{Bias}^2=\frac{1}{50}\sum_{\ell=1}^{50}\operatorname{Bias}^2_\ell,\quad\text{and}\quad\operatorname{MSE}=\frac{1}{50}\sum_{\ell=1}^{50}\operatorname{MSE}_\ell.$$

\begin{table}[!htb]
\caption{Simulation results for DGP1.}
\label{TabDGP1}
\centering
\begin{tabular}{l cc l cccc}
\toprule
  DGP  & $n$ & $m$& Method & MSE$_{\text{ratio}}$ & MSE & $\text{Bias}^2$ & Var \\
\midrule
  DGP1 & 50 & 15 & AYesCE & 1.00 & 0.161 & 0.135 & 0.025 \\
  DGP1 & 50 & 15 & AYes & 1.02 & 0.164 & 0.139 & 0.025 \\
  DGP1 & 50 & 15 & ANoCE & 1.38 & 0.222 & 0.199 & 0.023 \\
  DGP1 & 50 & 15 & ANo & 1.39 & 0.224 & 0.200 & 0.024 \\
  DGP1 & 50 & 15 & PACE & 10.49 & 1.685 & 0.259 & 1.426 \\
\midrule
  DGP1 & 50 & 30 & AYesCE & 1.00 & 0.136 & 0.112 & 0.024 \\
  DGP1 & 50 & 30 & AYes & 1.00 & 0.137 & 0.113 & 0.024 \\
  DGP1 & 50 & 30 & ANoCE & 1.48 & 0.202 & 0.173 & 0.029 \\
  DGP1 & 50 & 30 & ANo & 1.53 & 0.209 & 0.180 & 0.029 \\
  DGP1 & 50 & 30 & PACE & 5.19 & 0.707 & 0.131 & 0.576 \\
\midrule
  DGP1 & 100 & 15 & AYesCE & 1.00 & 0.131 & 0.112 & 0.018 \\
  DGP1 & 100 & 15 & AYes & 1.00 & 0.131 & 0.114 & 0.017 \\
  DGP1 & 100 & 15 & ANoCE & 1.58 & 0.207 & 0.191 & 0.017 \\
  DGP1 & 100 & 15 & ANo & 1.61 & 0.211 & 0.194 & 0.017 \\
  DGP1 & 100 & 15 & PACE & 8.74 & 1.145 & 0.154 & 0.991 \\
\midrule
  DGP1 & 100 & 30 & AYes & 1.00 & 0.125 & 0.108 & 0.017 \\
  DGP1 & 100 & 30 & AYesCE & 1.01 & 0.126 & 0.109 & 0.017 \\
  DGP1 & 100 & 30 & ANoCE & 1.36 & 0.170 & 0.146 & 0.023 \\
  DGP1 & 100 & 30 & ANo & 1.45 & 0.181 & 0.158 & 0.023 \\
  DGP1 & 100 & 30 & PACE & 3.59 & 0.448 & 0.123 & 0.325 \\
\bottomrule
\multicolumn{8}{l}{MSE$_{\text{ratio}}=\text{MSE}/\min(\text{MSE})$}
\end{tabular}
\end{table}

\begin{table}[!htb]
\caption{Simulation results for DGP2.}
\label{TabDGP2}
\centering
\begin{tabular}{l cc l cccc}
\toprule
DGP  & $n$ & $m$& Method & MSE$_{\text{ratio}}$ & MSE & $\text{Bias}^2$ & Var \\
\midrule
  DGP2 & 50 & 15 & AYesCE & 1.00 & 0.198 & 0.173 & 0.025 \\
  DGP2 & 50 & 15 & AYes & 1.04 & 0.207 & 0.179 & 0.027 \\
  DGP2 & 50 & 15 & PACE & 1.07 & 0.212 & 0.174 & 0.039 \\
  DGP2 & 50 & 15 & ANoCE & 1.14 & 0.227 & 0.203 & 0.023 \\
  DGP2 & 50 & 15 & ANo & 1.16 & 0.230 & 0.204 & 0.026 \\
\midrule
  DGP2 & 50 & 30 & AYesCE & 1.00 & 0.189 & 0.167 & 0.022 \\
  DGP2 & 50 & 30 & AYes & 1.01 & 0.192 & 0.169 & 0.023 \\
  DGP2 & 50 & 30 & PACE & 1.09 & 0.206 & 0.167 & 0.039 \\
  DGP2 & 50 & 30 & ANoCE & 1.14 & 0.215 & 0.188 & 0.027 \\
  DGP2 & 50 & 30 & ANo & 1.16 & 0.219 & 0.190 & 0.028 \\
\midrule
  DGP2 & 100 & 15 & AYesCE & 1.00 & 0.178 & 0.161 & 0.017 \\
  DGP2 & 100 & 15 & AYes & 1.01 & 0.180 & 0.162 & 0.018 \\
  DGP2 & 100 & 15 & PACE & 1.08 & 0.193 & 0.165 & 0.028 \\
  DGP2 & 100 & 15 & ANoCE & 1.20 & 0.213 & 0.198 & 0.015 \\
  DGP2 & 100 & 15 & ANo & 1.21 & 0.216 & 0.199 & 0.018 \\
\midrule
 DGP2 & 100 & 30 & AYesCE & 1.00 & 0.177 & 0.159 & 0.018 \\
  DGP2 & 100 & 30 & AYes & 1.03 & 0.181 & 0.162 & 0.020 \\
  DGP2 & 100 & 30 & PACE & 1.03 & 0.183 & 0.153 & 0.029 \\
  DGP2 & 100 & 30 & ANoCE & 1.07 & 0.189 & 0.167 & 0.023 \\
  DGP2 & 100 & 30 & ANo & 1.12 & 0.197 & 0.174 & 0.024 \\
\bottomrule
\multicolumn{8}{l}{MSE$_{\text{ratio}}=\text{MSE}/\min(\text{MSE})$}
\end{tabular}
\end{table}

The simulation study is implemented using the \textsf{R}-package \texttt{ReconstPoFD} which can be downloaded and installed from the second author's GitHub account.

Table \ref{TabDGP1} shows the simulation results for DGP1. The methods (ANo, ANoCE, AYes, AYesCE and PACE) are ranked according to their MSE$_{\text{ratio}}$ which is defined by the method's MSE-value relative to the lowest MSE-value within the comparison group. The rankings are stable for all sample sizes $m$ and $n$. The AYesCE reconstruction method shows the best performance. The AYes method, which uses integral scores instead of conditional expectation scores, is only marginally less efficient than AYesCE. Our non-alignment methods ANoCE and ANo are ranked third and fourth. The PACE method of \cite{Yao2005}, originally proposed for sparse functional data analysis, shows a rather poor performance. The reason for this is that PACE adds the variance of the measurement error to the diagonal of the discretized covariance matrix, which has a regularization effect on the generally ill-posed inversion problem. For DGP1, however, the variance of the error term is rather small which results in a too small regularization of the inverse.\\
Table \ref{TabDGP2} shows the simulation results for DGP2. DGP2 is equivalent to DGP1 except for a larger variance of the error term. Our alignment methods AYesCE and AYes still show the best performance. However, having a larger variance leads to a better regularization of the inverse problem involved in the PACE method, such that PACE is ranked third. Our non-alignment methods ANoCE and ANo are ranked fourth and fifth. Figures \ref{fig:simresults_DGP1} and \ref{fig:simresults_DGP2} in Appendix \ref{visual} of the supplementary paper \cite{KL_Suppl_19} provide graphical illustrations of the different reconstruction results as well as a visual impression of the different signal-to-noise ratios in DGP1 and DGP2.

\begin{table}[ht]
\caption{Simulation results for DGP3 and DGP4.}
\label{TabDGP34}
\centering
\begin{tabular}{l c l cccc}
\toprule
DGP  & $n$ & Method& MSE$_{\text{ratio}}$ & MSE & $\text{Bias}^2$ & Var \\
\midrule
  DGP3 & 50 & AYes & 1.00 & 0.168 & 0.131 & 0.037 \\
  DGP3 & 50 & PACE & 1.33 & 0.223 & 0.099 & 0.124 \\
  DGP3 & 50 & ANo & 1.40 & 0.234 & 0.178 & 0.056 \\
  DGP3 & 50 & KRAUS & 1.52 & 0.254 & 0.205 & 0.049 \\
\midrule
  DGP3 & 100 & AYes & 1.00 & 0.142 & 0.120 & 0.022 \\
  DGP3 & 100 & PACE & 1.26 & 0.179 & 0.081 & 0.098 \\
  DGP3 & 100 & KRAUS & 1.29 & 0.184 & 0.151 & 0.033 \\
  DGP3 & 100 & ANo & 1.36 & 0.194 & 0.158 & 0.035 \\
\midrule
\midrule
  DGP4 & 50 & AYes & 1.00 & 0.276 & 0.220 & 0.056 \\
  DGP4 & 50 & ANo & 1.11 & 0.307 & 0.247 & 0.060 \\
  DGP4 & 50 & KRAUS & 1.20 & 0.330 & 0.269 & 0.061 \\
  DGP4 & 50 & PACE & 41.93 & 11.564 & 0.313 & 11.252 \\
\midrule
  DGP4 & 100 & AYes & 1.00 & 0.232 & 0.202 & 0.030 \\
  DGP4 & 100 & KRAUS & 1.11 & 0.258 & 0.222 & 0.035 \\
  DGP4 & 100 & ANo & 1.12 & 0.261 & 0.227 & 0.034 \\
  DGP4 & 100 & PACE & 3.59 & 0.834 & 0.151 & 0.682 \\
\bottomrule
\multicolumn{7}{l}{MSE$_{\text{ratio}}=\text{MSE}/\min(\text{MSE})$}
\end{tabular}
\end{table}
Table \ref{TabDGP34} shows the simulation results for DGP3 and DGP4 comparing the methods ANo, AYes, PACE and KRAUS. Here, the alignment method AYes shows by far the best performance for all sample sizes and for both DGPs. The partially very bad performance of PACE is due to the missing measurement error in DGP3 and DGP4, which results in a missing regularization of the inverse problem involved in the PACE method. Furthermore, PACE is designed for the case where one observes only a few noisy discretization points per function, but these points should be distributed over the total domain $[a,b]$. For the considered DGPs, however, the discretization points are only observed within challenging small subdomains $[A_i,B_i]\subset[0,1]$. Graphical illustrations of the different reconstruction results for DGP3 and DGP4 are provided in Figures \ref{fig:simresults_DGP3} and \ref{fig:simresults_DGP4} in Appendix \ref{visual} of the supplementary paper \cite{KL_Suppl_19}.

Summing up, in all DGPs the best performing reconstruction method are our alignment methods AYesCE and AYes. For discretized functional data plus measurement errors it is advantageous to use the alignment method AYesCE with involves conditional expectation scores.

\section{Application}\label{sec:appl}
Our functional data point of view on electricity spot prices provides a practical framework that is useful for forecasting electricity spot prices \citep{Liebl2013,weron2014} and for testing price differences \citep{Liebl2018}. In the following, we focus on the problem of reconstructing the partially observed price-functions, which is highly relevant for practitioners who need complete price functions for doing comparative statics, i.e., a ceteris-paribus analysis of price effects with respect to changes in electricity demand \citep[cf.][]{weigt2009germany,hirth2013market}.

The data for our analysis come from three different sources. Hourly spot prices of the German electricity market are provided by the European Energy Power Exchange (EPEX) (\url{www.epexspot.com}), hourly values of Germany's gross electricity demand, $D_{ij}$, and net-imports of electricity from other countries, $N_{ij}$, are provided by the European Network of Transmission System Operators for Electricity (\url{www.entsoe.eu}), and German wind and solar power infeed data are provided by the transparency platform of the European energy exchange (\url{www.eex-transparency.com}). The data dimensions are given by $m=24$ hours and $n=241$ working days between March 15, 2012 and March 14, 2013. Very few ($0.4\%$) data pairs $(Y_{ij},U_{ij})$ with prices $Y_{ij}>120$ EUR/MWh and $U_{ij}>82000$ MW are considered as outliers and reset to $Y_{ij}=120$.
The German electricity market, like many other electricity markets, provides purchase guarantees for Renewable Energy Sources (RES). Therefore, the relevant variable for pricing at the energy exchange is electricity demand minus electricity infeeds from RES \citep{nicolosi2010wind}. Correspondingly, the relevant values of electricity demand $U_{ij}$ are defined as electricity demand minus infeeds from RES and plus net-imports from other countries, i.e., $U_{ij}:=D_{ij}-\texttt{RES}_{ij}+N_{ij}$, where $\texttt{RES}_{ij}=\texttt{Wind.Infeed}_{ij}+\texttt{Solar.Infeed}_{ij}$. The effect of further RES such as biomass is still negligible for the German electricity market.
\begin{figure}[!t]
\begin{minipage}{0.48\textwidth}
    \centering
    \includegraphics[width=1.1\textwidth]{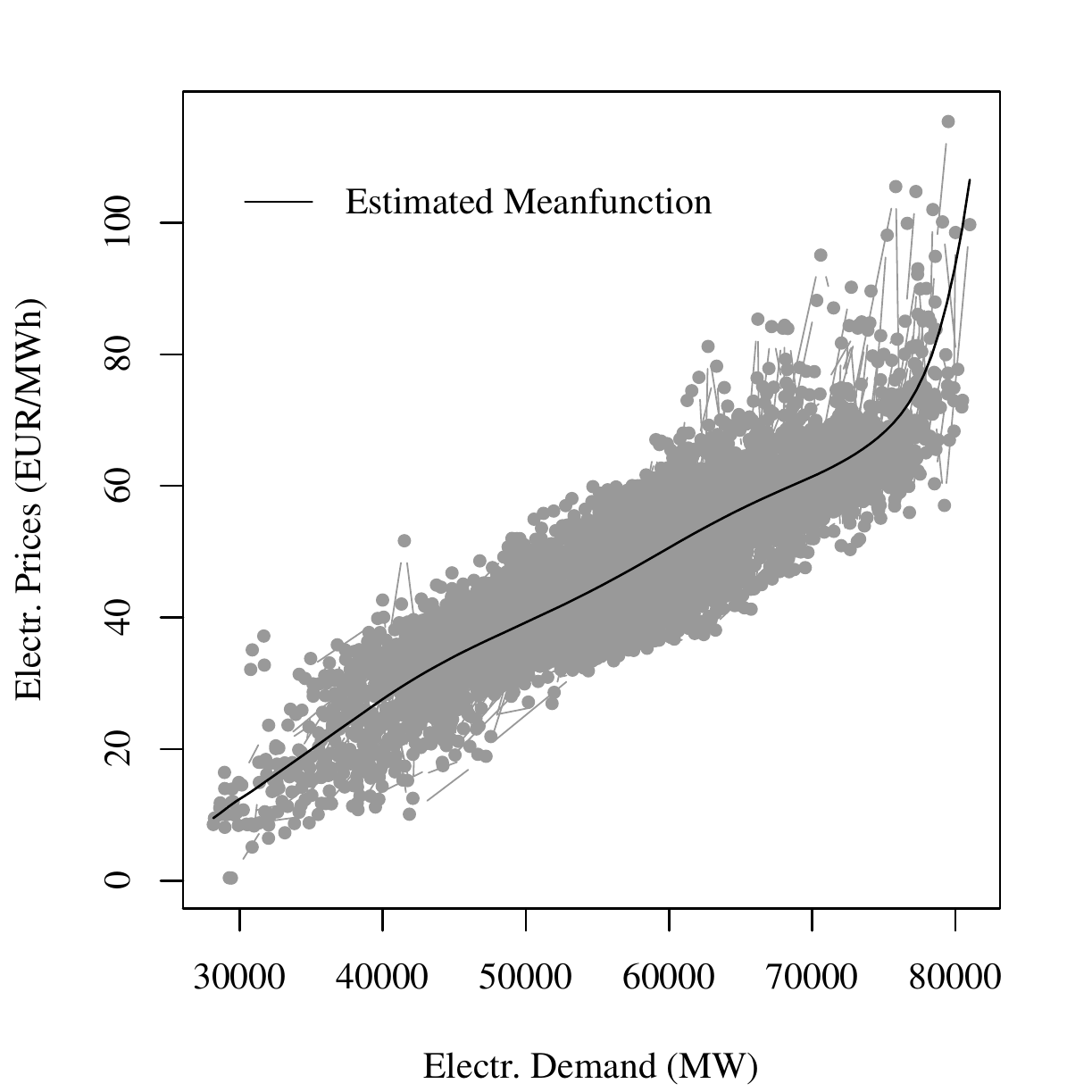}
\end{minipage}
\hspace*{1ex}
\begin{minipage}{0.48\textwidth}
    \centering
\includegraphics[width=1.1\textwidth]{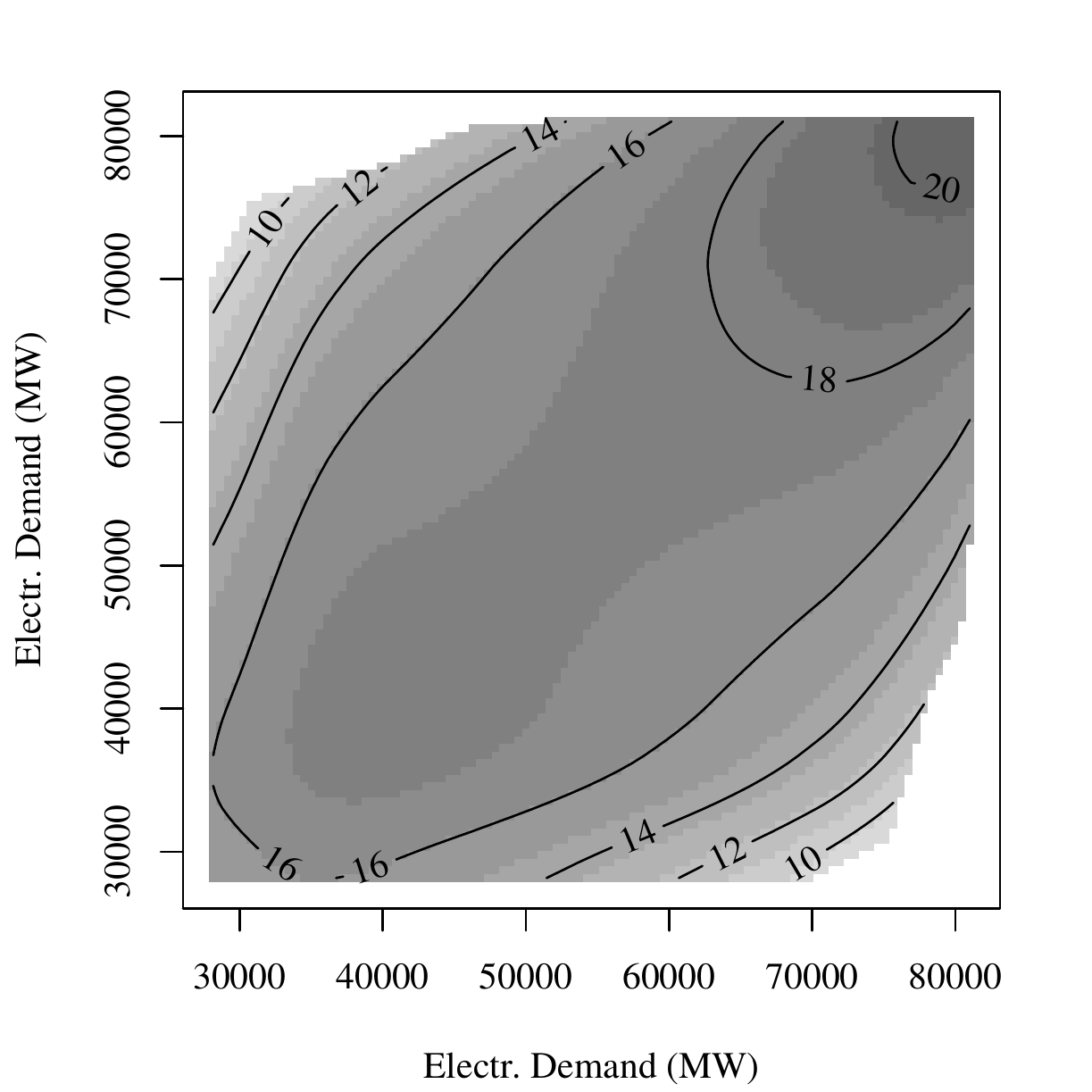}
\end{minipage}
\caption[]{{\sc Left Panel:} Estimated mean function plus a scatter plot of the data pairs $(Y_{ij},U_{ij})$. {\sc Right Panel:} Contour plot of the estimated covariance function. The white regions reflect the outer off-diagonal parts which are infeasible to estimate.}
\label{fig:emp_res}
\end{figure}

The estimated mean and covariance functions are shown in Figure \ref{fig:emp_res}. The outer off-diagonal parts of the covariance function $\gamma$ cannot be estimated, since these parts of the domain are not covered by data pairs $(U_{ij},U_{il})$, $j\neq l$. In order to reconstruct the entire missing parts $X_i^M$, we use the AYesCE estimator, which showed a very good performance in our simulation studies, and our iterative reconstruction Algorithm \ref{algo} implemented as follows. We use three iterations for each partially observed price function. In the first step, we use the information with respect to the original observations $\mathbb{X}_i^O$ in order to reconstruct the missing parts as far as possible. In the second step, we use the upper half of the reconstructed curve $\hat{\tilde{X}}_{i,1}$ and try to reconstruct possibly further missing upper fragments. In the final step we use the lower half of $\hat{\tilde{X}}_{i,1}$ and try to reconstruct possibly further missing lower fragments.

\begin{figure}[!t]
\centering
\includegraphics[width=.55\textwidth]{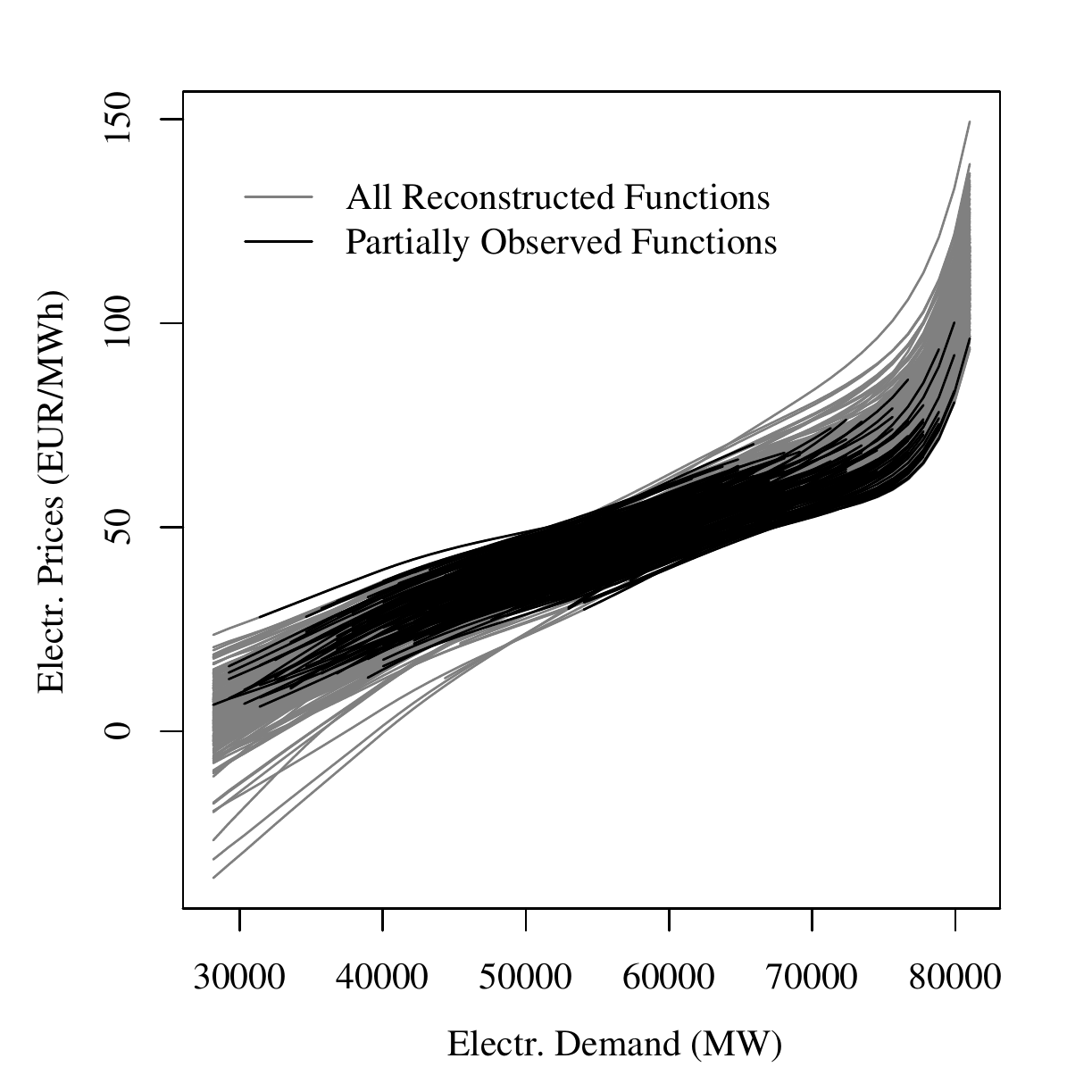}
\caption[]{Recovered functions (gray) and the original partially observed functions (black).}
\label{fig:emp_res_2}
\end{figure}

This approach allows us to recover 91\% of the price functions over the total support (Figure \ref{fig:emp_res_2}). Note that the price functions with negative electricity prices are perfectly plausible. Negative prices are an important market-feature of the EPEX \citep[see, for instance,][]{nicolosi2010wind,FGP2013,CHMG2014}. Electricity producers are willing to sell electricity at negative prices (i.e., to pay for delivering electricity) if shutting off and restarting their power plants is more expensive than selling their electricity at negative prices. That is, the reconstructed price functions are conform with the specific market design of the EPEX and may be useful for a variety of further subsequent analysis using classical methods of functional data analysis.

\section*{Acknowledgements}
We would like to thank the referees and the editors for their constructive feedback which helped to improve this research work. 

\begin{supplement}
\stitle{Supplemental Paper}
\slink[doi]{COMPLETED BY THE TYPESETTER}
\sdescription{The supplemental paper contains the proofs of our theoretical results.}
\end{supplement}

\bibliographystyle{imsart-nameyear}
\bibliography{bibfile}

\hfilneg\printaddressnum{1}

\newpage

\setcounter{page}{1}
\thispagestyle{plain}
\begin{center}
{\sc \Large Supplementary paper for:\\[2ex]
    {On the Optimal Reconstruction of Partially Observed Functional Data}}\\[2ex]
  by Alois Kneip and Dominik Liebl
\end{center}

\appendix

\pagenumbering{Roman}

\section*{Content}
In the following we give the proofs of our theoretical results. The main steps in our proofs of Theorems \ref{Theorem_UR} and \ref{Theorem_Main_Estim} are as in \cite{Yao2005}. Though, by contrast to \cite{Yao2005}, we allow for a time series context (see Assumption A1), impose more restrictive assumptions on the kernel function (see Assumption A5), and consider a different asymptotic setup (see Assumption A2). Appendix \ref{appendix:FE} contains further explanations and Appendix \ref{visual} contains visualizations of our simulation results.

\setcounter{section}{0}
\section{Proofs}\label{Proofs}\

\noindent\textbf{Proof of Theorem \ref{th:SRO}:}
For every linear operator $L:\mathbb{L}^2(O)\to\mathbb{L}^2(M)$ that is a reconstruction operator with respect to $X_i^O$ according to Def.~\ref{def:SRO}, we have that
\begin{align}\label{eq:VL_1}
\V(L(X_i^O)(u))=\sum_{k=1}^\infty \lambda_k^O \big(L(\phi_k^O)(u)\big)^2\quad\text{for every}\quad u\in M.
\end{align}
Existence: Writing $L(X_i^O)(u)$ as $L(X_i^O)(u)=\langle\alpha_u,X_i^O\rangle_H$ for some $\alpha_u\in H$ and computing again the variance of $L(X_i^O)(u)$ yields that
\begin{align}\label{eq:VL_2}
\V(L(X_i^O)(u))=\sum_{k=1}^\infty\lambda_k^O\left(\frac{\langle\alpha_u,\phi_k^O\rangle_2}{\lambda_k^O}\right)^2.
\end{align}
Since \eqref{eq:VL_1} and \eqref{eq:VL_2} must be equal, we have that $L(\phi_k^O)(u)=\langle\alpha_u,\phi_k^O\rangle_2/\lambda_k^O$ for all $k\geq 1$, which establishes that there exits a $\alpha_u\in H$ for every reconstruction $L(X_i^O)(u)$.\\
Uniqueness: Assume that there is an alternative $\tilde{\alpha}_u\in H$ such that $L(\phi_k^O)(u)=\langle\tilde{\alpha}_u,\phi_k^O\rangle_2/\lambda_k^O$ for all $k\geq 1$. Then
$\langle\alpha_u-\tilde{\alpha}_u,\phi_k^O\rangle_2/\lambda_k^O =0$ for all $k\geq 1$ or equivalently
$\langle\alpha_u,\phi_k^O\rangle_2-\langle\tilde{\alpha}_u,\phi_k^O\rangle_2=0$ for all $k\geq 1$ which shows that $\tilde{\alpha}_u-\alpha_u=0$.

\bigskip

\noindent\textbf{Proof of Theorem \ref{WellDef}, part (a):}
First, note that continuity of $\gamma(u,v)$ implies continuity of $\V(\mathcal{L}(X_i^O)(u))$. Second, note that for any $K$ and every $u\in M$, we have
\begin{align}\label{MSPE}
  0\leq \E\left(\Big(X_i^M(u)-\sum_{k=1}^K\xi^{O}_{ik}\tilde{\phi}^O_{k}(u)\Big)^2\right) =
\gamma(u,u)-\sum_{k=1}^K\lambda_k^{O}\tilde{\phi}_{k}^O(u)^2.
\end{align}
But this implies that $\V(\mathcal{L}_{u,K}(X_i^O))=\V(\sum_{k=1}^K\xi^{O}_{ik}\tilde{\phi}^O_{k}(u))=\sum_{k=1}^K\lambda_k^{O}\tilde{\phi}^O_{k}(u)^2$ converges to a fixed limit $0\leq \V(\mathcal{L}(X_i^O)(u))<\infty$ as $K\to\infty$ for all $u\in M$.

\smallskip

\noindent\textbf{Part (b):} Follows directly from observing that $\E(\mathcal{L}(X_i^O)(u))=0$ for all $u\in M$.

\smallskip

\noindent\textbf{Proof of Theorem \ref{OptimalPrediction}, part (a):}
For all $ v\in O \text{ and } u\in M$ we have that
\begin{align*}
 &\E\left(X_i^{O}(v)\mathcal{Z}_i(u)\right)=
  \E\left(X_i^{O}(v)\left(X_i^M(u)-\mathcal{L}(X_i^O)(u)\right)\right)=\\
 &=\E\left(\sum_{k=1}^\infty\xi^{O}_{ik}\phi^{O}_{k}(v)
 \Big(X_i^M(u)-\sum_{k=1}^\infty\xi^{O}_{ik}\tilde{\phi}_{k}^O(u)\Big)\right)=\\
 &=\sum_{k=1}^\infty \phi^{O}_{k}(v)\Big(\E(\xi^{O}_{ik}X_i(u))-\lambda_k \tilde{\phi}_k^O(u)\Big).
\end{align*}
From the definition $\tilde{\phi}_k^O(u)$ in  \eqref{eq:PEF} we get that $\E(\xi^{O}_{ik}X_i^M(u))=\lambda^{O}_k\,\tilde{\phi}_k^O(u)$, which leads to $\E(X_i^{O}(v)\mathcal{Z}_i(u))=0$ for all $u\in M$.
This proves  \eqref{Zcorr}, while  \eqref{ZVar} directly follows from the definition of $\mathcal{Z}_i(u)$.

\smallskip

\noindent\textbf{Part (b):}
By Theorem \ref{th:SRO} there exists a unique $b_u\in H$ such that
$$\ell(X_i^{O})(u)=\langle b_u, X_i^{O}\rangle_H.$$
By  \eqref{Zcorr} and the orthogonality property of the least squares projection we thus obtain
\begin{align*}
&\E\left(\left(X_i(u)-\ell(X_i^{O})(u)\right)^2\right)=\\
&=\E\left(\left(\mathcal{L}(X_i^O)(u)+\mathcal{Z}_i(u)-\langle b_u, X_i^{O}\rangle_H\right)^2\right)=\\
&=\E\left(\left(\mathcal{L}(X_i^O)(u)-\langle b_u, X_i^{O}\rangle_H\right)^2\right)+\E(\mathcal{Z}_i(u)^2)+\\
&+2\left(\E(\mathcal{L}(X_i^O)(u)\,\mathcal{Z}_i(u))-
\E\left(\langle b_u,X_i^{O}\rangle_H \mathcal{Z}_i(u)\right)\right)=\\
&=\E\left(\left(\mathcal{L}(X_i^O)(u)-\langle b_u, X_i^{O}\rangle_H\right)^2\right)
+\E(\mathcal{Z}_i^2(u))\geq \E(\mathcal{Z}_i^2(u)).
\end{align*}

\smallskip

\noindent\textbf{Part (c):}
Observe that $\V(\mathcal{Z}_{i}(u)-\mathcal{Z}_{j}(u))=\V(\mathcal{Z}_{i}(u))+\V(\mathcal{Z}_{j}(u))-2\,\Cov(\mathcal{Z}_{i}(u),\mathcal{Z}_{j}(u))=2\,\V(\mathcal{Z}_{i}(u))$ for all $u\in M$ and $i\neq j$. Rearranging and using that $\E(\mathcal{Z}_{i}(u))=\E(\mathcal{Z}_{j}(u))=0$ for all $u\in M$ and all $i,j\in\{1,\dots,n\}$ yields $\V(\mathcal{Z}_i(u))=\frac{1}{2}\,\E((\mathcal{Z}_i(u)-\mathcal{Z}_j(u))^2)$. From result (a) we know that $\mathcal{Z}_i(u)$ and $X_i^{O}(v)$ are orthogonal and therefore uncorrelated for all $u\in M$ and all $v\in O$, that is, $\E\left(X_i^{O}(v)\mathcal{Z}_i(u)\right)=\mathrm{Cov}\left(X_i^{O}(v),\mathcal{Z}_i(u)\right)=0$. Under the assumption of an independent Gaussian process, we have then independence between $\mathcal{Z}_i(u)$ and $X_i^{O}(v)$, such that
\begin{align*}
\V(\mathcal{Z}_i(u))&=\frac{1}{2}\,\E\left(\E\left((\mathcal{Z}_i(u)-\mathcal{Z}_j(u))^2\right)|X_i^{O}=X_j^{O}\right),
\end{align*}
where $X_i^{O}=X_j^{O}$ means that $X_i^{O}(u)=X_i^{O}(u)$ for all $u\in O$. For the two random functions $\mathcal{Z}_i(u)$ and $\mathcal{Z}_j(u)$ we can write $\mathcal{Z}_i(u)=X_i^M(u)-\mathcal{L}(X_i^O)(u)$ and $\mathcal{Z}_j(u)=X_j^M(u)-\mathcal{L}(X_j^O)(u)$. It follows from the definition of $\tilde{\phi}^O_k(u)$ in  \eqref{eq:PEF} that $\mathcal{L}(X_i^O)(u)=\mathcal{L}(X_j^O)(u)$ for all $u\in M$, if and only if $X_i^{O}=X_j^{O}$. Therefore,
\begin{align*}
\V(\mathcal{Z}_i(u))&=\frac{1}{2}\,\E\left(\E\left((X^M_i(u)-X^M_j(u))^2\right)|X_i^{O}=X_j^{O}\right),
\end{align*}
for all $u\in M$.

\bigskip

\noindent\textbf{Proof of Theorem \ref{Theorem_UR}}
Note that under Assumption A2, the asymptotic rates of convergence of the LLK estimators
$\widehat{X}_i^{O}$ (see \eqref{X_cent_Estim}),
$\hat{\mu}$ (see \eqref{Estimator_mu_paper}), and
$\hat{\gamma}_u$ (see \eqref{Estimator_gamma_paper}) are asymptotically equivalent to the scenario, where $m_1=m_2=\dots =m_n$. Therefore, we consider the simplified case of a common number of discretization points $m$.

For proofing the results in Theorem \ref{Theorem_UR} we make use of the following two lemmas:

\begin{lemma}\label{Lemma_UR}
Define
\begin{align}\label{Def_Psi}
\Psi_{q,nm}(u;h_\mu)&=\frac{1}{nmh_\mu}\sum_{ij}\kappa\left(\frac{U_{ij}-u}{h_\mu}\right)\psi_q\left(U_{ij}-u,Y_{ij}\right),
\end{align}
where
\begin{align*}
\psi_q\left(U_{ij}-u,Y_{ij}\right)&=\left\{
  \begin{array}{ll}
    \left(U_{ij}-u\right)^q&\text{for }q\in\{0,1,2\}\\
    Y_{ij}&\text{for }q=3\\
    \left(U_{ij}-u\right)Y_{ij}&\text{for }q=4.
  \end{array}
\right.\notag
\end{align*}

Then, under Assumptions A1-A5,
\begin{align*}
  \tau_{q,nm}&=\sup_{u\in [a,b]}\left|\Psi_{q,nm}(u;h_\mu)-m_{q}(u)\right|
  =\mathcal{O}_p\left(h_\mu^2+\frac{1}{\sqrt{nm\,h_\mu}}+\frac{1}{\sqrt{n}}\right),
\end{align*}
where $m_0(u)=f_U(u)$, $m_1(u)=0$,
$m_2(u)=f_U(u)\nu_2(\kappa)$,
$m_3(u)=\mu(u)f_U(u)=\E(Y_{ij}|U_{ij}=u)f_U(u)$, and $m_4(u)=0$.
\end{lemma}

\smallskip

\begin{lemma}\label{Lemma_UR_gamma}
Define
{\small\begin{align}\label{Def_Theta}
\Theta_{q,n\mathcal{M}}(u,v;h_\gamma)&=\frac{1}{n\mathcal{M}h_\gamma}\sum_{i,j\neq l}\kappa\left(\frac{U_{ij}-u}{h_\gamma}\right)\kappa\left(\frac{U_{il}-v}{h_\gamma}\right)\vartheta_q\left(U_{ij}-u,U_{il}-u,C_{ijl}\right),
\end{align}}\noindent
where
\begin{align*}
\vartheta_q\left(U_{ij}-u,U_{il}-v,C_{ijl}\right)&=\left\{
  \begin{array}{ll}
    \left(U_{ij}-u\right)^q\left(U_{il}-v\right)^q&\text{for }q\in\{0,1,2\}\\
    C_{ijl}&\text{for }q=3\\
    \left(U_{ij}-u\right)\left(U_{il}-v\right)C_{ijl}&\text{for }q=4.
  \end{array}
\right.\notag
\end{align*}

Then, under Assumptions A1-A5,
\begin{align*}
  \varrho_{q,n\mathcal{M}}&=\sup_{(u,v)\in
    [a,b]^2}\left|\Theta_{q,n\mathcal{M}}(u,v;h_\gamma)-\eta_{q}(u,v)\right|
=\mathcal{O}_p\left(h_\gamma^2+\frac{1}{\sqrt{n\mathcal{M}\,h^2_\gamma}}+\frac{1}{\sqrt{n}}\right),
\end{align*}
where $\eta_0(u,v)=f_{UU}(u,v)$, $\eta_1(u,v)=0$,
$\eta_2(u,v)=f_{UU}(u,v)(\nu_2(\kappa))^2$,
$\eta_3(u,v)=\gamma(u,v)f_{UU}(u,v)=\E(C_{ijl}|(U_{ij},U_{il})=(u,v))f_{UU}(u,v)$, and\\ $m_4(u,v)=0$.
\end{lemma}

\bigskip

\noindent\textbf{Proof of Lemma \ref{Lemma_UR}:}\
Remember that $\E(|\tau_{q,nm}|)=\mathcal{O}(\mathtt{rate}_{nm})$ implies that
$\tau_{q,nm}=\mathcal{O}_p(\mathtt{rate}_{nm})$, therefore, we focus in the
following on $\E(|\tau_{q,nm}|)$, where $\E(|\tau_{q,nm}|)=\E(\tau_{q,nm})$. Adding a zero and applying the
triangle inequality yields that $\E(\tau_{q,nm})=$
\begin{align}\label{UR_Ineq}
 \E(\sup_{u\in [a,b]}\left|\Psi_{q,nm}(u;h_\mu)-m_{q}(u)\right|)
  &\leq\sup_{u\in [a,b]}|\E(\Psi_{q,nm}(u;h_\mu))-m_p(u)|+\notag\\
  &+\E(\sup_{u\in [a,b]}\left|\Psi_{q,nm}(u;h_\mu)-\E(\Psi_{q,nm}(u;h_\mu))\right|).
\end{align}
Let us first focus on the second summand in \eqref{UR_Ineq}.
The next steps will make use of the Fourier transformation
of the kernel function $\kappa$ \citepappendix[see, e.g.,][Ch.~1.3]{tsybakov2009intro}:
\begin{align*}
\kappa^{\mathrm{ft}}(x):=\mathcal{F}[\kappa](x) =\int_{\mathbb{R}}\kappa(z)\exp(-\mathrm{i}zx)dz
=\int_{-1}^1\kappa(z)\exp(-\mathrm{i}zx)dz
\end{align*}
with $\mathrm{i}=\sqrt{-1}$. By Assumption A5, $\kappa(.)$
has a compact support $[-1,1]$. The inverse transform gives then
\begin{align*}
  \kappa\left(s\right)&=
  \frac{1}{2\pi}\int_{\mathbb{R}}\kappa^{\mathrm{ft}}(x)\exp\left(\mathrm{i}
    xs\right)dx
 =\frac{1}{2\pi}\int_{\mathbb{R}}\kappa^{\mathrm{ft}}(x)\exp\left(\mathrm{i}
    xs\right)dx\;\mathbbm{1}_{(|s|<1)}.
\end{align*}
Furthermore, we can use that \citepappendix[see][Ch.~1.3,  (1.34)]{tsybakov2009intro}
$\mathcal{F}[\kappa(./h_\mu)/h_\mu](x)=\mathcal{F}[\kappa](xh_\mu)=\kappa^{\mathrm{ft}}(xh_\mu)$
which yields
\begin{align}
  \kappa\left(s/h_\mu\right)/h_\mu
&= \frac{1}{2\pi}\int_{\mathbb{R}}\mathcal{F}[\kappa(./h_\mu)/h_\mu](x)\exp\left(\mathrm{i}xs\right)dx\;\mathbbm{1}_{(|s|<h_\mu)}\notag\\
  &=
  \frac{1}{2\pi}\int_{\mathbb{R}}\kappa^{\mathrm{ft}}(xh_\mu)\exp\left(\mathrm{i}xs\right)dx\;\mathbbm{1}_{(|s|<h_\mu)}.\label{kappa_transf}
\end{align}
Plugging \eqref{kappa_transf} into \eqref{Def_Psi} yields $\Psi_{q,nm}(u;h_\mu)=$
\begin{align*}
&=\frac{1}{nm}\sum_{ij}\kappa\left(\frac{U_{ij}-u}{h_\mu}\right)\frac{1}{h_\mu}\psi_q\left(U_{ij}-u,Y_{ij}\right)\\
&=\frac{1}{nm}\sum_{ij}\frac{1}{2\pi}\int_{\mathbb{R}}\kappa^{\mathrm{ft}}(xh_\mu)\exp\big(\mathrm{i}x(U_{ij}-u)\big)dx \;\mathbbm{1}_{(|U_{ij}-u|<h_\mu)}
\;\psi_q\left(U_{ij}-u,Y_{ij}\right)\\
&=\frac{1}{2\pi}\int_{\mathbb{R}}\left[\frac{1}{nm}\sum_{ij}\exp\big(\mathrm{i}xU_{ij}\big)\,\psi_q\left(U_{ij}-u,Y_{ij}\right)\;\mathbbm{1}_{(|U_{ij}-u|<h_\mu)}\right]\exp\big(\mathrm{i}xu\big)\kappa^{\mathrm{ft}}(xh_\mu)dx.
\end{align*}
Using that $|\exp(ixu)|\leq 1$ leads to
{\small\begin{align*}
  \E(\sup_{u\in
    [a,b]}\left|\Psi_{q,nm}(u;h_\mu)-\E(\Psi_{q,nm}(u;h_\mu))\right|)
\leq\frac{1}{2\pi}\E\left(\sup_{u\in [a,b]}\left|\int_{\mathbb{R}}\tilde\omega_{q,nm}(u,x)
  \cdot \kappa^{\mathrm{ft}}(xh_\mu)dx\right|\right),
\end{align*}}
where
\begin{align*}
\tilde\omega_{q,nm}(u,x)=&\frac{1}{nm}\sum_{ij}\big[\exp\big(\mathrm{i}xU_{ij}\big)\psi_q\left(U_{ij}-u,Y_{ij}\right)\,\mathbbm{1}_{(|U_{ij}-u|<h_\mu)}-\\
&\hspace*{1.0cm}\E\left(\exp\big(\mathrm{i}xU_{ij}\big)\psi_q\left(U_{ij}-u,Y_{ij}\right)\mathbbm{1}_{(|U_{ij}-u|<h_\mu)}\right)\big].
\end{align*}
Using further that $\kappa^{\mathrm{ft}}$ is symmetric, since $\kappa$ is
symmetric by Assumption A5, and that
$\exp\big(\mathrm{i}xU_{ij}\big)=\cos\big(xU_{ij}\big)+\mathrm{i}\sin\big(xU_{ij}\big)$
leads to
{\small\begin{align*}
\frac{1}{2\pi}\E\left(\sup_{u\in [a,b]}\left|\int_{\mathbb{R}}\tilde\omega_{q,nm}(u,x)
  \cdot \kappa^{\mathrm{ft}}(xh_\mu)dx\right|\right)
=\frac{1}{2\pi}\E\left(\sup_{u\in [a,b]}\left|\int_{\mathbb{R}}\omega_{q,nm}(u,x)
  \cdot \kappa^{\mathrm{ft}}(xh_\mu)dx\right|\right),
\end{align*}}
where
\begin{align}\label{omega_fun}
\omega_{q,nm}(u,x)=&\frac{1}{nm}\sum_{ij}\Big[\cos\big(xU_{ij}\big)\psi_q\left(U_{ij}-u,Y_{ij}\right)\,\mathbbm{1}_{(|U_{ij}-u|<h_\mu)}-\notag\\
&\hspace*{1.0cm}\E\Big(\cos\big(xU_{ij}\big)\psi_q\left(U_{ij}-u,Y_{ij}\right)\mathbbm{1}_{(|U_{ij}-u|<h_\mu)}\Big)\Big],
\end{align}
such that
\begin{align}
&\E(\sup_{u\in
    [a,b]}\left|\Psi_{q,nm}(u;h_\mu)-\E(\Psi_{q,nm}(u;h_\mu))\right|)\notag\\
&\leq\frac{1}{2\pi}\int_{\mathbb{R}}\E\left(\sup_{u\in [a,b]}\big|\omega_{q,nm}(u,x)
  \big|\right)\cdot\left|\kappa^{\mathrm{ft}}(xh_\mu)\right|dx\notag\\
&\leq\frac{1}{2\pi}\int_{\mathbb{R}}\sqrt{\E\left(\left(\sup_{u\in [a,b]}\big|\omega_{q,nm}(u,x)
  \big|\right)^2\right)}\cdot\left|\kappa^{\mathrm{ft}}(xh_\mu)\right|dx\notag\\
&=\frac{1}{2\pi}\int_{\mathbb{R}}\sqrt{\E\left(\sup_{u\in [a,b]}\;\big(\omega_{q,nm}(u,x)
  \big)^2\right)}\cdot\left|\kappa^{\mathrm{ft}}(xh_\mu)\right|dx.\label{Approx_1}
\end{align}


In order to simplify the notation we will denote
\begin{align*}
W_{ij}^q(x,u)&=\cos\big(xU_{ij}\big)\psi_q\left(U_{ij}-u,Y_{ij}\right),
\end{align*}
such that $\E\left(\sup_{u\in [a,b]}\,\left(\omega_{q,nm}(u,x)\right)^2\right)=$
\begin{align*}
\E\Big(\sup_{u\in [a,b]}\,\Big(\frac{1}{(nm)^2}\sum_{ij}&\left[W_{ij}^q(x,u)\mathbbm{1}_{(|U_{ij}-u|<h_\mu)}-\E(W_{ij}^q(x,u)\mathbbm{1}_{(|U_{ij}-u|<h_\mu)})\right]^2+\\
\frac{1}{(nm)^2}\sum_{(i,j)\neq(r,l)}&\big[(W_{ij}^q(x,u)\mathbbm{1}_{(|U_{ij}-u|<h_\mu)}-\E(W_{ij}^q(x,u)\mathbbm{1}_{(|U_{ij}-u|<h_\mu)}))\cdot\\
\hspace*{0cm}\cdot&
(W_{rl}^q(x,u)\mathbbm{1}_{(|U_{rl}-u|<h_\mu)}-\E(W_{rl}^q(x,u)\mathbbm{1}_{(|U_{rl}-u|<h_\mu)}))\big]\Big)\Big).
\end{align*}

As $u$ takes only values within the compact interval $[a,b]$, there
exist constants $C_1$ and $C_2$ such that, uniformly for all $u\in[a,b]$,
$\Prob(|U_{ij}-u|<h_\mu)\leq C_1 h_\mu<\infty$, for all $i,j$, and
$\Prob(|U_{ij}-u|<h_\mu\text{ AND }|U_{rl}-u|<h_\mu)\leq C_2
h_\mu^2<\infty$, for all $(i,j)\neq (r,l)$. Together with the
triangle inequality, this yields that
$\E\left(\sup_{u\in [a,b]}\,\left(\omega_{q,nm}(u,x)\right)^2\right)\leq$
{\small\begin{align*}
&\frac{C_1 h_\mu}{(nm)^2}\sum_{ij}\E\left(\sup_{u\in [a,b]}\,\left[W_{ij}^q(x,u)-\E(W_{ij}^q(x,u))\right]^2\right)+\\
&\frac{C_2 h_\mu^2}{(nm)^2}\sum_{(i,j)\neq (r,l)}\E\left(\sup_{u\in [a,b]}\,\big[(W_{ij}^q(x,u)-\E(W_{ij}^q(x,u)))
(W_{rl}^q(x,u)-\E(W_{rl}^q(x,u)))\big]\right)
\end{align*}}

From our moment assumptions (Assumption
A1) and the fact that $[a,b]$ is compact, we can conclude that there must
exist a constant $C_3$ such that, point-wise for every $x\in\mathbb{R}$,
\begin{align}\label{VAR_approx}
\E\Big(\big(\sup_{u\in
  [a,b]}\big|W_{ij}^q(x,u)-\E(W_{ij}^q(x,u))\big|\big)^2\Big)\leq
C_3<\infty
\end{align}
for all $i$ and $j$.

``Within function'' dependencies: By the same reasoning there
must exist a constant $C_4$ such that, point-wise for every $x\in\mathbb{R}$,
\begin{align}\label{Within_approx}
&\E\Big(\sup_{u\in
  [a,b]}\big|W_{ij}^q(x,u)-\E(W_{ij}^q(x,u))\big|\cdot\sup_{u\in
  [a,b]}\big|W_{il}^q(x,u)-\E(W_{il}^q(x,u))\big|\Big)\notag\\
& \leq C_4<\infty
\end{align}
for all $j\neq l$ and all $i$.

``Between function'' dependencies: Our weak dependency assumption (Assumption A1) and the fact that $[a,b]$
is compact yields that point-wise for every $x\in\mathbb{R}$
\begin{align}\label{Between_approx}
&\E\Big(\sup_{u\in
  [a,b]}\big|W_{ij}^q(x,u)-\E(W_{ij}^q(x,u))\big|\cdot\sup_{u\in
  [a,b]}\big|W_{rl}^q(x,u)-\E(W_{rl}^q(x,u))\big|\Big)\notag\\
& \leq c_1 \iota_1^{|i-r|}
\end{align}
for all $j,l$ and $|i-r|\geq 1$, where $0<c_1<\infty$ and $0<\iota_1<1$.

Eq.s \eqref{VAR_approx}, \eqref{Within_approx}, and
\eqref{Between_approx} yield that
$\E\left(\sup_{u\in [a,b]}\left(\omega_{q,nm}(u,x)\right)^2\right)\leq$
\begin{align*}
\leq&\frac{C_1h_\mu}{(nm)^2}\sum_{ij} C_3+
\frac{C_2h^2_\mu}{(nm)^2}\sum_{i,j\neq l} C_4+
\frac{C_2h^2_\mu}{(nm)^2}\sum_{i\neq r,jl} c_1 \iota_1^{|i-r|}
\\
=&\mathcal{O}\left(\frac{h_\mu}{nm}+\frac{h_\mu^2(m-1)}{nm}+\frac{h_\mu^2}{n}\right)
=\mathcal{O}\left(\frac{h_\mu}{nm}+\frac{h_\mu^2}{n}\right)
,
\end{align*}
such that
\begin{align}
\sqrt{\E\left(\sup_{u\in [a,b]}\left(\omega_{q,nm}(u,x)\right)^2\right)}=\mathcal{O}\left(\sqrt{\frac{h_\mu}{nm}}+\frac{h_\mu}{\sqrt{n}}\right).\label{Approx_2}
\end{align}

Plugging
\eqref{Approx_2} into \eqref{Approx_1} and integration by substitution
leads to
\begin{align}
  &\E(\sup_{u\in [a,b]}\left|\Psi_{q,nm}(u;h_\mu)-\E(\Psi_{q,nm}(u;h_\mu))\right|)\leq\notag\\
  &\frac{1}{2\pi}\int_{\mathbb{R}}\sqrt{\E\left(\sup_{u\in [a,b]}\left(\omega_{q,nm}(u,x)\right)^2\right)}\cdot\left|\kappa^{\mathrm{ft}}(xh_\mu)\right|dx
=\mathcal{O}\left(\frac{1}{\sqrt{nm\,h_\mu}}+\frac{1}{\sqrt{n}}\right).\label{UR_result1}
\end{align}

Let us now focus on the first summand in \eqref{UR_Ineq}.
From standard arguments in nonparametric statistics (see, e.g., \cite{ruppert1994}) we know that
\begin{align*}
  \E(\Psi_{q,nm}(u;h_\mu))-m_q(u)=\mathcal{O}(h_\mu^2)
\end{align*}
for each $u\in [a,b]$ and for all
$q\in\{0,\dots,4\}$. Under our smoothness Assumption A3, the ``$\mathcal{O}(h_\mu^2)$''
term becomes uniformly valid for all $u\in [a,b]$ and all $q\in\{0,1,2,4\}$, since all of the involved functions have uniformly bounded second order derivatives.
We can conclude with respect to the first summand
in \eqref{UR_Ineq} that
\begin{align}\label{UR_result2}
\sup_{u\in [a,b]}|\E(\Psi_{q,nm}(u;h_\mu))-m_p(u)|=\mathcal{O}(h_\mu^2) \quad\text{for
all}\quad q\in\{0,\dots,4\}.
\end{align}

Finally, plugging our results \eqref{UR_result1} and
\eqref{UR_result2} into \eqref{UR_Ineq} leads to
\begin{align}\label{UR_result3}
\tau_{q,nm}&=\sup_{u\in
  [a,b]}\left|\Psi_{q,nm}(u;h_\mu)-m_{q}(u)\right|=\mathcal{O}_p\left(h_\mu^2+\frac{1}{\sqrt{nm\,h_\mu}}+\frac{1}{\sqrt{n}}\right)
\end{align}
for all $q\in\{0,\dots,4\}$.

\bigskip

\noindent\textbf{Proof of Lemma \ref{Lemma_UR_gamma}:}
Analogously to that of Lemma \ref{Lemma_UR}.

\bigskip


\noindent\textbf{Proof of Theorem \ref{Theorem_UR}, part (a):}
Let us rewrite the estimator $\hat{\mu}$ using matrix notation
as in \cite{ruppert1994}, i.e.,
\begin{align}
\hat{\mu}(u;h_{\mu})=e_1^{\top}\left([\mathbf{1},\mathbf{U}_{u}]^{\top}\mathbf{W}_{\mu,u}[\mathbf{1},\mathbf{U}_{u}]\right)^{-1}[\mathbf{1},\mathbf{U}_{u}]^{\top}\mathbf{W}_{\mu,u}\mathbf{Y},\label{Estimator_mu}
\end{align}
where $e_1=(1,0)^{\top}$, $[\mathbf{1},\mathbf{U}_{u}]$ is a $nm\times
2$ dimensional data matrix with typical rows $(1,U_{ij}-u)$, the
$nm\times nm$ dimensional diagonal weighting matrix
$\mathbf{W}_{\mu,u}$ holds the kernel weights
$K_{\mu,h}(U_{ij}-u)=h^{-1}_{\mu}\,\kappa(h^{-1}_{\mu}(U_{ij}-u))$. The objects $\mathbf{U}_{u}$ and
$\mathbf{W}_{\mu,u}$ are filled in correspondence with the $nm$
dimensional vector
$\mathbf{Y}=(Y_{11},Y_{12},\dots,Y_{n,m-1},Y_{n,m})^{\top}$.

This way we can decompose the estimator $\hat{\mu}(u;h_\mu)$ as
\begin{align}\label{estimator_mean}
\hat{\mu}(u;h_\mu)&=e_1^{\top}L_{1,nm,u}^{-1} L_{2,nm,u},
\end{align}
with $2\times 2$ matrix
\begin{align}
L_{1,nm,u}&=(nm)^{-1}[\mathbf{1},\mathbf{U}_{u}]^{\top}\mathbf{W}_{\mu,u}[\mathbf{1},\mathbf{U}_{u}]\notag\\
&=\left(\begin{matrix}
      \frac{1}{nmh_\mu}\sum_{ij}\kappa\left(\frac{U_{ij}-u}{h_\mu}\right)&
\frac{1}{nmh_\mu}\sum_{ij}\kappa\left(\frac{U_{ij}-u}{h_\mu}\right)\left(U_{ij}-u\right)\\
    \frac{1}{nmh_\mu}\sum_{ij}\kappa\left(\frac{U_{ij}-u}{h_\mu}\right)\left(U_{ij}-u\right)&
\frac{1}{nmh_\mu}\sum_{ij}\kappa\left(\frac{U_{ij}-u}{h_\mu}\right)\left(U_{ij}-u\right)^2
    \end{matrix}\right),\notag
\end{align}
and $2\times 1$ vector
\begin{align}
L_{2,nm,u}&=(nm)^{-1}[\mathbf{1},\mathbf{U}_{u}]^{\top}\mathbf{W}_{\mu,u}\mathbf{Y}
=\left(\begin{matrix}
\frac{1}{nmh_\mu}\sum_{ij}\kappa\left(\frac{U_{ij}-u}{h_\mu}\right)Y_{ij}\\
\frac{1}{nmh_\mu}\sum_{ij}\kappa\left(\frac{U_{ij}-u}{h_\mu}\right)\left(U_{ij}-u\right)Y_{ij}\\
\end{matrix}\right).\notag
\end{align}

Using the notation and the results from Lemma \ref{Lemma_UR} we have that
\begin{align}
  L_{1,nm,u}
&=\left(\begin{matrix}
    \Psi_{0,nm}(u;h_\mu)&\Psi_{1,nm}(u;h_\mu)\\
    \Psi_{1,nm}(u;h_\mu)&\Psi_{2,nm}(u;h_\mu)
  \end{matrix}
\right)\notag\\
&=\left(\begin{matrix}
    f_U(u) & 0\\
    0      & f_U(u)\nu_2(\kappa)
  \end{matrix}
\right)+\mathcal{O}^{\operatorname{Unif}}_p\left(h_\mu^2+\frac{1}{\sqrt{nm\,h_\mu}}+\frac{1}{\sqrt{n}}\right)\quad\text{and}
\label{S1}\\
L_{2,nm,u}
  &=\left(\begin{matrix}
      \Psi_{3,nm}(u;h_\mu)\\
      \Psi_{4,nm}(u;h_\mu)
    \end{matrix}
  \right)
  =\left(\begin{matrix}
      \mu(u)f_U(u)\\
      0
    \end{matrix}
  \right)+\mathcal{O}^{\operatorname{Unif}}_p\left(h_\mu^2+\frac{1}{\sqrt{nm\,h_\mu}}+\frac{1}{\sqrt{n}}\right),\label{S2}
\end{align}
where we write
$\Psi_{q,nm}(u;h_\mu)-m_q(u)=\mathcal{O}^{\operatorname{Unif}}_p(\texttt{rate})$
in order to denote that
$\sup_{u\in[a,b]}|\Psi_{q,nm}(u;h_\mu)-m_q(u)|=\mathcal{O}_p(\texttt{rate})$. Taking
the inverse of \eqref{S1} gives
\begin{align}\label{S1_inv}
  L_{nm,u}^{-1}
  &=
  \left(\begin{matrix}
      1/f_U(u) & 0\\
      0      & 1/(f_U(u)\nu_2(\kappa))
    \end{matrix}\right)
  +\mathcal{O}_p^{\operatorname{Unif}}\left(h_\mu^2+\frac{1}{\sqrt{nm\,h_\mu}}+\frac{1}{\sqrt{n}}\right).
\end{align}

Plugging \eqref{S1_inv} and \eqref{S2} into \eqref{estimator_mean} leads to


\begin{align*}
\sup_{u\in [a,b]}|\hat{\mu}(u;h_\mu)-\mu(u)|&=\mathcal{O}_p\left(h_\mu^2+\frac{1}{\sqrt{nm\,h_\mu}}+\frac{1}{\sqrt{m}}\right).
\end{align*}

\smallskip

\noindent\textbf{Proof of Theorem \ref{Theorem_UR}, part (\~a):}
Observe that
\begin{align*}
&\sup_{u\in O}|\widehat{X}_i^O(u;h_\mu,h_X)-X_i^O(u)|\leq\\
&\sup_{u\in O}|\hat{X}_i^{c,O}(u;h_X)-(X_i^O(u)-\mu(u))|+
\sup_{u\in O}|\hat{\mu}(u;h_\mu)-\mu(u)|.
\end{align*}
From Theorem \ref{Theorem_UR}, part (a), we have that $\sup_{u\in O}|\hat{\mu}(u;h_\mu)-\mu(u)|=\mathcal{O}_p(r_\mu)$ with $r_\mu=h_\mu^2+1/\sqrt{nm\,h_\mu}+1/\sqrt{m}$. From a simplified version of the proof of Theorem \ref{Theorem_UR}, part (a), with $n=1$, it follows that
$$
\sup_{u\in O}|\hat{X}_i^{c,O}(u;h_X)-(X_i^O(u)-\mu(u))|=\mathcal{O}_p\left(h_X^2+\frac{1}{\sqrt{mh_X}}\right)
$$

\bigskip

\noindent\textbf{Proof of Theorem \ref{Theorem_UR}, part (b):}


Let us rewrite the estimator $\hat{\gamma}$ using matrix notation
as in \cite{ruppert1994}, i.e., $\hat{\gamma}(u,v;h_{\gamma})=$
\begin{equation}\begin{array}{c}
=e_1^{\top}\left([\mathbf{1},\mathbf{U}_{u},\mathbf{U}_{v}]^{\top}\mathbf{W}_{\gamma,u,v}[\mathbf{1},\mathbf{U}_{u},\mathbf{U}_{v}]\right)^{-1}[\mathbf{1},\mathbf{U}_{u},\mathbf{U}_{v}]^{\top}\mathbf{W}_{\gamma,u,v}\hat{\mathbf{C}},\label{Estimator_gamma}
\end{array}\end{equation}
where $e_1=(1,0,0)^{\top}$, $[\mathbf{1},\mathbf{U}_{u},\mathbf{U}_{v}]$ is a
$n\mathcal{M}\times 3$ dimensional data matrix with typical rows
$(1,U_{ij}-u,U_{il}-v)$, the $n\mathcal{M}\times n\mathcal{M}$ dimensional diagonal weighting
matrix $\mathbf{W}_{\gamma,u,v}$ holds the bivariate kernel weights
$K_{\gamma,h}(U_{ij}-u,U_{il}-v)$. For the bivariate kernel weights
$K_{\gamma,h}(z_1,z_2)=h^{-2}_\gamma\kappa_{\gamma}(z_1,z_2)$ we use a
multiplicative kernel function
$\kappa_{\gamma}(z_1,z_2)=\kappa(z_1)\kappa(z_2)$ with $\kappa$ as defined
above. The usual kernel constants are then
$\nu_{2}(\kappa_\gamma):=\left(\nu_{2}(\kappa)\right)^2$ and
$R(\kappa_\gamma):=R(\kappa)^2$. The rows of the matrices
$[\mathbf{1},\mathbf{U}_{u},\mathbf{U}_{v}]$ and $\mathbf{W}_{\gamma,u,v}$ are
filled in correspondence with the $n\mathcal{M}$ elements of the vector of raw-covariances
$\hat{\mathbf{C}}=(\dots,\hat{C}_{ijl}, \dots)^{\top}$.

Let us initially consider the infeasible estimator $\hat\gamma_C$ that is
based on the infeasible ``clean'' raw-covariances $C_{ijl}=(Y_{ij} - \mu(U_{ij}))(Y_{il}
- \mu(U_{il}))$ instead of the estimator $\hat\gamma$ in
\eqref{Estimator_gamma} that is based on the ``dirty'' raw-covariances
 $\hat{C}_{ijl}=(Y_{ij} - \hat{\mu}(U_{ij}))(Y_{il}
- \hat{\mu}(U_{il}))$.
Equivalently to the estimator $\hat{\mu}$ above, we can write the estimator $\hat\gamma_C$ as
\begin{align}\label{estimator_gamma}
\hat{\gamma}_C(u,v;h_\gamma)&=e_1^{\top}\tilde{S}_{1,n\mathcal{M},(u,v)}^{-1}\tilde{S}_{2,n\mathcal{M},(u,v)},
\end{align}
with
{\small\begin{align}
  \tilde{S}_{1,n\mathcal{M},(u,v)}^{-1}
&=\left(\begin{matrix}
    \Theta_{0,n\mathcal{M}}(u,v;h_\gamma)&\Theta_{1,n\mathcal{M}}(u,v;h_\gamma)\\
    \Theta_{1,n\mathcal{M}}(u,v;h_\gamma)&\Theta_{2,n\mathcal{M}}(u,v;h_\gamma)
  \end{matrix}
\right)^{-1}\notag\\
&=\left(\begin{matrix}
    1/f_{UU}(u,v) & 0\\
    0      & 1/f_{UU}(u,v)(\nu_2(\kappa))^2
  \end{matrix}
\right)+\mathcal{O}^{\operatorname{Unif}}_p\left(h_\gamma^2+\frac{1}{\sqrt{n\mathcal{M}\,h^2_\gamma}}+\frac{1}{\sqrt{n}}\right)
\label{S1_inv_tilde}
\end{align}}
and $\tilde{S}_{2,n\mathcal{M},(u,v)}=$
{\small\begin{align}
  &=\left(\begin{matrix}
      \Theta_{3,n\mathcal{M}}(u,v;h_\gamma)\\
      \Theta_{4,n\mathcal{M}}(u,v;h_\gamma)
    \end{matrix}
  \right)
  =\left(\begin{matrix}
      \gamma(u,v)f_{UU}(u,v)\\
      0
    \end{matrix}
  \right)+\mathcal{O}^{\operatorname{Unif}}_p\left(h_\gamma^2+\frac{1}{\sqrt{n\mathcal{M}\,h^2_\gamma}}+\frac{1}{\sqrt{n}}\right),\label{S2_tilde}
\end{align}}
where we use the notation and the results from Lemma \ref{Lemma_UR_gamma},
and where we write
$\Theta_{q,n\mathcal{M}}(u,v;h_\gamma)-\eta_q(u,v)=\mathcal{O}^{\operatorname{Unif}}_p(\texttt{rate})$
in order to denote that
$\sup_{(u,v)\in [a,b]^2}|\Theta_{q,n\mathcal{M}}(u,v;h_\gamma)-\eta_q(u,v)|=\mathcal{O}_p(\texttt{rate})$.

Plugging \eqref{S1_inv_tilde} and \eqref{S2_tilde} into \eqref{estimator_gamma} leads to
\begin{align}\label{gamma_interm_res}
\sup_{(u,v)\in [a,b]^2}|\hat{\gamma}_C(u,v;h_\mu)-\gamma(u,v)|&=\mathcal{O}_p\left(h_\gamma^2+\frac{1}{\sqrt{n\mathcal{M}\,h^2_\gamma}}+\frac{1}{\sqrt{n}}\right).
\end{align}

It remains to consider the additional estimation error, which comes from
using the ``dirty'' response variables $\hat{C}_{ijl}$ instead of
``clean'' dependent variables $C_{ijl}$.
Observe that we can expand $\hat{C}_{ijl}$ as following:
\begin{align*}
\hat{C}_{ijl} &= C_{ijl} + (Y_{ij} - \mu(U_{ij}))(\mu(U_{il}) - \hat{\mu}(U_{il}))\\
&+ (Y_{il} - \mu(U_{il}))(\mu(U_{ij}) - \hat{\mu}(U_{ij}))\\
&+ (\mu(U_{ij}) - \hat{\mu}(U_{ij}))(\mu(U_{il})-\hat{\mu}(U_{il})).
\end{align*}
Using our finite moment assumptions on $Y_{ij}$ (Assumption A1) and
our result in Theorem
\ref{Theorem_UR}, part (a), we have that
{\small\begin{align*}
\hat{C}_{ijl} &= C_{ijl} + \mathcal{O}_p(1)\mathcal{O}_p\left(h_\mu^2+\frac{1}{\sqrt{nm\,h_\mu}}+\frac{1}{\sqrt{n}}\right)\\
&+ \mathcal{O}_p(1)\mathcal{O}_p\left(h_\mu^2+\frac{1}{\sqrt{nm\,h_\mu}}+\frac{1}{\sqrt{n}}\right)\\
&+
\left(\mathcal{O}_p\left(h_\mu^2+\frac{1}{\sqrt{nm\,h_\mu}}+\frac{1}{\sqrt{n}}\right)\right)^2
=C_{ijl} + \mathcal{O}_p\left(h_\mu^2+\frac{1}{\sqrt{nm\,h_\mu}}+\frac{1}{\sqrt{n}}\right),
\end{align*}}
\noindent
uniformly for all $j\neq l\in\{1,\dots,m\}$ and $i\in\{1,\dots,n\}$. Therefore
\begin{align*}
\sup_{(u,v)\in [a,b]^2}|\hat{\gamma}(u,v;h_\mu)-\gamma(u,v)|&=\mathcal{O}_p\left(h_\gamma^2+h_\mu^2+\frac{1}{\sqrt{n\mathcal{M}\,h^2_\gamma}}+\frac{1}{\sqrt{nm\,h^2_\mu}}+\frac{1}{\sqrt{n}}\right).
\end{align*}

\bigskip

\noindent\textbf{Proof of Theorem \ref{Theorem_UR}, parts (c) and (d):}
Part (c) follows directly from inequality $\sup_{k\geq 1}|\hat\lambda_k^O - \lambda_k^O|\leq ||\hat\gamma - \gamma||_2$; see inequality (4.43) in \cite{bosq2000linear}. Part (d) follows directly from Lemma 4.3 in \cite{bosq2000linear}. 

In the following let $O:=O_i=[A_i,B_i]$ for some  $i\in {1,\dots,n}$. By assumption of Theorem \ref{Theorem_Main_Estim} we have $B_i-A_i\geq \ell_{\min}$, and recall that by Assumption (A1) the structure of a function $X_i$, to be observed on $O_i$, does not depend on the specific interval $O_i$.

For the proof of Theorem \ref{Theorem_Main_Estim} we need some additional lemmas. Generally note that under the assumed choice of bandwidths we have $r_\mu+r_\gamma\asymp r_{mn}$ for $r_{mn}=\frac{1}{\min\{n^{1/2},(n\mathcal{M})^{1/3}\}}$, since for all $n$ and $m\asymp n^\theta$ sufficiently large we have $(mn)^{2/5}\geq \min\{n^{1/2},(n\mathcal{M})^{1/3}\}$.

Recall that $\hat{\tilde{\phi}}^O_{k}(u)=
\frac{\langle\hat\phi_k^{O},\hat\gamma_u\rangle_2}{\hat\lambda_k^O}$ and $\tilde{\phi}^O_{k}(u)=
\frac{\langle \phi_k^{O},\gamma_u\rangle_2}{\lambda_k^O}$ for $u\in O\cup M$, where in the particular case of $u\in O$ we additionally have  $\hat{\tilde{\phi}}^O_{k}(u)=\hat\phi_k^O(u)$ and $\tilde{\phi}_{k}^O(u)=\phi_k^{O}(u)$. Also recall that $\bar{K}_{mn}^{a_O+3/2}r_{mn}=O(1)$ and that by (A6) we have $\delta_k^O=O(k^{-a_O-1})$ as well as $1/\delta_k^O=O(k^{a_O+1})$.

\begin{lemma}\label{Lemma_phitilde}
 Under the assumptions of Theorem \ref{Theorem_Main_Estim} we have for all  $1\leq k\leq K\leq \bar{K}_{mn}$
\begin{align}\label{phi_tilde_rate}
\sup_{u\in O\cup M}\sup_{1\leq k\leq K} \delta_k^O|\hat{\tilde{\phi}}^O_k(u)-\tilde{\phi}^O_k(u)|=
\mathcal{O}_p\left(K^{1/2}r_{mn}\right),
\end{align}
\end{lemma}

\noindent\textbf{Proof of Lemma \ref{Lemma_phitilde}:}\
Using results (b) and  (c) of Theorem \ref{Theorem_UR} we obtain
\begin{align}
\hat{\tilde{\phi}}^O_{k}(u)=\frac{\langle\hat\phi_k^{O},\gamma_u\rangle_2}{\hat\lambda_k^O}+\mathcal{R}_{1,k}(u),
 \ \sup_{u\in O\cup M}\sup_{1\leq k\leq K} \lambda_k^O \mathcal{R}_{1,k}(u)=\mathcal{O}_p\left(r_{mn}\right).
\label{phitilde1}
\end{align}
But by the established properties (in particular \eqref{Zcorr} in Theorem \ref{OptimalPrediction}) of our operator we have $\gamma(u,v)=\sum_{j=1}^\infty
\lambda_j^O \tilde{\phi}^O_j(u)\phi^O_j(v)$ for all $u\in O\cup M$ and $v\in O$.
Hence
\begin{align}
\frac{\langle\hat\phi_k^{O},\gamma_u\rangle_2}{\hat\lambda_k^O}= \frac{1}{\hat\lambda_k^O} \sum_{j=1}^\infty
\lambda_j^O \tilde{\phi}^O_j(u)\langle\hat\phi_k^{O},\phi_j^{O}\rangle_2
\label{phitilde2}
\end{align}
Now note that for all $j\geq 1$
\begin{align}
\lambda_j^O \langle\hat\phi_k^{O},\phi_j^{O}\rangle_2
&=\int_{O^2} \gamma(u,v)\hat\phi_k^{O}(u)\phi_j^{O}(v)dudv \nonumber \\
&=\hat\lambda_k^O\langle\hat\phi_k^{O},\phi_j^{O}\rangle_2  + \int_{O^2} (\gamma(u,v)-\hat\gamma(u,v))\hat\phi_k^{O}(u)\phi_j^{O}(v)dudv \label{phitilde2a}
\end{align}
Let $v_{\hat\gamma,\gamma,k}(u):= \int_{O}(\gamma(u,v)-\hat\gamma(u,v))\hat\phi_k^{O}(u)du$. By the orthonormality of the system $\phi_1^{O},\phi_2^{O},\dots$ of eigenfunctions, the Cauchy-Schwarz inequality, and Theorem \ref{Theorem_UR} we have
 \begin{align}
&\sup_{1\leq k\leq K} \sum_{j=1}^\infty \left(\int_{O^2} (\gamma(u,v)-\hat\gamma(u,v))\hat\phi_k^{O}(u)\phi_j^{O}(v)dudv\right)^2 \nonumber \\
&= \sup_{1\leq k\leq K} \sum_{j=1}^\infty \langle v_{\hat\gamma,\gamma,k}, \phi_j^{O}\rangle_2^2\leq
\sup_{1\leq k\leq K} \langle v_{\hat\gamma,\gamma,k}, v_{\hat\gamma,\gamma,k}\rangle_2^2
 \nonumber \\ &\leq
\sup_{1\leq k\leq K} \sup_{u,v\in O} |\gamma(u,v)-\hat\gamma(u,v)|^2=
\mathcal{O}_p\left(r_{mn}^2\right).
\label{phitilde4a}
\end{align}
And since by Assumption (A7) $\sup_j\sup_{u\in O\cup M} |\tilde{\phi}^O_j(u)|\leq D_O<\infty$ the Cauchy-Schwarz inequality yields
\begin{align}
\sup_{u\in O\cup M}\sup_{1\leq k\leq K} \sum_{j=1}^K |\tilde{\phi}^O_j(u)\int_{O^2} (\hat\gamma(u,v)-\gamma(u,v))\hat\phi_k^{O}(u)\phi_j^{O}(v)dudv|=
\mathcal{O}_p\left(K^{1/2}r_{mn}\right),
\label{phitilde4}
\end{align}

By (c) of Theorem \ref{Theorem_UR}, \eqref{phitilde2a}, \eqref{phitilde4}, $\langle\phi_k^{O}-\phi_k^{O},\phi_j^{O}\rangle_2=0$ for $j\neq k$, and $\langle\phi_k^{O}-\phi_k^{O},\phi_k^{O}\rangle_2=1$ relation \eqref{phitilde2} can thus be rewritten in the form
\begin{align} \label{phitilde5}
\frac{\langle\hat\phi_k^{O},\gamma_u\rangle_2}{\hat\lambda_k^O}
&= \tilde{\phi}^O_k(u)+  \sum_{j=1}^{K}
 \tilde{\phi}^O_j(u)\langle\hat\phi_k^{O}-\phi_k^{O},\phi_j\rangle_2 \nonumber \\
 & \quad  +
\frac{1}{\lambda_k^O}\sum_{j=K+1}^\infty\lambda_j^O \tilde{\phi}^O_j(u)\langle\hat\phi_k^{O}-\phi_k^{O},\phi_j^{O}\rangle_2
+\mathcal{R}_{2,k}(u), \\
&\text{ where }
  \sup_{u\in O\cup M}\sup_{1\leq k\leq K} \lambda_k^O \mathcal{R}_{2,k}(u)=\mathcal{O}_p\left(K^{1/2}r_{mn}\right),
\nonumber
\end{align}
Result (d) of Theorem \ref{Theorem_UR} additionally implies
 \begin{align}
\sup_{1\leq k\leq K} (\delta_k^O)^2 \sum_{j=1}^\infty \langle\hat\phi_k^{O}-\phi_k^{O},\phi_j^{O}\rangle_2^2\leq
(\delta_k^O)^2\langle\hat\phi_k^{O}-\phi_k^{O},\hat\phi_k^{O}-\phi_k^{O}\rangle_2^2=
\mathcal{O}_p\left(r_{mn}^2\right),
\label{phitilde6a}
\end{align}
and, similar to \eqref{phitilde4}, the Cauchy-Schwarz inequality leads to
\begin{align}
\sup_{u\in O\cup M}\sup_{1\leq k\leq K}  \delta_k^O |\sum_{j=1}^K
 \tilde{\phi}^O_j(u)\langle\hat\phi_k^{O}-\phi_k^{O},\phi_j\rangle_2 |&=\mathcal{O}_p\left(K^{1/2}r_{mn}\right).
 \label{phitilde7}
\end{align}
By our assumptions on the sequence of eigenvalues we have for all $u\in O\cup M$
\begin{align*}
\sum_{j=K+1}^\infty (\lambda_j^O)^2  \tilde{\phi}^O_j(u)^2 &\leq D_O^2\sum_{j=K+1}^\infty (\lambda_j^O)^2\\
&=O(\sum_{j=K+1}^\infty j^{-2a_O})=O(K^{-2a_O+1})
=O(K(\lambda_K^O)^2)
\end{align*}
 When combining this result with \eqref{phitilde6a}, a further application of the Cauchy-Schwarz inequality yields
\begin{align}
\sup_{u\in O\cup M}&\sup_{1\leq k\leq K}| \delta_k^O\frac{1}{\lambda_k^O}\sum_{j=K+1}^\infty\lambda_j^O \tilde{\phi}^O_j(u)\langle\hat\phi_k^{O}-\phi_k^{O},\phi_j^{O}\rangle_2|\nonumber \\
&\leq  \delta_k^O \frac{1}{\lambda_k^O} (D_O^2 \sum_{K+1}^\infty (\lambda_j^O)^2)^{1/2}(\sum_{j=k+1}^\infty \langle\hat\phi_k^{O}-\phi_k^{O},\phi_j^{O}\rangle_2^2)^{1/2}
=\mathcal{O}_p\left(K^{1/2}r_{mn}\right).
\label{phitilde8}
\end{align}
Since  $\frac{\delta_k^O}{\lambda_k^O}\rightarrow 0$ as $k\rightarrow \infty$, the desired result is an immediate consequence of \eqref{phitilde1} - \eqref{phitilde8}.
\bigskip

A technical difficulty in the proof of Theorem \ref{Theorem_Main_Estim} consists in the fact that $\hat{\phi}^O_k$ and the observations $(Y_{ij},U_{ij})$ corresponding to the selected $i\in {1,\dots,n}$. But let $\hat\gamma_{-i}(t,s)\equiv\hat\gamma_{-i}(t,s,h_\gamma)$ denote the estimate of the covariance matrix when eliminating the $m$ observations $\{(Y_{ij},U_{ij})\}_{j=1,\dots,m}$  from the sample, and let $\hat{\lambda}_{k,-i}^O$ and $\hat{\phi}_{k,-i}^O$, $k=1,2,\dots$, denote eigenvalues and eigenfunctions of the corresponding covariance operator.
Although in our time series context there may still exist dependencies between $X_i$ and $\hat{\phi}_{k,-i}^O$, our assumptions imply that then all $\hat{\phi}_{k,-i}^O$ are independent of the particular samples $\{\epsilon_{ij}\}_{j=1,\dots,m}$ and $\{U_{ij}\}_{j=1,\dots,m}$. The following Lemma now provides bounds for the differences between $\hat{\tilde{\phi}}^O_{k}(u)=  \frac{\langle\hat\phi_k^{O},\hat\gamma_u\rangle_2}{\hat\lambda_k^O}$ and $\hat{\tilde{\phi}}^O_{k,-i}(u)=  \frac{\langle\hat\phi_{k,-i}^{O},\hat\gamma_{-i;u}\rangle_2}{\hat\lambda_{k,-i}^O}$, where $\hat\gamma_{-i;u}(v):=\hat\gamma_{-i}(u,v;h_\gamma)$.
\begin{lemma}\label{Lemma_phitildecv}
Under the assumptions of Theorem \ref{Theorem_Main_Estim} we have for all $1\leq k\leq K\leq \bar{K}_{mn}$
\begin{itemize}
\item[a)] $\sup_{(u,v)\in[a,b]^2}\Big|\hat\gamma(u,v;h_\gamma) - \hat\gamma_{-i}(u,v;h_\gamma)\Big|
=\mathcal{O}_p\left(\frac{1}{n^{1/2}r_{mn}}\right)$.
 \item[b)] $\sup_{k\leq K}
{\delta}_{k}^O\Vert \hat{\phi}_{k,-i}^O-\hat{\phi}_k^O\Vert_2=\mathcal{O}_p\left(\frac{1}{n^{1/2}r_{mn}}\right)$,
$\sup_{k\leq K}|\hat{\lambda}_{k,-i}^O-\hat\lambda_k^O|=\mathcal{O}_p\left(\frac{1}{n^{1/2}r_{mn}}\right)$
\item[c)]
$\sup_{u\in O}\sup_{1\leq k\leq K} \delta_k^O|\hat{\phi}_{k,-i}(u)-\hat{\phi}_k(u)|=
\mathcal{O}_p\left(\frac{K^{a_O+3/2}}{n^{1/2}r_{mn}}\right)=\mathcal{O}_p(\frac{1}{n^{1/2}})$
\end{itemize}
\end{lemma}

\noindent\textbf{Proof of Lemma \ref{Lemma_phitildecv}:}\
Based on the definitions and techniques introduced in the proof of Theorem \ref{Theorem_UR} and Lemma \ref{Lemma_UR_gamma} it is immediately seen that uniform rates of convergence of $\hat\gamma(u,v;h_\gamma) - \hat\gamma_{-i}(u,v;h_\gamma)$ can be derived by considering the following difference: 
\begin{align*}
\sup_{(u,v)\in[a,b]^2}&\Big|\Theta_{q,n\mathcal{M}}(u,v;h_\gamma) - \Theta_{q,n\mathcal{M}}^{-i}(u,v;h_\gamma)\Big|\leq\\
\sup_{(u,v)\in[a,b]^2}&\Big|\frac{1}{n\mathcal{M}h_\gamma^2}\sum_{j\neq l}\kappa\left(\frac{U_{ij}-u}{h_\gamma}\right)\kappa\left(\frac{U_{il}-v}{h_\gamma}\right)
\vartheta_q\left(U_{ij}-u,U_{i^\ast l}-u,C_{i^\ast jl}\right)\Big|+\\
\sup_{(u,v)\in[a,b]^2}&\Big|\frac{1}{n(n-1)}\frac{1}{\mathcal{M}h_\gamma^2}\sum_{k\neq i,j\neq l}\kappa\left(\frac{U_{kj}-u}{h_\gamma}\right)\kappa\left(\frac{U_{kl}-v}{h_\gamma}\right)\vartheta_q
\left(U_{kj}-u,U_{kl}-u,C_{kjl}\right)\Big|.
\end{align*}
Using similar arguments as in the proof of Lemma
\ref{Lemma_UR}, leads to
\begin{align*}
\sup_{(u,v)\in[a,b]^2}\Big|\Theta_{q,n\mathcal{M}}(u,v;h_\gamma) - \Theta_{q,n\mathcal{M}}^{-i}(u,v;h_\gamma)\Big|=\mathcal{O}_p\left(\frac{1}{\sqrt{n^2\mathcal{M}h_\gamma^2}} + \frac{1}{\sqrt{n^2}}\right)
\end{align*}
which implies
\begin{align*}
\sup_{(u,v)\in[a,b]^2}\Big|\hat\gamma(u,v;h_\gamma) - \hat\gamma_{-i}(u,v;h_\gamma)\Big|=\mathcal{O}_p\left(\frac{1}{\sqrt{n^2\mathcal{M}h_\gamma^2}} + \frac{1}{\sqrt{n^2}}\right),
\end{align*}
and assertion a) of the Lemma is an immediate consequence.

The inequalities used to prove (c) and (d) of Theorem \ref{Theorem_UR} now lead to
$\sup_{k\geq 1}|\hat{\lambda}_{k,-i}^O-\hat\lambda_k^O|=\mathcal{O}_p\left(n^{-1/2}r_{mn}\right)$ and
$\sup_{1\leq k\leq K} \hat{\delta}_{k,-i}^O\Vert \hat{\phi}_{k,-i}^O-\hat{\phi}_k^O\Vert_2=\mathcal{O}_p\left(n^{-1/2}r_{mn}\right)$, where
$\hat{\delta}_{k,-i}^O:=\min_{ j\neq k}\{\hat{\lambda}_{j}^O-\hat{\lambda}_{k}^O\}$.
By (c) of Theorem \ref{Theorem_UR}, our assumptions on $\lambda_j$, and $k\leq K\leq \bar{K}_{mn}$ with
$\bar{K}_{mn}^{a_O+3/2}r_{mn}=O(1)$ we additionally have $\hat{\delta}_{k,-i}^O(\delta_{k,-i}^O)^{-1}=1+o_p(1)$. This proves assertion b) of the Lemma.

Furthermore, we have $\hat{\phi}_{k}^O(u)=\frac{\langle\hat{\phi}_{k}^O,\hat\gamma_{u}\rangle_2}{\hat\lambda_{k}^O}$ as well as
$\hat{\phi}^O_{k,-i}(u)= \frac{\langle\hat{\phi}_{k,-i}^O,\hat\gamma_{-i;u}\rangle_2}{\hat\lambda_{k,-i}^O}$ for $u\in O$.
The difference can be rewritten in the form
\begin{align}
\hat{\phi}_{k}^O(u)-\hat{\phi}^O_{k,-i}(u)=&
\frac{\langle\hat\phi_k^{O}-\hat{\phi}_{k,-i}^O,\gamma_u\rangle_2}{\hat\lambda_k^O}
+\frac{\langle\hat\phi_k^{O}-\hat{\phi}_{k,-i}^O,\hat\gamma_u-\gamma_u\rangle_2}{\hat\lambda_k^O}
\nonumber \\
&+\frac{\langle\hat{\phi}_{k,-i}^O,\hat\gamma_u-\hat\gamma_{-i;u}\rangle_2}{\hat\lambda_{k}^O}
+\frac{(\hat\lambda_{k,-i}^O-\hat\lambda_k^O)\langle\hat{\phi}_{k,-i}^O,\hat\gamma_{-i;u}\rangle_2}
{\hat\lambda_k^O\hat\lambda_{k,-i}^O}.
\label{phitildecv2}
\end{align}
When analyzing the terms in \eqref{phitildecv2} first note that by assertions a) and b) of the lemma, and by (c) of Theorem \ref{Theorem_UR}  \\
 \begin{align}
&\sup_{u\in O}\sup_{1\leq k\leq K}\left|\frac{\langle\hat{\phi}_{k,-i}^O,\hat\gamma_u-\hat\gamma_{-i;u}\rangle_2}{\hat\lambda_{k,-i}^O}\right|
=\mathcal{O}_p\left(\frac{1}{\lambda_K^O}\frac{r_{mn}}{n^{1/2}}\right)
=\mathcal{O}_p\left(\frac{K^{a_O}r_{mn}}{n^{1/2}}\right)\nonumber \\
&\sup_{u\in O}\sup_{1\leq k\leq K}\left|\frac{\langle\hat\phi_k^{O}-\hat{\phi}_{k,-i}^O,\hat\gamma_u-\gamma_u\rangle_2}{\hat\lambda_k^O}\right|
=\mathcal{O}_p\left(\frac{K^{2a_O+1}r_{mn}^2}{n^{1/2}}\right)=o_p\left(\frac{K^{a_O}r_{mn}}{n^{1/2}}\right)
\label{phitildecv3}
\end{align}
The  differences between the eigenfunctions  $\hat{\phi}_{k,-i}^O(u)$ and $\hat{\phi}_{k}^O(u)$ reflect the elimination of one single curve, and it is immediately clear that the convergence results of Theorem \ref{Theorem_UR} and all arguments of Lemma \ref{Lemma_phitilde} remain valid when considering estimated covariances $\hat\gamma_{-i}(u,v)$, eigenvalues $\hat{\lambda}_{k,-i}^O$, and eigenfunctions $\hat{\phi}_{k,-i}^O(u)$ of the reduced sample. It thus follows from Lemma \ref{Lemma_phitilde} and our assumption on $K\leq \bar{K}_{mn}$  that $\frac{\langle\hat{\phi}_{k,-i}^O,\hat\gamma_{-i;u}\rangle_2}
{\hat\lambda_{k,-i}^O}$ is asymptotically uniformly bounded over all $u\in O\cup M$ and $k\le K$. Hence,
\begin{align}
&\sup_{u\in O\cup M}\sup_{1\leq k\leq K}\left|\frac{(\hat\lambda_{k,-i}^O-\hat\lambda_k^O)\langle\hat{\phi}_{k,-i}^O,\hat\gamma_{-i;u}\rangle_2}
{\hat\lambda_k^O\hat\lambda_{k,-i}^O}\right|
=\mathcal{O}_p\left(\frac{1}{\lambda_K^O}\frac{r_{mn}}{n^{1/2}}\right)
=\mathcal{O}_p\left(\frac{K^{a_O}r_{mn}}{n^{1/2}}\right)
\label{phitildecv4}
\end{align}
The first term on the right side of \eqref{phitildecv2} can now be analyzed by generalizing the arguments of
Lemma \ref{Lemma_phitilde}. Similar to  \eqref{phitilde2} we obtain
 \begin{align*}
\frac{\langle\hat\phi_k^{O}-\hat{\phi}_{k,-i}^O,\gamma_u\rangle_2}{\hat\lambda_k^O}= \frac{1}{\hat\lambda_k^O} \sum_{j=1}^\infty
\lambda_j^O \tilde{\phi}^O_j(u)\langle\hat\phi_k^{O}-\hat{\phi}_{k,-i}^O,\phi_j^{O}\rangle_2,
\end{align*}
while \eqref{phitilde2a} becomes
\begin{align*}
&\lambda_j^O \langle\hat\phi_k^{O}-\hat{\phi}_{k,-i}^O,\phi_j^{O}\rangle_2
=(\hat\lambda_k^O-\hat\lambda_{k,-i}^O)\langle\hat\phi_k^{O}-\hat{\phi}_{k,-i}^O,\phi_j^{O}\rangle_2
\nonumber \\
&+ \int_{O^2} (\hat\gamma(u,v)-\gamma(u,v))(\hat\phi_k^{O}(u)-\hat{\phi}_{k,-i}^O(u))\phi_j^{O}(v)dudv
\nonumber \\
&+\int_{O^2} (\hat\gamma(u,v)-\hat\gamma_{-i}(u,v))\hat{\phi}_{k,-i}^O(u)\phi_j^{O}(v)dudv
\end{align*}
Using result b),  a straightforward generalization of the arguments given by \eqref{phitilde4a} and  \eqref{phitilde4} then leads to
\begin{align*}
\frac{\langle\hat\phi_k^{O}-\hat{\phi}_{k,-i}^O,\gamma_u\rangle_2}{\hat\lambda_k^O}
&= \frac{(\hat\lambda_k^O-\hat\lambda_{k,-i}^O)}{\hat\lambda_k^O}  \sum_{j=1}^{K}
 \tilde{\phi}^O_j(u)\langle\hat\phi_k^{O}-\hat{\phi}_{k,-i}^O,\phi_j\rangle_2 \nonumber \\
 & \quad  +
\frac{1}{\lambda_k^O}\sum_{j=K+1}^\infty\lambda_j^O \tilde{\phi}^O_j(u)\langle\hat\phi_k^{O}-\hat{\phi}_{k,-i}^O,\phi_j^{O}\rangle_2
+\mathcal{R}_{1,k}^{(-i)}(u), \\
&\text{ where }
  \sup_{u\in O\cup M}\sup_{1\leq k\leq K}  \mathcal{R}_{1,k}^{(-i)}(u)=\mathcal{O}_p\left(\frac{K^{a_O}r_{mn}}{n^{1/2}}\right),
\end{align*}
and proceeding similar to \eqref{phitilde6a} -  \eqref{phitilde8} we can conclude that
\begin{align}
&\sup_{u\in O\cup M}\sup_{1\leq k\leq K}\left|\frac{\langle\hat\phi_k^{O}-\hat{\phi}_{k,-i}^O,\gamma_u\rangle_2}{\hat\lambda_k^O}\right|
=\mathcal{O}_p\left(\frac{K^{a_O+3/2}r_{mn}}{n^{1/2}}\right)
\label{phitildecv5}
\end{align}
Assertion c) of the lemma is now an immediate consequence of \eqref{phitildecv2} - \eqref{phitildecv5}.

\bigskip

\noindent\textbf{Proof of Theorem \ref{Theorem_Main_Estim}.}\\
We have to consider the asymptotic behavior of
\begin{align*}
\widehat{\mathcal{L}}_{K}(\mathbb{X}_i^O)(u)&=\hat{\mu}(u;h_\mu)
+\sum_{k=1}^{K}\hat{\xi}_{ik}^O\hat{\tilde{\phi}}^O_k(u),\quad u\in O\cup  M.
\end{align*}
Rates of convergence of $\hat{\mu}(u;h_\mu)$ are given by Theorem \ref{Theorem_UR} (a), while Lemma \ref{Lemma_phitilde} provides rates of convergence for $\hat{\tilde{\phi}}^O_k(u)$. We therefore additionally have to consider convergence of the PC scores
\begin{align*}
&\hat{\xi}^O_{ik}=\sum_{j=2}^{m}\hat{\phi}_k^O(U_{i(j)})
(Y_{i(j)}-\hat{\mu}(U_{i(j)};h_\mu))(U_{i(j)}-U_{i,(j-1)}),
\end{align*}
where $U_{i(j)}$, $j=1,\dots,m$ is the order sample of observation points. By our assumption on $h_\mu$ result (a) of Theorem \ref{Theorem_UR} directly implies that
\begin{align}
\hat{\xi}^O_{ik}&= \sum_{j=2}^{m}\hat{\phi}_k^O(U_{i(j)})
(Y_{i(j)}-\mu(U_{i(j)}))(U_{i(j)}-U_{i,(j-1)})+\mathcal{R}_{1,k,i},
\nonumber \\
&\text{where } \sup_{1\leq k\leq K} |\mathcal{R}_{1,k,i}| =\mathcal{O}_p\left((mn)^{-2/5}\right)
\label{TME1}
\end{align}
A technical difficulty of subsequent analysis consists in the fact that $\hat{\phi}_k^O$ and $\{(Y_{ij},U_{ij})\}_{j=1,\dots,m}$ are correlated. As defined above, we thus eliminate the $m$ observations representing curve $X_i$ and consider the eigenfunction $\hat{\phi}^O_{k,-i}(u)$ of the reduced sample. We can then infer from  result c) of Lemma \ref{Lemma_phitildecv} that
\begin{align}
\mathcal{R}_{2,k,i}&:= |\sum_{j=2}^{m}(\hat{\phi}_{k,-i}^O(U_{i(j)})-\hat{\phi}_{k}^O(U_{i(j)}))
(Y_{i(j)}-\mu(U_{i(j)}))(U_{i(j)}-U_{i,(j-1)})|\nonumber \\
&\text{satisfies } \sup_{1\leq k\leq K} |\mathcal{R}_{2,k,i}| =\mathcal{O}_p(n^{-1/2})
\label{TME2}
\end{align}
Recall that $Y_{ij}=X_i(U_{ij})+\epsilon_{ij}$ and that $\hat{\phi}_{k,-i}^O$ is independent of $\epsilon_{ij}, U_{ij}$. Since $k\le K\leq \bar{K}_{mn}$, result c) of Lemma \ref{Lemma_phitildecv} also implies that $\hat{\phi}_{k,-i}^O$, $k\le K$, are asymptotically uniformly bounded over all $u\in O\cup M$ and $k\le K$. By our assumptions on the error term we can now immediately infer that
\begin{align}
 &\sum_{j=2}^{m}\hat{\phi}_{k,-i}^O(U_{i(j)})
(Y_{i(j)}-\mu(U_{i(j)})(U_{i(j)}-U_{i,(j-1)}) \nonumber \\
&=\sum_{j=2}^{m}\hat{\phi}_{k,-i}^O(U_{i(j)})
(X_i(U_{i(j)})-\mu(U_{i(j)}))(U_{i(j)}-U_{i,(j-1)})+\mathcal{R}_{3,k,i},
\nonumber \\
& \quad \text{ with } \E\left(\mathcal{R}_{3,k,i}^2|\ \hat\gamma_{-i}\right)\le D_{k,i} \frac{1}{m}
\text{ and }  \sup_{1\leq k\leq K} |D_{k,i}|=\mathcal{O}_p(1)
\label{TME3}
\end{align}

Let $X_i^*:=X_i-\mu$, and let $F_{U|O}$ denote the distribution function of $U_{ij}$.
It is then  well-known that the random variables $V_{ij}:=F_{U|O}(U_{ij})$ are $U(0,1)$-distributed, and
$\E(V_{i(j)})-V_{i(j-1)}))=\frac{1}{m+1}$ while $\V(V_{i(j)})-V_{i(j-1)}))=\frac{m}{(m+1)^2(m+2)}$. By our assumptions the density $f_{U|O}$ ($O\equiv O_i$) of $U_{ij}$ a Taylor expansion now yields  \\
$F_{U|O}^{-1}(V_{i(j)})-F_{U|O}^{-1}(V_{i(j-1)})=\frac{1}{mf_{U|O}(F_{U|O}^{-1}(V_{i(j)}))}+\mathcal{R}_{4,k,i}^*$ with
$\E(\mathcal{R}_{4,k,i}^*)=\mathcal{O}_p(1/m^2)$. This implies
\begin{align*}
 &\sum_{j=2}^{m}\hat{\phi}_{k,-i}^O(U_{i(j)})
X_i^*(U_{i(j)})(U_{i(j)}-U_{i,(j-1)})\\
&=\sum_{j=2}^{m}\hat{\phi}_{k,-i}^O(F^{-1}(V_{i(j)}))X_i^*(F^{-1}(V_{i(j)}))
(F^{-1}(V_{i(j)})-F^{-1}(V_{i(j-1)}))\\
&= \sum_{j=2}^{m}\hat{\phi}_{k,-i}^O(U_{i(j)}))X_i^*(U_{i(j)})\frac{1}{m f(U_{i(j)})}
+\mathcal{R}_{4,k,i}^{**}
\end{align*}
with $ \E\left(\mathcal{R}_{4,k,i}^{**}|\ \hat\gamma_{-i}\right)\le D_{ik}^* \frac{1}{m}$
for some
$D_{ik}^*<\infty$ satisfying $\sup_{1\leq k\leq K} |D_{k,i}^*|=\mathcal{O}_p(1)$. Obviously,
\begin{align*}
&\E\left(\sum_{j=2}^{m}\hat{\phi}_{k,-i}^O(U_{i(j)}))X_i^*(U_{i(j)})\frac{1}{m f(U_{i(j)})}|\ \hat{\phi}_{k,-i}^O,X_i\right)=\frac{m-1}{m} \int_O \hat{\phi}_{k,-i}^O(u)
X_i^*(u)du,
\end{align*}
and independent of $k$ the conditional variance of this random variable can be bounded by $1/m$. We therefore arrive at
\begin{align}
 &\sum_{j=2}^{m}\hat{\phi}_{k,-i}^O(U_{i(j)}))
X_i^*(U_{i(j)})(U_{ij}-U_{i,j-1})
=\int_O \hat{\phi}_{k,-i}^O(u)
X_i^*(u)du+\mathcal{R}_{4,k,i},
\nonumber \\
& \text{ with } \E\left(\mathcal{R}_{4,k,i}^2|\ \hat\gamma_{-i}\right)\le D_{ik}^{**} \frac{1}{m},
\text{ and }  \sup_{1\leq k\leq K} |D_{k,i}^{**}|=\mathcal{O}_p(1).
\label{TME4}
\end{align}

Additionally note that Lemma \ref{Lemma_phitilde} together with $K^{a_O+3/2}r_{mn}=\mathcal{O}(1)$ implies that
$\sup_{u\in O\cup M}\sup_{1\leq k\leq K} \hat{\tilde{\phi}}^O_k(u)=\mathcal{O}(1)$, while
(a) of Theorem \ref{Theorem_UR} yields $\sup_{u\in O\cup M}|\mu(u)-\hat\mu(u)|=\mathcal{O}_p((mn)^{-2/5})$ We can therefore infer from \eqref{TME1} - \eqref{TME4} that for $u\in O\cup M$
\begin{align}
\widehat{\mathcal{L}}_{K}(\mathbb{X}_i^O)(u)=\mu(u)+\sum_{k=1}^K (\int_O \hat{\phi}_{k,-i}^O(u)
X_i^*(u)du)\hat{\tilde{\phi}}^O_k(u)+\mathcal{O}_p\left(\frac{K}{m^{1/2}}+\frac{K}{n^{1/2}}\right).
\label{TME4a}
\end{align}

Since $\mathcal{L}_{K}(X_i^O)(u)=\mu(u) + \sum_{k=1}^K \xi^O_k\tilde{\phi}_j^O(u)$ the next step is to consider the errors $\xi^O_{ik}-\int_O \hat{\phi}_{k,-i}^O(u)X_i^*(u)du$ for $k=1,\dots,K$.

The differences between the eigenfunctions $\hat{\phi}_{k,-i}^O(u)$ and $\hat{\phi}_{k}^O(u)$ reflect the elimination of one single curve, and it is immediately clear that the convergence results of Theorem \ref{Theorem_UR} and all arguments of Lemma \ref{Lemma_phitilde} remain valid when considering estimated covariances $\hat\gamma_{-i}(u,v)$, eigenvalues $\hat{\lambda}_{k,-i}^O$, and eigenfunctions $\hat{\phi}_{k,-i}^O(u)$ of the reduced sample. Furthermore, note that since
$\Vert \hat{\phi}_{k,-i}^O\Vert_2=\Vert \phi_j^{O}\Vert_2=1$ implies $\langle\hat\phi_{k,-i}^{O},\phi_k^{O}\rangle_2=1-\frac{1}{2}\Vert \hat\phi_{k,-i}^{O}-\phi_k^{O}\Vert_2^2$, and recall the Karhunen-Lo\`eve decomposition $X_i^*(u)=\sum_{j=1}^\infty \xi^O_{ik}
\phi^O_{k}(u)$.

 We can thus infer from \eqref{phitilde2a} in Lemma \ref{Lemma_phitilde} that
\begin{align}
&\int_O \hat{\phi}_{k,-i}^OX_i^*(u)(u)=\sum_{j=1}^\infty
\xi^O_{ij} \langle\hat{\phi}_{k,-i}^O,\phi_j^{O}\rangle_2=\xi^O_{ik} \langle\hat{\phi}_{k,-i}^O,\phi_k^{O}\rangle_2+\sum_{j\ne k}
\xi^O_{ij} \langle\hat{\phi}_{k,-i}^O,\phi_j^{O}\rangle_2
 \nonumber\\
&=\xi^O_{ik}+\sum_{j=1}^{k-1}
\xi^O_{ij} \frac{\hat{\lambda}_{k,-i}^O}{\lambda_j^O}\langle\hat{\phi}_{k,-i}^O-\phi_k^{O},\phi_j^{O}\rangle_2
+\sum_{j=k+1}^\infty
\xi^O_{ik} \langle\hat{\phi}_{k,-i}^O-\phi_k^{O},\phi_j^{O}\rangle_2 \nonumber\\
&+\sum_{j=1}^{k-1} \xi^O_{ik}\frac{1}{\lambda_j^O}\int_{O^2} (\hat\gamma_{-i}(u,v)-\gamma(u,v))\hat\phi_k^{O}(u)\phi_j^{O}(v)dudv
-\frac{\xi^O_{ik}}{2}\Vert \hat\phi_{k,-i}^{O}-\phi_k^{O}\Vert_2^2.
\label{TME5b}
\end{align}
 By our assumptions on the sequence of eigenvalues we have
{\small \begin{align}
Q_{i,k}:=\sum_{j=1}^{k-1}\frac{(\xi^O_{ik})^2}{(\lambda_j^O)^2}=\mathcal{O}_p\left(\sum_{j=1}^{k-1}
\frac{\E(( \xi^O_{ik})^2)}{(\lambda_j^O)^2}\right)=
\mathcal{O}_p\left(\sum_{j=1}^{k-1}\frac{1}{\lambda_j^O} \right)=\mathcal{O}_p\left(\sum_{j=1}^{k-1} j^{a_O}\right)=\mathcal{O}_p(k^{a_O+1})
\label{TME5c}
\end{align}}
Using the Cauchy-Schwary inequality, \eqref{TME5c}, inequality \eqref{phitilde6a} in Lemma \ref{Lemma_phitilde}, and (d) of Theorem \ref{Theorem_UR} lead to
\begin{align}
&|\sum_{j=1}^{k-1}
\xi^O_{ij} \frac{\hat{\lambda}_{k,-i}^O}{\lambda_j^O}\langle\hat{\phi}_{k,-i}^O-\phi_k^{O},\phi_j^{O}\rangle_2|
\leq \frac{\lambda_k^OQ_{i,k}^{1/2}}{\delta_k^O} R_{5,k,i}, \nonumber\\
& \qquad\text{ where } R_{5,k,i}\ge 0
\text{ with } \sup_{1\leq k\leq K}R_{5,k,i} =\mathcal{O}_p(r_{mn})
\label{TME5d}
\end{align}
for all $1\leq k\leq K$. Hence,
{\small  \begin{align}
\sum_{k=1}^{K} |\hat{\tilde{\phi}}^O_k(u)| |\sum_{j=1}^{k-1}\xi^O_{ij} \frac{\hat{\lambda}_{k,-i}^O}{\lambda_j^O}\langle\hat{\phi}_{k,-i}^O-\phi_k^{O},\phi_j^{O}\rangle_2|
\le & \left(\sum_{k=1}^{K} \frac{(\lambda_k^O)^2 Q_{i,k}}{(\delta_k^O)^2}\right)^{1/2}
\left(\sum_{j=k}^{K} R_{5,k,i}^2 \hat{\tilde{\phi}}^O_j(u)^2\right)^{1/2}\nonumber\\
&=\mathcal{O}_p\left(K^{a_O/2+5/2}r_{mn}\right)
\label{TME5e}
\end{align}}

Since $\E(Q_{i,k}^*):=\E(\sum_{j=k+1}^\infty (\xi^O_{ij})^2)=\sum_{j=k+1}^\infty \lambda_j^O=O(k^{-a_O+1})$,
similar arguments based on inequalities \eqref{phitilde6a} and \eqref{phitilde4a} in Lemma \ref{Lemma_phitilde}
yield \\
$\sum_{k=1}^{K} |\hat{\tilde{\phi}}^O_k(u)| |\sum_{j=k+1}^\infty
\xi^O_{ij} \langle\hat{\phi}_{k,-i}^O-\phi_k^{O},\phi_j^{O}\rangle_2|=\mathcal{O}_p\left(K^{a_O/2+5/2}r_{mn}\right)$ as well as \\
$\sum_{k=1}^{K} |\hat{\tilde{\phi}}^O_k(u)| |\sum_{j=1}^{k-1} \xi^O_{ij}\frac{1}{\lambda_j^O}\int_{O^2} (\hat\gamma_{-i}(u,v)-\gamma(u,v))\hat\phi_k^{O}(u)\phi_j^{O}(v)dudv|\allowbreak
=\mathcal{O}_p\left(K^{a_O/2+3/2}r_{mn}\right)$,\\
while by (d) of Theorem \ref{Theorem_UR} $\sum_{k=1}^{K} |\hat{\tilde{\phi}}^O_k(u)| \frac{|\xi^O_{ik}|}{2}\Vert \hat\phi_{k,-i}^{O}-\phi_k^{O}\Vert_2^2 =\mathcal{O}_p\left(K^{3a_O/2+3}r_{mn}^2\right)$.
Together with \eqref{TME4a} and $K^{a_O+3/2}r_{mn}=\mathcal{O}(1)$ we therefore arrive at
\begin{align}
&\widehat{\mathcal{L}}_{K}(\mathbb{X}_i^O)(u)=\mu(u)+\sum_{k=1}^K \xi^O_{ik}\hat{\tilde{\phi}}^O_k(u)+\mathcal{O}_p\left(K\left(\frac{1}{m^{1/2}}+K^{a_O/2+3/2}r_{mn}\right)\right)
\nonumber \\
&= \mathcal{L}_{K}(X_i^O) +\sum_{k=1}^K \xi^O_{ik}(\hat{\tilde{\phi}}^O_k(u)-\tilde{\phi}^O_k(u))+
\mathcal{O}_p\left(K\left(\frac{1}{m^{1/2}}+K^{a_O/2+3/2}r_{mn}\right)\right)
\label{TME6}
\end{align}
By Lemma \ref{Lemma_phitilde} we have $\sum_{k=1}^K (\delta_k^O)^2(\hat{\tilde{\phi}}^O_k(u)-\tilde{\phi}^O_k(u))^2
=\mathcal{O}(K^2r_{mn}^2)$, and the Cauchy-Schwarz inequality implies
\begin{align*}
\sum_{k=1}^K \xi^O_{ik}(\hat{\tilde{\phi}}^O_k(u)-\tilde{\phi}^O_k(u))&=
\mathcal{O}_p\left( \left(\sum_{k=1}^K  \frac{\E((\xi^O_{ik})^2)}{(\delta_k^O)^2} \right)^{1/2}
\left(\sum_{k=1}^K (\delta_k^O)^2(\hat{\tilde{\phi}}^O_k(u)-\tilde{\phi}^O_k(u))^2\right)^{1/2}
\right)\nonumber \\
&=\mathcal{O}_p\left( \left( \sum_{k=1}^K  k^{a_O+2} \right)^{1/2}
\left(K^2r_{mn}^2\right)^{1/2}
\right)=\mathcal{O}_p\left(K^{a_O/2+5/2}r_{mn}\right)
\end{align*}
and therefore
\begin{align}
\widehat{\mathcal{L}}_{K}(\mathbb{X}_i^O)(u)=\mathcal{L}_{K}(X_i^O)(u)+
\mathcal{O}_p\left(K\left(\frac{1}{m^{1/2}}+K^{a_O/2+3/2}r_{mn}\right)\right)
\label{TME7}
\end{align}

We finally have to consider the truncation error.
Recall that it is assumed that there is a constant $D_O<\infty$ such that
$\sup_{u\in O\cup M} \sup_{k\geq 1}  |\tilde{\phi}^O_k(u)|\leq D_O $ Since $\mathcal{L}_{K}(X_i^O)(u)=\mu(u)+\sum_{k=1}^{\infty}\xi_{ik}^O\tilde{\phi}^O_k(u)$ we have
\begin{align}\label{TME8}
\E\left(\bigl(\mathcal{L}_{K}(X_i^O)(u)-\mu(u)-\sum_{k=1}^{K}\xi_{ik}^O\tilde{\phi}^O_k(u)\bigr)^2\right)
=\sum_{k=K+1}^{\infty}\lambda_k \tilde{\phi}^O_k(u)^2\leq D_O^2 \sum_{k=K+1}^{\infty}\lambda_k^O\nonumber \\
=\mathcal{O}\left(\sum_{k=K+1}^{\infty} k^{-a_O} \right) =\mathcal{O}\left(K^{-a_O+1} \right)
\end{align}
Result \eqref{AME1} now follows from \eqref{TME7} and \eqref{TME8}. When additionally noting that  standard arguments imply that the local linear estimator of $X_i$  with bandwidth $h_X\asymp m^{-1/5}$ satisfies
 $|\widehat{X}_i^{O}(\vartheta_u;h_X)-X_i(\vartheta_u)|=\mathcal{O}\left(m^{-2/5}\right)$, result \eqref{AME2} follows from \eqref{AME1} and definition of
$\widehat{\mathcal{L}}^*_{K}(\mathbb{X}_i^O)$.

\newpage

\noindent\textbf{Proof of Proposition \ref{pro:IPA}:}\\
For $u\in O_2$ the optimal linear reconstruction of $X_i^{O_2}(u)$, given $X_i^{O_1}$, is $\mathcal{L}(X_i^{O_1})(u)$, such that $X_i^{O_2}(u)=\mathcal{L}(X_i^{O_1})(u)+\mathcal{Z}_i^{O_2}(u)$, with $\tilde{X}_i^{O_2}(u)=\mathcal{L}(X_i^{O_1})(u)$. From result (a) of Theorem \ref{OptimalPrediction}, we know that $\mathcal{Z}_i^{O_2}(u)$ and $X_i^{O_1}(v)$ are uncorrelated for all $u\in O_2$ and all $v\in O_1$. Consequently, by linearity of $\mathcal{L}$, also
\begin{equation*}
\E(\mathcal{L}(X_i^{O_1})(u)\,\mathcal{L}(\mathcal{Z}_i^{O_2})(u))=0\quad\text{for all}\quad u\in O_2,
\end{equation*}
where in $\mathcal{L}(\mathcal{Z}_i^{O_2})(u)$ we are using that $\mathcal{L}$ is also well defined as a linear mapping from $\mathbb{L}^2(O_2)$ to $O_2$; see remark to  \eqref{eq:PEFexp}.

Therefore,
\begin{align*}
&0\leq
\E\left(\left(X_i^{M_2}(u)-\mathcal{L}(X_i^{O_1})(u)-\mathcal{L}(\mathcal{Z}_i^{O_2})(u)\right)^2\right)=\\
&
=\E\left(\left(X_i^{M_2}(u)-\mathcal{L}(X_i^{O_1})(u)\right)^2\right)
-2\E\left(X_i^{M_2}(u)\mathcal{L}(\mathcal{Z}_i^{O_2})(u)\right)
+\E\left(\left(\mathcal{L}(\mathcal{Z}_i^{O_2})(u)\right)^2\right)\\
&\Rightarrow
-\E\left(\left(X_i^{M_2}(u)-\mathcal{L}(X_i^{O_1})(u)\right)^2\right)\leq\\
&\leq - 2\E\left(X_i^{M_2}(u)\mathcal{L}(\mathcal{Z}_i^{O_2})(u)\right)
+\E\left(\left(\mathcal{L}(\mathcal{Z}_i^{O_2})(u)\right)^2\right).
\end{align*}
But then,
\begin{align*}
&\E\left(\left(X_i^{M_2}(u)-\mathcal{L}(\mathcal{Z}_i^{O_2})(u)\right)^2\right)=\\
&=\E\left(\left(X_i^{M_2}(u)\right)^2\right)
-2\E\left(X_i^{M_2}(u)\mathcal{L}(\mathcal{Z}_i^{O_2})(u)\right)
+\E\left(\left(\mathcal{L}(\mathcal{Z}_i^{O_2})(u)\right)^2\right)\\
&\geq \E\left(\left(X_i^{M_2}(u)\right)^2\right) -
\E\left(\left(X_i^{M_2}(u)-\mathcal{L}(X_i^{O_1})(u)\right)^2\right)
\end{align*}

On the other hand, we have that $\E(\mathcal{L}(\tilde{X}_i^{O_2})(u)\mathcal{L}(\mathcal{Z}_i^{O_2})(u))=0$ for all $u\in M_2$, which follows by the same reasoning as used above, since $\tilde{X}_i^{O_2}(u)=\mathcal{L}(X_i^{O_1})(u)$, with $u\in O_2$, is just another linear transformation of $X_i^{O_1}$ and $X_i^{O_1}(v)$ is known to be uncorrelated with $\mathcal{Z}_i^{O_2}(u)$ for all $u\in O_2$ and $v\in O_1$ by result (a) of Theorem \ref{OptimalPrediction}. Therefore, using also the latter inequality,
\begin{align*}
&\E\left(\left(X_i^{M_2}(u)-\mathcal{L}(X_i^{O_2})(u)\right)^2\right)=
\E\left(\left(X_i^{M_2}(u)-\mathcal{L}(\tilde{X}_i^{O_2})(u)-\mathcal{L}(\mathcal{Z}_i^{O_2})(u)\right)^2\right)=\\
&\E\left(\left(X_i^{M_2}(u)-\mathcal{L}(\mathcal{Z}_i^{O_2})(u)\right)^2\right)
-2\E\left(X_i^{M_2}(u)\,\mathcal{L}(\tilde{X}_i^{O_2})(u)\right)
+\E\left(\left(\mathcal{L}(\tilde{X}_i^{O_2})(u)\right)^2\right)\\
&\geq
\E\left(\left(X_i^{M_2}(u)\right)^2\right) -
\E\left(\left(X_i^{M_2}(u)-\mathcal{L}(X_i^{O_1})(u)\right)^2\right)\\
&-2\E\left(X_i^{M_2}(u)\,\mathcal{L}(\tilde{X}_i^{O_2})(u)\right)
+\E\left(\left(\mathcal{L}(\tilde{X}_i^{O_2})(u)\right)^2\right)=\\
&=\E\left(\left(X_i^{M_2}(u)-\mathcal{L}(\tilde{X}_i^{O_2})(u)\right)^2\right)-
\E\left(\left(X_i^{M_2}(u)-\mathcal{L}(X_i^{O_1})(u)\right)^2\right)\\[1ex]
&\Rightarrow
\E\left(\left(X_i^{M_2}(u)-\mathcal{L}(\tilde{X}_i^{O_2})(u)\right)^2\right)\leq \\
&\leq\E\left(\left(X_i^{M_2}(u)-\mathcal{L}(X_i^{O_1})(u)\right)^2\right)
+\E\left(\left(X_i^{M_2}(u)-\mathcal{L}(X_i^{O_2})(u)\right)^2\right)
\end{align*}

\newpage

\section{Further explanations}\label{appendix:FE}
\subsection{Discontinuity of $\mathcal{L}$}\label{appendix:FE1} 
In the second footnote in the introduction we claim that the optimal linear functional $\mathcal{L}(X_i^O)(u)$ may not be a continuous functional $\mathbb{L}^2(O)\rightarrow \mathbb{R}$. This possible discontinuity may occur even though we are considering continuous functions $X_i$. In order to clarify this, we give here an example where a small $\mathbb{L}^2$-distance $||f-g||_2$ goes along with a very large pointwise distance $|f(\vartheta)-g(\vartheta)|$ using the additional requirement that $f$ and $g$ are both absolute continuous functions.

Consider the functional  $\mathcal{L}:\mathbb{L}^2( O)\rightarrow\mathbb{R}$ defined by the point evaluation $\mathcal{L}(f)=f(\vartheta)$ for some $\vartheta\in O$ with $O\subset\mathbb{R}$. This functional is discontinuous (and unbounded). Smoothness does not help, since $\mathcal{L}$ remains discontinuous (and unbounded) when restricting attention to the subclass of all functions $f\in C^\infty(O)\subset \mathbb{L}^2( O)$ with infinitely many derivatives.

This is easily seen by the following construction: Let $f\in C^\infty(O)$ possess infinitely many derivatives, and for $\sigma>0$ consider the functions
$$
g_\sigma(u):=f(u)+\Big(\frac{1}{2\pi \sigma}\Big)^{1/4}\exp\Big(-\frac{(u-\vartheta)^2}{4\sigma^2}\Big),\quad u\in O
$$
Obviously, $g_\sigma(u)\in C^\infty(O)$ for all $\sigma>0$. Moreover,
\begin{align*}
\Vert f-g_\sigma\Vert_2&=\left(\int_O \Big(\frac{1}{2\pi \sigma}\Big)^{1/2}\exp\Big(-\frac{(u-\vartheta)^2}{2\sigma^2}\Big)du\right)^{1/2}\\
& \leq \left(\int_{-\infty}^\infty \Big(\frac{1}{2\pi \sigma}\Big)^{1/2}\exp\Big(-\frac{(u-\vartheta)^2}{2\sigma^2}\Big)du\right)^{1/2}= \sigma^{1/4}
\end{align*}
Hence for arbitrary $\epsilon>0$ we have
$$
\Vert f-g_\sigma\Vert_2\leq \epsilon \quad\text{for all } \sigma\leq \epsilon^4
$$
On the other hand,
\begin{align*}
&\sup_{w\in C^\infty(O);\Vert f-w\Vert_2\leq \epsilon}|f(\vartheta)-w(\vartheta)|
 \geq \sup_{\sigma\leq \epsilon^4 }|f(\vartheta)-g_\sigma(\vartheta)|
=\sup_{\sigma\leq \epsilon^4 } \Big(\frac{1}{2\pi \sigma}\Big)^{1/4}=\infty.
\end{align*}

\newpage

\subsection{Functional linear regression and PACE}\label{appendix:FE2}
In the following we discuss the equivalence of the functional linear regression model of \cite{yao2005functional} and the PACE method of \cite{Yao2005} when used to reconstruct functional data from its own irregular and noise contaminated measurements. 

\cite{yao2005functional} consider the function-on-function linear regression model
\begin{equation*}
\E(Y(t)|X)=\mu_Y(t)+\int_{\mathcal{S}}\beta(t,s)(X(s)-\mu_X(s))ds,  
\end{equation*}
where $Y\in\mathbb{L}^2(\mathscr{T})$ denotes the response function, $X\in\mathbb{L}^2(\mathscr{S})$ denotes the predictor function, and $\beta\in\mathbb{L}^2(\mathscr{T}\times\mathscr{S})$ denotes the parameter function. Estimation of $\beta$ is based on a truncated series expansion of $\beta$ which leads to the following approximative model:
\begin{align}
\E(Y(t)|X)
\approx&\mu_Y(t)+\int_{\mathcal{S}}\sum_{k=1}^K\sum_{m=1}^M\frac{\E(\langle\psi_m,X^c\rangle_2\langle\phi_k,Y^c\rangle_2)}{\E(\langle\psi_m,X^c\rangle_2^2)}\psi_m(s)\phi_k(t)X^c(s)ds\notag\\
=&\mu_Y(t)+\sum_{k=1}^K\sum_{m=1}^M\frac{\E(\zeta_m\xi_k)}{\E(\xi_m^2)}\zeta_m\phi_k(t)\notag\\
=&\mu_Y(t)+\sum_{k=1}^K\sum_{m=1}^M\frac{\sigma_{mk}}{\rho_m}\zeta_m\phi_k(t),\label{YMWa}
\end{align}
with $X^c:=X-\mu_X$ and $Y^c:=Y-\mu_Y$, where $\mu_X$ and $\mu_Y$ denote the mean functions of $X$ and $Y$, $(\psi_m)_{1\leq m \leq M}$ and $(\phi_k)_{1\leq k \leq K}$ denote the eigenfunctions associated with the decreasingly ordered eigenvalues of the covariance operators $\E(X\otimes X)$ and $\E(Y\otimes Y)$, $\sigma_{km}$ denotes the covariance of the functional principal component scores $\zeta_m=\langle\psi_m,X^c\rangle_2$ and $\xi_k=\langle\phi_k,Y^c\rangle_2$, and $\rho_m$ denotes the $m$th ordered eigenvalue of the covariance operator $\E(X\otimes X)$. The authors propose to estimate $\E(Y(t)|X)$ by plugging estimates $\hat\mu_Y$, $\hat\sigma_{mk}$, $\hat\rho_m$, $\hat\zeta_m^\ast$, $\hat\phi_k$ to be obtained for sparse functional data of $X$ and $Y$.

In our context of reconstructing partially observed functions, we have $Y=X$, such that $\E(X(t)|X)=X(t)$ as well as $\mu_Y(t)=\mu_X(t)$ for all $t\in\mathscr{S}$. Moreover, $\zeta_k=\xi_k$, with $\E(\zeta_m\zeta_k)=0$ for all $m\neq k$. So, \eqref{YMWa} becomes the truncated Karhunen-Lo\'eve decomposition of $X$
\begin{align}\label{YMWb}
X(t)\approx&\mu_X(t)+\sum_{k=1}^K\zeta_k\phi_k(t).
\end{align}

Let's consider the case of sparse functional data, where the functions $X$ are not fully observed, but only at a few irregular measurements $(U_{l},S_{l})_{1\leq l\leq L}$, with $U_{l}=X(S_{l})+\varepsilon_{l}$. This case prevents the direct computation of the functional principal component scores $\zeta_k$, $k=1,\dots,K$. Therefore, \cite{Yao2005} propose to predict the scores $\zeta_k$ using the conditional expectations, $\tilde\zeta_k$, of $\zeta_k$ given the irregular measurements $(U_{l},S_{l})_{1\leq l\leq L}$, of $X$ \citep[see Equation (4) in][]{Yao2005}.

The empirical version of \eqref{YMWb} for the case of sparse functional data is given by
\begin{align}\label{YMW2}
\widehat{X}_{K}(t)=&\hat\mu_X(t)+\sum_{k=1}^K\hat\zeta^\ast_k\hat\phi_k(t),
\end{align}
where $\hat{\mu}_X$ and $\hat\phi_k$ denote the consistent estimators of the meanfunction and the $k$th eigenfunction as described in \cite{Yao2005}, and where $\hat{\zeta}^\ast_k$ denotes the estimator of the conditional expectation, $\tilde\zeta_k^\ast$, of $\zeta_k$, given the irregular measurements $(U_{l},S_{l})_{1\leq l\leq L}$ of $X$ \citep[see Equation (5) in][]{Yao2005}.

Equation \eqref{YMW2} is just the PACE method as proposed in \cite{Yao2005} for predicting the trajectory of $X$ given its own irregular measurements \citep[see Equation (6) in][]{Yao2005}. Both articles, \cite{Yao2005} and \cite{yao2005functional}, use the same nonparametric estimators for the equivalent model components and predict the principal component scores using conditional means as original proposed by \cite{Yao2005}. So, the functional linear regression model of \cite{yao2005functional} is equivalent to the PACE method of \cite{Yao2005} when used to reconstruct functional data from its own irregular and noise contaminated measurements.

\newpage

\section{Visualizations of simulation results}\label{visual}

\begin{figure}[!ht]
\centering
\includegraphics[width=1\textwidth]{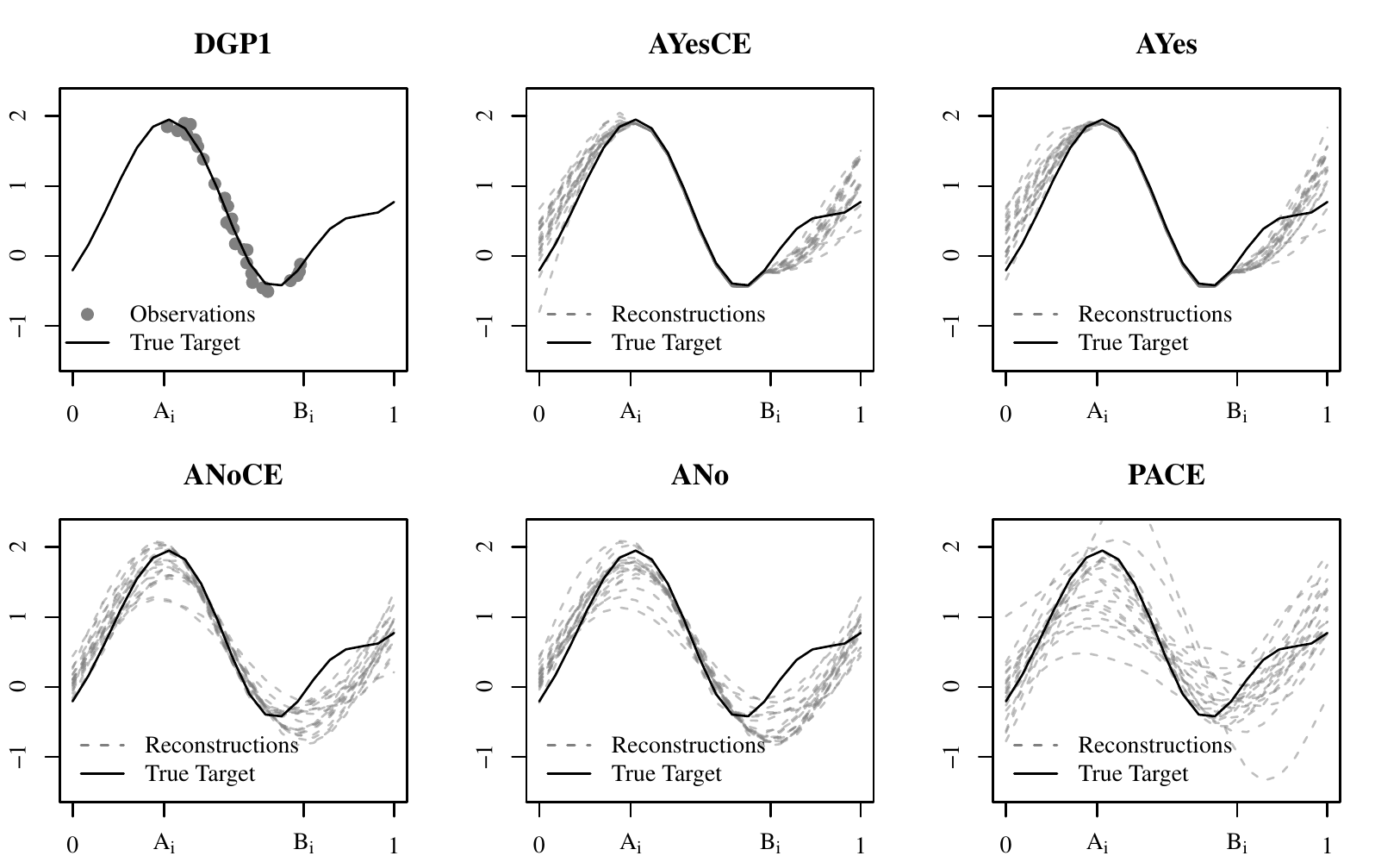}
\caption[]{Reconstruction results for DGP1 with $n=100$ and $m=30$.}
\label{fig:simresults_DGP1}
\end{figure}
\begin{figure}[!hb]
\centering
\includegraphics[width=1\textwidth]{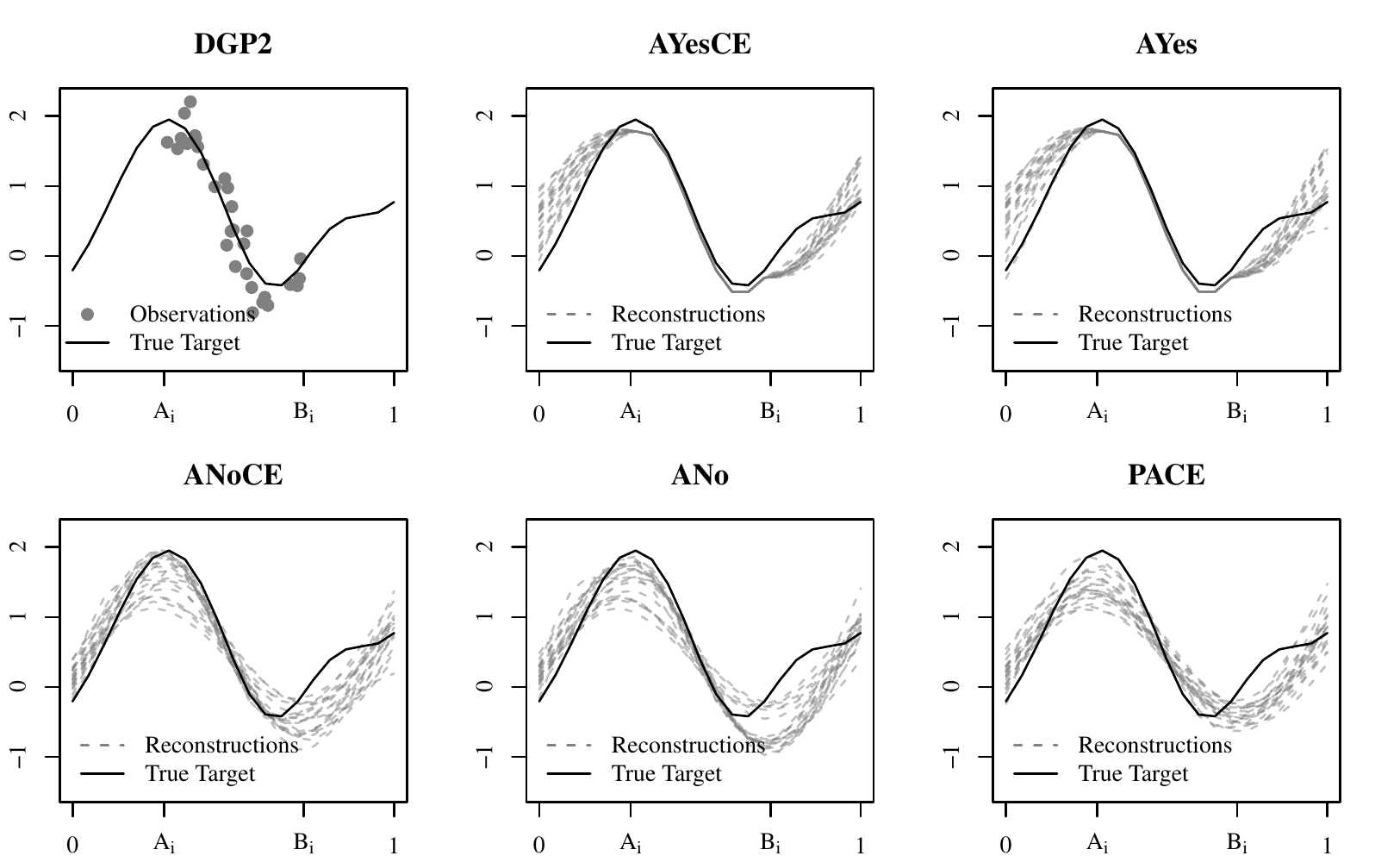}
\caption[]{Reconstruction results for DGP2 with $n=100$ and $m=30$.}
\label{fig:simresults_DGP2}
\end{figure}
\begin{figure}[!ht]
\centering
\includegraphics[width=1\textwidth]{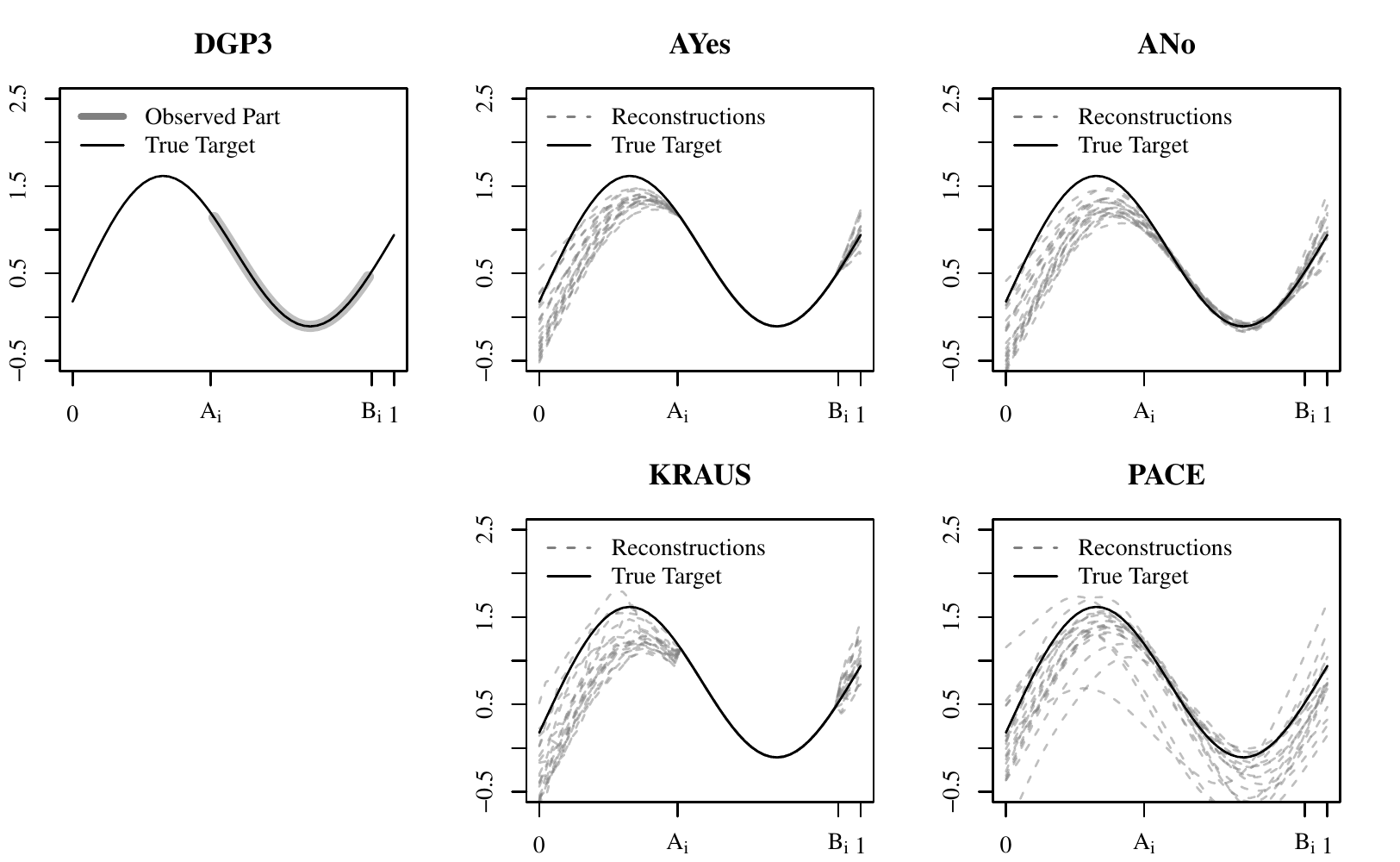}
\caption[]{Reconstruction results for DGP3 with $n=100$.}
\label{fig:simresults_DGP3}
\end{figure}
\begin{figure}[!hb]
\centering
\includegraphics[width=1\textwidth]{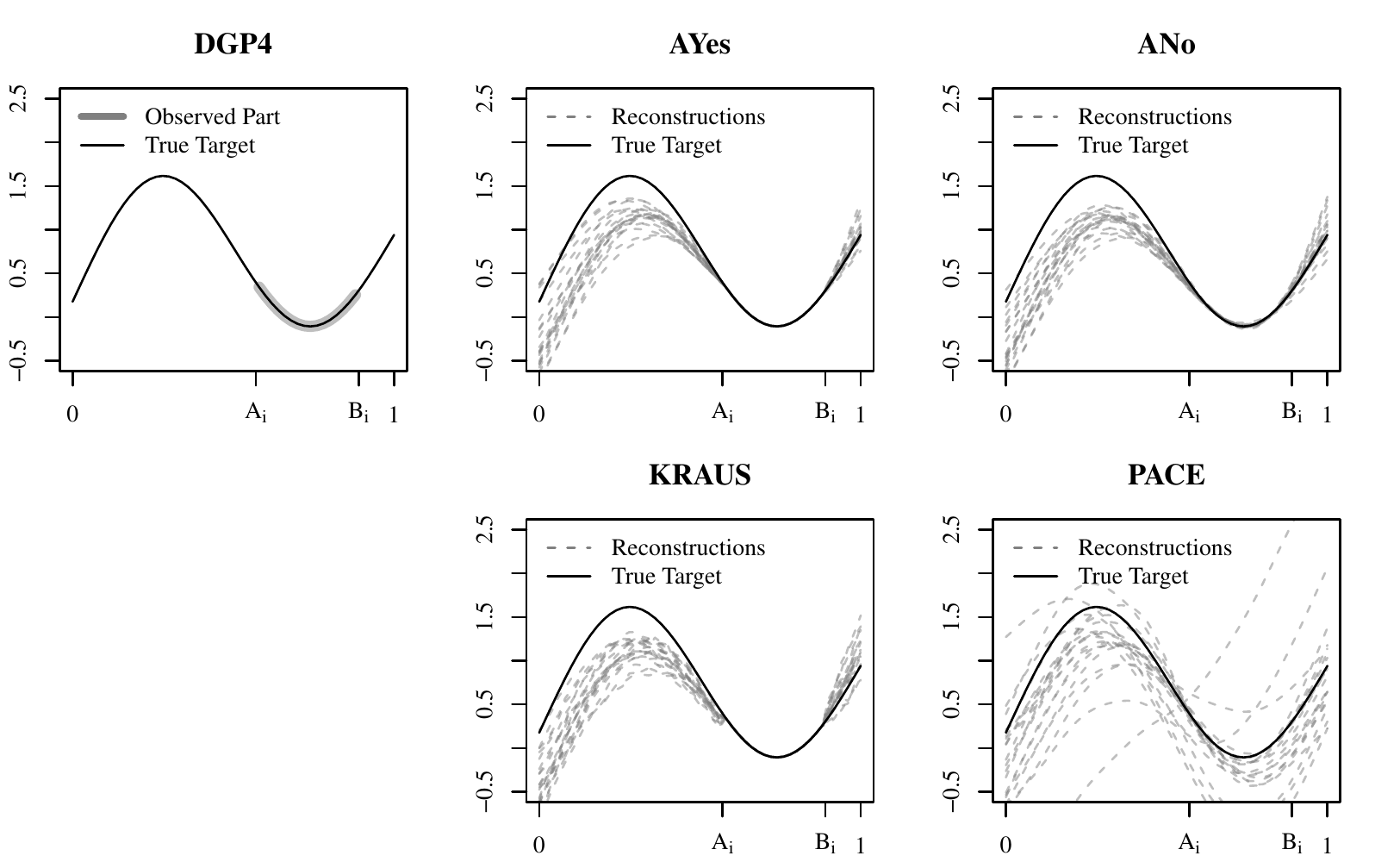}
\caption[]{Reconstruction results for DGP4 with $n=100$.}
\label{fig:simresults_DGP4}
\end{figure}

\newpage

\bibliographystyleappendix{imsart-nameyear}
\bibliographyappendix{bibfile}

\end{document}